\documentclass[1 [leqno,11pt]{amsart}
\usepackage{amssymb, amsmath}
\newcommand{\andf}{\quad\hbox{and}\quad}
\newcommand{\with}{\quad\hbox{with}\quad}
\def\Supp{\mathop{\rm Supp}\nolimits\ }
\newcommand{\newcom}{\newcommand}

\def\inte#1{
\displaystyle\mathop{#1\kern0pt}^\circ }

\def\fa{\frak{a}}
\def\fb{\frak{b}}
\newcom{\al}{\alpha}
\newcom{\de}{\delta}
\newcom{{\Th}}{\Phi}
\newcom{\be}{\beta}
\newcom{\s}{\sigma}
\newcom{\eps}{\epsilon}
\newcom{\ve}{\varepsilon}
\newcom{\ga}{\gamma}
\newcom{\Ga}{\Gamma}
\newcom{\ka}{\kappa}
\newcom{\Lam}{\Lambda}
\newcom{\lam}{\lambda}
\newcom{\vp}{\varphi}
\newcom{\om}{\omega}
\newcom{\Sig}{\Sigma}
\newcom{\sig}{\sigma}
\newcom{\tht}{\theta}
\newcom{\tri}{\triangle}
\newcom{\oo}{\infty}
\newcom{\h}{{\rm h}}
\newcom{\vphi}{\varphi}
\newcom{\cB}{{\mathcal B}}
\newcom{\cC}{{\mathcal C}}
\newcom{\cD}{{\mathcal D}}
\newcom{\cF}{{\mathcal F}}
\newcom{\cL}{{\mathcal L}}
\newcom{\cM}{{\mathcal M}}
\newcom{\cP}{{\mathcal P}}
\newcom{\cS}{{\mathcal S}}
\newcom{\cQ}{{\mathcal Q}}
\newcom{\cT}{{\mathcal T}}
\newcom{\cY}{{\mathcal Y}}
\newcom{\cZ}{{\mathcal Z}}
\newcom{\R}{\Bbb R}
\newcom{\T}{\Bbb T}
\newcom{\N}{\Bbb N}
\newcom{\Z}{\Bbb Z}
\newcom{\C}{\Bbb C}
\newcom{\E}{\Bbb E}
\let\wh=\widehat
\def\dive{\mathop{\rm div}\nolimits}
\let\e=\varepsilon

\def\ve{v^\varepsilon}

\def\vp{v_{\Psi}}
\def\fk{\frak{K}}
\def\ekt{e^{\frak{K}t}}
\def\ektp{e^{\frak{K}t'}}

\def\ueps{u^\epsilon}
\def\veps{v^\epsilon}
\def\fw{\frak{W}}
\def\fs{\frak s}

\newcom{\f}{\frac}
\newcom{\dint}{\displaystyle\int}
\newcom{\dsum}{\displaystyle\sum}
\newcom{\dlim}{\displaystyle\lim}
\newcom{\ov}{\overline}
\newcom{\wt}{\widetilde}
\newcom{\pa}{\partial}
\newcom{\p}{\partial}
\newcom{\cK}{\mathcal{K}}
\newcom\na{\nabla}
\newcom{\D}{\Delta}
\newcom\rto{\rightarrow}
\newcom\lto{\leftarrow}
\newcom\mto{\mapsto}
\newcom{\disp}{\displaystyle}
\newcom{\non}{\nonumber}
\newcom{\no}{\noindent}
\newcom{\QED}{$\square$}

\def\eqdefa{\buildrel\hbox{\footnotesize def}\over =}

\newcommand{\beq}{\begin{equation}}
\newcommand{\eeq}{\end{equation}}
\newcommand{\ben}{\begin{eqnarray}}
\newcommand{\een}{\end{eqnarray}}
\newcommand{\beno}{\begin{eqnarray*}}
\newcommand{\eeno}{\end{eqnarray*}}

\newtheorem{Def}{Definition}[section]
\newtheorem{thm}{Theorem}[section]
\newtheorem{lem}{Lemma}[section]
\newtheorem{rmk}{Remark}[section]

\renewcommand{\theequation}{\thesection.\arabic{equation}}

\begin{document}
\title[hydrostatic approximation of  hyperbolic Navier-Stokes  system]
{Global hydrostatic approximation of  hyperbolic Navier-Stokes system with small Gevrey class 2 data}

\author{Marius Paicu}
\address{Universit\'e  Bordeaux \\
 Institut de Math\'ematiques de Bordeaux\\
F-33405 Talence Cedex, France} \email{marius.paicu@math.u-bordeaux.fr}

\author{Ping Zhang}%
\address{\footnote{Corresponding author} Academy of
Mathematics $\&$ Systems Science and  Hua Loo-Keng Key Laboratory of
Mathematics, The Chinese Academy of Sciences, Beijing 100190, China, and School of Mathematical Sciences,
University of Chinese Academy of Sciences, Beijing 100049, China.
} \email{zp@amss.ac.cn}

\date{\today}

\begin{abstract}
We investigate the hydrostatic approximation of  a hyperbolic version of  Navier-Stokes equations, which is  obtained by using  Cattaneo type law instead of Fourier law, evolving  in a thin strip $\R\times (0,\varepsilon)$. The formal limit of these equations is a hyperbolic Prandtl  type equation. We first prove the global existence  of  solutions to these equations under a uniform smallness assumption on the data in Gevrey $2$ class. Then we justify the limit globally-in-time from the anisotropic hyperbolic Navier-Stokes system to the hyperbolic Prandtl system
with such Gevrey $2$ class data. Compared with \cite{PZZ2} for the hydrostatic approximation of 2-D classical Navier-Stokes system with analytic data,
here the initial data belong to the Gevrey $2$ class, which is very sophisticated even for the well-posedness  of the classical
Prandtl system (see \cite{DG19,WWZ1}), furthermore, the estimate of the pressure term in the hyperbolic Prandtl system arises additional
difficulties.
\end{abstract}

\maketitle

\vskip-0.3cm

\centerline{\it{ In  memory of Professor Genevi\`eve Raugel}}

\vskip0.5cm


\qquad\qquad hyperbolic Prandtl system.


\renewcommand{\theequation}{\thesection.\arabic{equation}}
\setcounter{equation}{0}


\section{Introduction}\label{sect1}

The Navier-Stokes system governing the evolution of Newtonian incompressible viscous
fluids reads
\begin{equation*}
	(NS)\qquad
	\left\{\;
	\begin{aligned}
		& \partial_t U + U\cdot \nabla U  -\nu \Delta U + \nabla P= 0, \\
		&\text{ div } U = 0,
	\end{aligned}
	\right.
\end{equation*}
which is obtained from  the momentum equations
 \begin{equation}\label{S1mom}
\partial_t U+U\cdot\nabla U=\dive\tau,
\end{equation}
coupled with the Fourier law  to describe the stress tensor
\begin{equation}\label{S1str}
\tau(t) = -P {\rm Id} + \nu \left(\nabla U + (\nabla U)^\top\right)(t),
\end{equation}
where $U$ stands for the velocity field, $P$  the scalar pressure function and $\nu$ is the viscosity coefficient of the fluid.

 Nevertheless, the system $(NS)$ gives rise to a physical paradox coming from the fact that the equation of the velocity  has infinite propagation speed. In order to avoid this non-physical aspect,  Cattaneo \cite{Cattaneo} (see also Vernotte \cite{Vernotte}) proposed  to modify
  the heat equation to a hyperbolic version known as Cattaneo's heat transfer law. More precisely, they proposed to replace the Fourier law \eqref{S1str} by the following hyperbolic model
  $$ \frac{1}{c^2} \pa_t^2 \theta + \frac{1}{\beta}\pa_t \theta - \Delta \theta = 0. $$
 The resulting system has  finite propagation speed and constitutes a satisfactory physical model  being  compatible with the principle of relativity and with the second law of thermodynamics. It is therefore natural to consider a hyperbolic Navier-Stokes system by adding the term $ \eta\p_t^2u$ to the classical Navier-Stokes system $(NS),$ where $\eta= \frac{1}{c^2}$ is a small parameter. The hyperbolic Navier-Stokes system \eqref{HNS} has been extensively  studied in the literature and one may check \cite{BA2019, BNP,  CHR,  IHachicha,  RS2000,  MR2007, MR2008} and the references therein. The justification of this system by using the Cattaneo's law follows by considering that the stress tensor is given by the solution of the retarded equation
\begin{equation}
 \tau(t+\eta) = -P {\rm Id} + \nu \left(\nabla U + (\nabla U)^\top\right)(t).
\end{equation}
By using the fact that $\tau(t+\eta)\simeq \tau(t) + \eta\pa_t \tau(t)$ and applying the operator $ (\eta \pa_t + Id)$ to the momentum equations
\eqref{S1mom},
 we obtain the following quasilinear hyperbolic version of Navier-Stokes equations:
\begin{equation}
	\label{HPNS}
	\left\{\;
	\begin{aligned}
		&\eta\pa_t^2 U + \partial_t U + U\cdot \nabla U+\eta(U\cdot \nabla U)_t  -\nu \Delta U + \nabla P= 0, \\
		&\text{ div } U = 0,  \\
		&U|_{t=0} = U_0, \text{ \ \ \ } \pa_t U|_{t=0} =U_1.
	\end{aligned}
	\right.
\end{equation}
For simplicity,  we  neglect the term $\eta(U\cdot\nabla U)_t$ in \eqref{HPNS}, and then we get the hyperbolic version of  Navier-Stokes system:
\begin{equation}
	\label{HNS}
	\left\{\;
	\begin{aligned}
		&\eta\pa_t^2 U+ \partial_t U + U\cdot \nabla U  -\nu \Delta U + \nabla P= 0, \\
		&\text{ div }  U= 0,\\
		&U|_{t=0} = U_0, \text{ \ \ \ } \pa_t U|_{t=0} =U_1.
	\end{aligned}
	\right.
\end{equation}

We remark that the system \eqref{HNS} can also  be obtained by the relaxation of the Euler's equations under a diffusive scaling, as it is proposed in \cite{BNP}.  Indeed, Brenier, Natalini and Puel \cite{BNP} proved the  global existence and uniqueness of solutions to the following system:
 \begin{equation}
	\label{relaxed-Euler}
	\left\{\;
	\begin{aligned}
		&\pa_t U+ \nabla\cdot V =\nabla P, \\
		&\sqrt{\eta}\partial_t V +\frac{\nu}{\sqrt{\eta}}\nabla U
=\frac{1}{\sqrt{\eta}}(U \otimes U -V ), \\
		&\dive U= 0, \\
		&U|_{t=0} = U_0, \text{ \ \ \ } V|_{t=0} =V_0,
	\end{aligned}
	\right.
\end{equation}	
 with initial data in ${H}^2(\Bbb{T}^2)^2\times {H}^1(\Bbb{T}^2)^4 $, where $\Bbb{T}^2$ is the periodic  square $\mathbb{R}^2/\mathbb{Z}^2$.
Moreover they proved the convergence of such solutions to a  smooth solution of the classical Navier-Stokes system. In fact, one can show
that as $\eta\to 0+,$ the equations governing the leading order terms in \eqref{relaxed-Euler} are again \eqref{HNS}.

Paicu and Raugel in \cite{MR2007,MR2008} obtained the global existence and uniqueness result  for the
 system \eqref{HNS} with significantly improved regularity for the initial data by using the Strichartz inequalities for  this dispersive model.
  Hachicha \cite{IHachicha} obtained the global  existence and uniqueness result of the perturbed Navier-Stokes system \eqref{HNS} under suitable smallness assumptions on the initial data in the space $ {H}^{\frac{n}{2} + \delta} \times {H}^{\frac{n}{2}-1 + \delta}(\mathbb{R}^n)$ for $n=2,3$. Moreover,  she proved the local in time convergence  of such solutions  towards solution of the classical Navier-Stokes system  with initial data in $ {H}^{\frac{n}{2}-1 + s}(\mathbb{R}^n)$ for $s>0$.  Very recently, Coulaud, Hachicha and  Raugel  \cite{CHR} investigated the well-posedness of
   a hyperbolic quasi-linear version of the Navier-Stokes equation in $\mathbb{R}^2.$
   In particular,
   under smallness assumptions on the data,  they proved the global well-posedness of the system.

 On the other hand, we recall that in order to describe geophysical flows on the Earth, it is usually assumed that vertical scale is much smaller than horizontal one as the fluid layer depth is small compared to the radius of the earth. In order to to take into account this anisotropy between horizontal and vertical directions, we shall suppose that the fluid is evolving in a thin striped domain  and for a vanishing viscosity so that
  they are good approximation of the atmospheric flow and oceanic flow.

The purpose of this paper is first to show the global existence solutions to \eqref{HNS} in the thin domain $\mathcal{S}_\varepsilon=\R\times(0,\varepsilon)$ and for a vanishing viscosity $\nu=\varepsilon^2,$
 we denote the corresponding solutions as $(U^\e, P^\e),$
and then to justify  the limit system when the parameter $\varepsilon$ goes to zero.  For simplicity, here we take $\eta=1$ in
\eqref{HNS}. The case when $\eta\to 0+$ also presents  interest and will be considered separately.

We complement the system \eqref{HNS} with the non-slip boundary condition
\beno
U|_{y=0}=U|_{y=\e}=0,
\eeno
and the initial condition
\beno
U|_{t=0}=\left(u_0\bigl(x,\f{y}\e\bigr), \e v_0\bigl(x,\f{y}\e\bigr) \right)=U_0^\e \quad  \mbox{in} \ \ \cS^\e.
\eeno

As in \cite{BL, Lions96, PZZ2} for the classical Navier-Stokes system, we  write
\beq \label{S1eq1.7}
U(t,x,y)=\Bigl(u^\e\bigl(t,x,\f{y}\e\bigr), \e v^\e\bigl(t,x,\f{y}\e\bigr) \Bigr) \andf P(t,x,y)=p^\e\bigl(t,x,\f{y}\e\bigr).
\eeq
Let $\cS\eqdefa \bigl\{(x,y)\in \R^2:\ 0<y<1\bigr\}$. After a natural change of scale (see \cite{PZZ2}),  the system \eqref{HNS} becomes the following scaled anisotropic hyperbolic Navier-Stokes system:
\begin{equation}\label{S1eq1.8}
 \quad\left\{\begin{array}{l}
\displaystyle \p^2_tu^\e+ \p_tu^\e+u^\e\p_x u^\e+v^\e\p_yu^\e-\e^2\p_x^2u^\e-\p_y^2u^\e+\p_xp^\e=0\ \ \mbox{in} \ \cS\times ]0,\infty[,\\
\displaystyle \e^2\left(\p_t^2v^\e+\p_tv^\e+u^\e\p_x v^\e+v^\e\p_yv^\e-\e^2\p_x^2v^\e-\p_y^2v^\e\right)+\p_yp^\e=0,\\
\displaystyle \p_xu^\e+\p_yv^\e=0,\\
\displaystyle \left(u^\e, v^\e\right)|_{t=0}=\left(u_0, v_0\right)\andf \left( \p_tu^\e,  \p_tv^\e\right)|_{t=0}=\left(u_1, v_1\right),
\end{array}\right.
\end{equation}
together with the boundary condition
\beq \label{S1eq5}
\left(u^\e, v^\e\right)|_{y=0}=\left(u^\e, v^\e\right)|_{y=1}=0.
\eeq

Formally taking $\e\to 0$ in the system \eqref{S1eq1.8}, we obtain the hyperbolic Prandtl equations:
\begin{equation}\label{S1eq1}
 \quad\left\{\begin{array}{l}
\displaystyle \p_t^2u+\p_tu+u\p_x u+v\p_yu-\p_y^2u+\p_xp=0\ \ \mbox{in} \ \cS\times ]0,\infty[,\\
\displaystyle\p_yp=0,\\
\displaystyle \p_xu+\p_yv=0,\\
\displaystyle u|_{t=0}=u_0,\quad u_t|_{t=0}=u_1,
\end{array}\right.
\end{equation}
together with the boundary condition
\beq \label{S1eq2}
\left(u, v\right)|_{y=0}=\left(u, v\right)|_{y=1}=0.
\eeq

The goal of this paper is to justify the limit from the system \eqref{S1eq1.8} to the system \eqref{S1eq1} with initial data
in the Gevrey class 2.
Before we present the main results of this paper, we remark that similar to the classical Prandtl equation $(P),$ the nonlinear term $v\p_yu$ in \eqref{S1eq1} will lead to one derivative loss in the $x$ variable in the process of energy estimates. Thus, it is natural to work with analytic data in order to overcome this difficulty if we don't impose extra structural assumptions on the initial data (see \cite{GD, Renard09}). Indeed,  for the data which is analytic in $x,y$ variables,
Sammartino and Caflisch \cite{Caf} established the local
well-posedness result of the system $(P)$ in the upper half space. Later, the analyticity in $y$
variable was removed by Lombardo, Cannone and Sammartino in
\cite{Can}. Lately, Paicu and Zhang \cite{PZ5} proved the global well-posedness of the system $(P)$ with small analytic data.
 We  mention that the optimal local
 well-posedness result of Prandtl system is obtained by Dietert and G\'erard-Varet in \cite{DG19} (see also \cite{GM}), where the authors proved
 the local well-posedness of $(P)$ with Gevrey class 2 data (one may check \cite{WWZ1} for the corresponding global result).
 We also mention that
 for a class of convex data, G\'erard-Varet, Masmoudi and Vicol \cite{GMV} proved the well-posedness of the  hydrostatic Navier-Stokes system in the Gevrey class.

Our first result in this paper is concerned with the global well-posedness of the hyperbolic Prandtl system \eqref{S1eq1} with
small data in Gevrey class 2 for $x$ variable.

\begin{thm}\label{th1.2}
{\sl \begin{enumerate}
\item[(1)] Let $a>0.$ We assume that the initial data satisfies
\beq\label{S1eq3}
\begin{split} \frak{H}^0_{\f12}(u_0,u_1)\leq  c_0 \with 
\frak{H}^0_s(u_0,u_1)&\eqdefa 
\bigl\|e^{a|D_x|^{\f12}}(u_0,u_1,\p_yu_0)\bigr\|_{\cB^{s}} \\ &+\sqrt{a\fk}\bigl\|e^{a|D_x|^{\f12}}u_0\bigr\|_{\cB^{s+\f14}}+a\fk\bigl\|e^{a|D_x|^{\f12}}u_0\bigr\|_{\cB^{s+\f12}},
\end{split}
\eeq
for some $c_0$ sufficiently small and there holds the compatibility condition $\int_0^1u_0dy=\int_0^1u_1dy=0$. Then the system (\ref{S1eq1}) has a  unique global solution $u$ so that for any $t>0,$ there holds
\beq\label{S4eq10}
\begin{split}
{\rm E}_{\f12}(u)(t)\leq  C\frak{H}^0_{\f12}(u_0,u_1),
\end{split}
\eeq  where
\beq\label{S1eq10}
\begin{split}
{\rm E}_{s}(u)(t)\eqdefa &\bigl\|\ektp \bigl(u_\Phi, \p_yu_\Phi, (\p_tu)_\Phi\bigr)\bigr\|_{\wt{L}^\infty_t(\cB^{s})}
+\sqrt{a\fk}\|e^{\f34\fk t'} u_\Phi\|_{\wt{L}^\infty_{t}(\cB^{s+\f14})}\\
&+a\fk\|e^{\f12\fk t'} u_\Phi\|_{\wt{L}^\infty_{t}(\cB^{s+\f12})}
+\bigl\|\ektp\bigl((\p_t u)_\Phi,\p_y u_\Phi\bigr)\bigr\|_{\wt{L}^2_t(\cB^{s})}.
\end{split}
\eeq
$u_\Phi$ will be determined by \eqref{S4eq1}, the constant $\frak{K}\eqdefa \min\left(\f1{6},\f1{4(1+\cK)}\right)$ with
$\cK$ being determined by Poincar\'e inequality on the strip $\cS$ (see \eqref{S4eq5f}), and the functional spaces will be presented in Appendix \ref{sect2}.

\item[(2)] If we assume moreover that for $s>-\f14,$ $\frak{H}^0_{s}(u_0,u_1)<\infty,$ then one has
\beq\label{S1eq12}  {\rm E}_{s}(u)(t)\leq C\frak{H}^0_{s}(u_0,u_1)\quad \forall\ t\in\R^+.
 \eeq

\item[(3)] If we assume in addition that
\beq \label{S1eq13}
\begin{split}
\frak{H}^0_{\f12}(u_0,u_1)&+\frak{H}^0_{\f34}(u_0,u_1)+\frak{H}^0_{\f12}(\p_yu_0,\p_yu_1)\\
&+\bigl\|e^{a|D_x|^{\f12}}u_0\bigr\|_{\cB^{\f12}}
\bigl\|e^{a|D_x|^{\f12}}\p_yu_0\bigr\|_{\cB^{\f12}}\bigl\|e^{a|D_x|^{\f12}}u_0\bigr\|_{\cB^{\f32}}^2 \leq c_1,
\end{split}
\eeq
for some $c_1$ sufficiently small, then for any $\fs\geq\f12,$ one has
\beq \label{S1eq14}
\begin{split}
  {\rm E}_{\fs}(\p_y u)(t)&+\|\ektp(\p_t^2u)_\Phi\|_{\wt{L}_t^2(\cB^{\fs})}\leq C\Bigl(\frak{H}^0_{\fs}(u_0,u_1)+\frak{H}^0_{\fs+\f14}(u_0,u_1)\\
  &+\frak{H}^0_{\fs}(\p_yu_0,\p_yu_1)+\bigl\|e^{a|D_x|^{\f12}}u_0\bigr\|_{\cB^{\fs}}
\bigl\|e^{a|D_x|^{\f12}}\p_yu_0\bigr\|_{\cB^{\fs}}\bigl\|e^{a|D_x|^{\f12}}u_0\bigr\|_{\cB^{\fs+1}}^2\Bigr).
\end{split}
  \eeq
\end{enumerate}
}
\end{thm}

\begin{rmk} Compared with the wellposedness result of the hydrostatic Navier-Stokes system with initial data
in Gevrey class in \cite{GMV}, here we do not need the convex assumption on  the initial data.
Compared with the proof of the well-posedness result of the Prandtl system in \cite{DG19, WWZ1} with initial data in the Gevrey class 2,
our proof here is much simplified. The main idea is to use multi time-weighted norms like $\|u\|_{\wt{L}^2_{t,\dot\tht(t)}(\cB^{s+\f14})},$
$\|u\|_{\wt{L}^2_{t,\dot\tht^2(t)}(\cB^{s+\f12})}$ and $\|u\|_{\wt{L}^2_{t,\dot\tht^3(t)}(\cB^{s+\f34})}$ (see \eqref{Esu} and Definition \ref{def1.1} below).
\end{rmk}

The second result  is about the global well-posedness of the system \eqref{S1eq1.8} with small Gevrey class 2 data in $x$ variable. The main interesting point is that the smallness of data is independent of $\e$ and there holds the global uniform estimate \eqref{S3eq18} with respect to the parameter $\e$.

\begin{thm}\label{thm1.1}
{\sl Let $a>0.$ We assume that the initial data satisfies
\beq\label{S1eq8}
\begin{split}
\bigl\|&e^{a|D_x|^{\f12}}(u_0,\e v_0, u_1,\e v_1,\p_y(u_0,\e v_0))\bigr\|_{\cB^{\f12}}+\bigl\|e^{a|D_x|^{\f12}}(u_0,\e v_0)\bigr\|_{\cB^{\f34}\cap\cB^{1}}\leq c_2
\end{split}
\eeq
for some $c_2$ sufficiently small. Then the system \eqref{S1eq1.8} has a unique global solution $(u^\e,v^\e)$ so that for any $t>0,$ there holds
\beq\label{S1eq16}
{\rm E}^1_{\fk}(u^\e,v^\e)(t)\leq C\frak{H}^1_\e(u_0,v_0),
\eeq
where
\beq\label{S3eq18}
\begin{split}
{\rm E}^1_{\fk}(u^\e,v^\e)(t)\eqdefa &\bigl\|\ektp \bigl((u^\e,\e v^\e)_\Phi,\e\p_x(u^\e,\e v^\e)_\Phi, \p_y(u^\e,\e v^\e)_\Phi, (\p_tu^\e,\e\p_tv^\e)_\Phi\bigr)\bigr\|_{\wt{L}^\infty_t(\cB^{\f12})}\\
&+\sqrt{a\fk}\|e^{\f34\fk t'} (u^\e,\e v^\e)_\Phi\|_{\wt{L}^\infty_{t}(\cB^{\f34})}+a\fk\|e^{\f12\fk t'} (u^\e,\e v^\e)_\Phi\|_{\wt{L}^\infty_{t}(\cB^{1})}\\
&+\bigl\|\ektp\bigl(\e\p_x (u^\e,\e v^\e),\p_y (u^\e,\e v^\e), \p_t(u^\e, \e v^\e) \bigr)_\Phi\|_{\wt{L}^2_t(\cB^{\f12})},
\end{split}
\eeq and
\beq\label{S1eq17}
\begin{split}
  \frak{H}^1_\e(u_0,v_0)\eqdefa &\bigl\|e^{a|D_x|^{\f12}}(u_0,\e v_0, \e\p_x(u_0,\e v_0), \p_y(u_0,\e v_0), u_1,\e v_1)\bigr\|_{\cB^{\f12}}\\
&+\sqrt{a\fk}\bigl\|e^{a|D_x|^{\f12}}(u_0,\e v_0)\bigr\|_{\cB^{\f34}}+a\fk\bigl\|e^{a|D_x|^{\f12}}(u_0,\e v_0)\bigr\|_{\cB^{1}}.
\end{split}
\eeq
where $(u_\Phi^\e,v_\Phi^\e)$ will be given by \eqref{S4eq1}.  }
\end{thm}

The third result is concerned with the convergence of the solutions from the scaled anisotropic  hyperbolic Navier-Stokes system \eqref{S1eq1.8}  to the hyperbolic Prandtl system \eqref{S1eq1}, which corresponds to \cite{PZZ2} for hydrostatic approximation of Navier-Stokes system in a thin strip with
analytic data in the tangential variable.
Compared with \cite{PZZ2}, here our initial data belong to Gevrey class 2.

\begin{thm}\label{thm3}
{\sl Let $a>0$ and $(u_0, v_0)$ satisfy \eqref{S1eq8}. Let $(u_0,u_1)$ verify \eqref{S1eq3},
 \eqref{S1eq13} and $\frak{H}_3^0(u_0,u_1)<\infty,$ and the compatibility condition
$\int_0^1u_0\,dy=\int_0^1u_1\,dy=0.$ Let $(u^\e,v^\e)$ and $u$ be the global solutions of the system \eqref{S1eq1.8} and \eqref{S1eq1}
obtained respectively in Theorems \ref{th1.2} and \ref{thm1.1}. Then we have
\beq\label{S1eq15}
\begin{split}
{\rm E}_0^1(w^{1}_{\e},w_{\e}^2)(t)
&\leq C\bigl(\frak{H}^1_\e(w^1_{0,\e},w^2_{0,\e})+
M\e\bigr),
\end{split}
\eeq
Here the energy functionals ${\rm E}^1_0(w^{1}_{\e},w_{\e}^2)(t)$ and $\frak{H}^1_\e(w^1_{0,\e},w^2_{0,\e})$ are determined
 respectively by \eqref{S3eq18} and \eqref{S1eq17}, $w^1_\e\eqdefa u^\e-u,\ w^2_\e\eqdefa v^\e-v$ and $v$ is determined from $u$ via $\p_xu+\p_yv=0$ and $v|_{y=0}=v|_{y=1}=0.$
}
\end{thm}

\begin{rmk}
 Similar convergence result from hyperbolic  Navier-Stokes system to hyperbolic hydrostatic equations has been rigorously justified in the
analytic functional framework by  Aarach \cite{NA}.
\end{rmk}

\medskip

The structure of this paper lists as follows: in Section \ref{Sect2}, we prove the global well-posedness of the system \eqref{S1eq1}
with initial data in the Gevrey class 2 for the tangential variable. Let us point out that the multi time-weighted Chemin-Lerner nroms  plays an important role. These ideas are common to the proof of the three theorems in this paper.

Section \ref{Sect3} is devoted to the proof of the propagation of regularity for $\p_y u_\Phi,$ which arises additional
difficulty in the estimate of the pressure function compared with that in \cite{PZZ2}.

Section \ref{Sect4} is devoted to the proof of the global well-posedness of the  system \eqref{S1eq1.8}, namely the proof of
Theorem \ref{thm1.1}

Finally we present the proof of Theorem \ref{thm3} in Section \ref{Sect5}.

In the appendix \ref{sect2}, we shall collect some basic tools on the functional framework used in this paper.

\medbreak
We end this introduction by the  notations that we shall use in this context. For~$a\lesssim b$, we mean that there is a
uniform constant $C,$ which may be different on different lines,
such that $a\leq Cb$.  We denote by $(a|b)_{L^2}$ the $L^2(\cS)$
inner product of $a$ and $b$.
We designate
by $L^p_T(L^q_{\rm h}(L^r_{\rm v}))$ the space $L^p(]0,T[;
L^q(\R_{x};L^r(\R_y))).$ Finally, we denote by $(d_k)_{k\in\Z}$ (resp. $(d_k(t))_{k\in\Z}$) to be a generic
element of $\ell^1(\Z)$ so that $\sum_{k\in\Z}d_k=1$ (resp. $\sum_{k\in\Z}d_k(t)=1$).

\renewcommand{\theequation}{\thesection.\arabic{equation}}
\setcounter{equation}{0}
\section{Global well-posedness of the system \eqref{S1eq1}}\label{Sect2}

In this section, we study the global well-posedness of the hyperbolic Prandtl equations \eqref{S1eq1} with small Gevrey class 2 data£¬
namely, we are going to present the proof of part (1) and part (2) of Theorem \ref{th1.2}.

As in \cite{Ch04, CGP, PZ5, PZZ2,  mz1, mz2, ZZ}, especially motivated by \cite{WWZ1},
we define
\beq \label{S4eq1}
u_\Phi(t,x,y)\eqdefa \cF_{\xi\to
x}^{-1}\bigl(e^{\Phi(t,\xi)}\widehat{u}(t,\xi,y)\bigr) \with \Phi(t,\xi)\eqdefa (a-\lam \tht(t))|\xi|^{\f12},
\eeq
where the quantity $\tht(t)$ describes the evolution of the loss of
Gevrey radius  of $u,$ which is determined by
\begin{equation}\label{S4eq-1}
 \dot{\tht}(t)=\de^{\f12} e^{-\f{\frak{K}}2t}\with \tht|_{t=0}=0 \with  \de\eqdefa \left(\f{a\fk}{4\lam}\right)^2,
\end{equation}
where $\lam$ is a large enough constant to be determined later on.  Indeed we observe from \eqref{S4eq-1} that
\beno
\tht(t)=\f{2\de^{\f12}}{\fk}\left(1-e^{-\f{\frak{K}}2t}\right),
\eeno
so that for all $t\in\R^+,$
\beno
a-\lam \tht(t)=\f{a}2\left(1+e^{-\f{\frak{K}}2t}\right)
>\f{a}2. \eeno Hence by virtue of \eqref{S4eq1}, $\Phi(t,\xi)$ verifies
 the
following convex inequality \beq\label{eq4.4} \Phi(t,\xi)\leq
\Phi(t,\xi-\eta)+\Phi(t,\eta)\quad\mbox{for}\quad \forall\ t\in\R^+ \andf
\xi,\eta\in \R. \eeq

\begin{proof}[Proof of  Theorem \ref{th1.2}] {\bf Part (1)}. For simplicity, we just present the {\it a priori} estimates
for smooth enough solutions of \eqref{S1eq1}.
Indeed in  view of \eqref{S1eq1} and \eqref{S4eq1}, we observe
that $u_\Phi$ verifies
\beq\label{S4eq2}
\begin{split}
&\p_t(\p_tu)_\Phi+\lam\dot{\tht}(t)|D_x|^{\f12}(\p_tu)_\Phi+(\p_tu)_\Phi+(u\pa_x u)_\Phi+(v\pa_yu)_\Phi-\p_{y}^2u_\Phi+\pa_xp_\Phi=0,
\end{split} \eeq
where $|D_x|^{s}$ denotes the Fourier multiplier in the horizontal variable with symbol $|\xi|^s.$

By applying $\D_k^\h$ to \eqref{S4eq2} and taking $L^2$ inner product of the resulting equation with $\D_k^\h (\p_tu)_\Phi,$ we find
\beq \label{S4eq5}
\begin{split}
\f12\f{d}{dt}&\|\D_k^{\rm h}(\p_tu)_\Phi(t)\|_{L^2}^2+\lam\dot{\tht}\||D_x|^{\f14}\D_k^{\rm
h}(\p_tu)_\Phi\|^2_{L^2}+\|\D_k^{\rm h}(\p_tu)_\Phi\|_{L^2}^2\\
&-\bigl(\D_k^\h\p_y^2 u_\Phi | \D_k^{\rm h}(\p_tu)_\Phi\bigr)_{L^2}
=-\bigl(\D_k^\h\left(u\p_xu\right)_\Phi | \D_k^{\rm h}(\p_tu)_\Phi\bigr)_{L^2}\\
&-\bigl(\D_k^\h\left(v\p_yu\right)_\Phi | \D_k^{\rm h}(\p_tu)_\Phi\bigr)_{L^2}-\bigl(\D_k^\h \p_x p_\Phi  | \D_k^{\rm h}(\p_tu)_\Phi\bigr)_{L^2}.
\end{split}
\eeq
Thanks to \eqref{S1eq2} and $\p_xu+\p_yv=0,$ we get, by using integration by parts, that
\beno
\begin{split}
\big(\D_k^\h\pa_xp_\Phi | \D_k^\h (\p_tu)_\Phi\big)_{L^2}=&-\big(\D_k^\h p_\Phi | \D_k^\h\pa_x(\p_tu)_\Phi\big)_{L^2}\\
=&\big(\D_k^\h p_\Phi | \D_k^\h\pa_y(\p_tv)_\Phi\big)_{L^2}
=-\big(\D_k^\h \pa_yp_\Phi | \D_k^\h (\p_tv)_\Phi\big)_{L^2}=0.
\end{split}
\eeno
Whereas observing from \eqref{S4eq1} that \beq \label{S4eq5p}
(\p_tu)_\Phi=\p_tu_\Phi+\lam\dot{\tht}|D_x|^{\f12}u_\Phi,\eeq one has
\begin{align*}
\lam\dot{\tht}\||D_x|^{\f14}\D_k^{\rm
h}(\p_tu)_\Phi\|^2_{L^2}=&\lam\dot{\tht}\||D_x|^{\f14}\D_k^{\rm
h}\p_tu_\Phi\|^2_{L^2}+\lam^3\dot{\tht}^3\||D_x|^{\f34}\D_k^{\rm
h}u_\Phi\|^2_{L^2}\\
&+2\lam^2\dot{\tht}^2\bigl(|D_x|^{\f14}\D_k^{\rm
h}\p_tu_\Phi\ |\ |D_x|^{\f34}\D_k^{\rm
h}u_\Phi\bigr)_{L^2},
\end{align*}
and
\begin{align*}
\bigl(|D_x|^{\f14}\D_k^{\rm
h}\p_tu_\Phi(t) | |D_x|^{\f34}\D_k^{\rm
h}u_\Phi(t)\bigr)_{L^2}=\f12\f{d}{dt}\||D_x|^{\f12}\D_k^{\rm h}u_\Phi(t)\|_{L^2}^2.
\end{align*}
Similarly, we find
\begin{align*}
-\bigl(\D_k^\h\p_y^2u_\Phi | \D_k^{\rm h}(\p_tu)_\Phi\bigr)_{L^2}=&-\bigl(\D_k^\h\p_y^2 u_\Phi | \D_k^{\rm h}(\p_tu_\Phi+\lam\dot{\tht}|D_x|^{\f12}u_\Phi
\bigr)_{L^2}\\
=&\f12\f{d}{dt}\|\D_k^{\rm h}\p_yu_\Phi(t)\|_{L^2}^2+\lam\dot{\tht}\||D_x|^{\f14}\D_k^{\rm
h} \p_yu_\Phi\|^2_{L^2}.
\end{align*}
By inserting the above estimates into \eqref{S4eq5}, we obtain
\beq\label{S4eq5a}
\begin{split}
\f{d}{dt}&\Bigl(\f12\|\D_k^{\rm h}(\p_tu)_\Phi(t)\|_{L^2}^2+\f12\|\D_k^{\rm h}\p_yu_\Phi(t)\|_{L^2}^2
+\lam^2\dot{\tht}^2(t)\||D_x|^{\f12}\D_k^{\rm h}u_\Phi(t)\|_{L^2}^2\Bigr)\\
&+\lam\dot{\tht}\bigl(\||D_x|^{\f14}\D_k^{\rm
h}\p_tu_\Phi\|^2_{L^2}+\||D_x|^{\f14}\D_k^{\rm
h} \p_yu_\Phi\|^2_{L^2}\bigr)\\
&-2\lam^2\dot{\tht}\ddot{\tht}\||D_x|^{\f12}\D_k^{\rm h}u_\Phi\|_{L^2}^2+\lam^3\dot{\tht}^3\||D_x|^{\f34}\D_k^{\rm
h}u_\Phi\|^2_{L^2}+\|\D_k^{\rm h}(\p_tu)_\Phi\|_{L^2}^2\\
=&-\bigl(\D_k^\h\left(u\p_xu\right)_\Phi | \D_k^{\rm h}(\p_tu)_\Phi\bigr)_{L^2}-\bigl(\D_k^\h\left(v\p_yu\right)_\Phi | \D_k^{\rm h}(\p_tu)_\Phi\bigr)_{L^2}.
\end{split}
\eeq

On the other hand, we get,
by applying $\D_k^\h$ to \eqref{S4eq2} and then taking $L^2$ inner product of the resulting equation with $\D_k^\h u_\Phi,$ that
\beq \label{S4eq5b}
\begin{split}
\bigl(&\D_k^{\rm h}\p_t(\p_tu)_\Phi | \D_k^\h u_\Phi\bigr)_{L^2}+\lam\dot{\tht}\bigl(|D_x|^{\f12}\D_k^{\rm
h}(\p_tu)_\Phi  | \D_k^\h u_\Phi\bigr)_{L^2}\\
&+\bigl(\D_k^{\rm h}(\p_tu)_\Phi  | \D_k^\h u_\Phi\bigr)_{L^2}-\bigl(\D_k^\h\p_y^2u_\Phi | \D_k^{\rm h} u_\Phi\bigr)_{L^2}\\
=&-\bigl(\D_k^\h\left(u\p_xu\right)_\Phi | \D_k^{\rm h} u_\Phi\bigr)_{L^2}-\bigl(\D_k^\h\left(v\p_yu\right)_\Phi | \D_k^{\rm h} u_\Phi\bigr)_{L^2}-\bigl(\D_k^\h \p_x p_\Phi  | \D_k^{\rm h} u_\Phi\bigr)_{L^2}.
\end{split}
\eeq
By using integration by parts, we find
\begin{align*}
\bigl(&\D_k^{\rm h}\p_t(\p_tu)_\Phi | \D_k^\h u_\Phi\bigr)_{L^2}=\f{d}{dt}\bigl(\D_k^{\rm h}(\p_tu)_\Phi | \D_k^\h u_\Phi\bigr)_{L^2}
-\bigl(\D_k^{\rm h}(\p_tu)_\Phi | \D_k^\h \p_tu_\Phi\bigr)_{L^2}\\
&=\f{d}{dt}\bigl(\D_k^{\rm h}(\p_tu)_\Phi | \D_k^\h u_\Phi\bigr)_{L^2}-\|\D_k^\h(\p_tu)_\Phi\|_{L^2}^2+\lam\dot{\tht}
\bigl(\D_k^\h (\p_tu)_\Phi\ |\ |D_x|^{\f12}\D_k^{\rm
h}u_\Phi\bigr)_{L^2},
\end{align*}
and
\beno
\bigl( \D_k^\h (\p_tu)_\Phi | |D_x|^{\f12}\D_k^{\rm
h}u_\Phi\bigr)_{L^2}=\f12\f{d}{dt}\||D_x|^{\f14}\D_k^{\rm h}u_\Phi(t)\|_{L^2}^2+\lam\dot{\tht}\||D_x|^{\f12}\D_k^{\rm
h}u_\Phi\|^2_{L^2}.
\eeno
Similarly, one has
\begin{align*}
\bigl(\D_k^{\rm h}(\p_tu)_\Phi  | \D_k^\h u_\Phi\bigr)_{L^2}=\f12\f{d}{dt}\|\D_k^{\rm h}u_\Phi(t)\|_{L^2}^2
+\lam\dot{\tht}\||D_x|^{\f14}\D_k^{\rm
h}u_\Phi\|^2_{L^2}.
\end{align*}
By substituting the above estimates into \eqref{S4eq5b}, we obtain
\beq
\label{S4eq5c}
\begin{split}
\f{d}{dt}&\Bigl(\bigl(\D_k^{\rm h}(\p_tu)_\Phi(t) | \D_k^\h u_\Phi(t)\bigr)_{L^2}+\lam\dot{\tht}(t)\||D_x|^{\f14}\D_k^{\rm h}u_\Phi(t)\|_{L^2}^2+\f12\|\D_k^{\rm h}u_\Phi(t)\|_{L^2}^2\Bigr)\\
&-\|\D_k^\h(\p_tu)_\Phi\|_{L^2}^2+\lam\bigl(\dot{\tht}-\ddot{\tht}\bigr)\||D_x|^{\f14}\D_k^{\rm
h}u_\Phi\|^2_{L^2}+\|\D_k^{\rm
h} \p_yu_\Phi\|^2_{L^2}\\
&+2\lam^2\dot{\tht}^2\||D_x|^{\f12}\D_k^{\rm h}u_\Phi\|_{L^2}^2
=-\bigl(\D_k^\h\left(u\p_xu\right)_\Phi | \D_k^{\rm h} u_\Phi\bigr)_{L^2}-\bigl(\D_k^\h\left(v\p_yu\right)_\Phi | \D_k^{\rm h} u_\Phi\bigr)_{L^2}.
\end{split}
\eeq

By summing up \eqref{S4eq5a} with $\f12\times$\eqref{S4eq5c}, we achieve
\beq\label{S4eq5d}
\begin{split}
\f{d}{dt}&\frak{F}_0(t)+\f\lam2\bigl(\dot{\tht}-\ddot{\tht}\bigr)\||D_x|^{\f14}\D_k^{\rm
h}u_\Phi\|^2_{L^2}+\lam\dot{\tht}\bigl(\||D_x|^{\f14}\D_k^{\rm
h}\p_tu_\Phi\|^2_{L^2}+\||D_x|^{\f14}\D_k^{\rm
h} \p_yu_\Phi\|^2_{L^2}\bigr)\\
&+\lam^2\bigl(\dot{\tht}^2-2\dot{\tht}\ddot{\tht}\bigr)\||D_x|^{\f12}\D_k^{\rm h}u_\Phi\|_{L^2}^2+\lam^3\dot{\tht}^3\||D_x|^{\f34}\D_k^{\rm
h}u_\Phi\|^2_{L^2}\\
&+\f12\bigl(\|\D_k^{\rm h}(\p_tu)_\Phi\|_{L^2}^2+\|\D_k^{\rm
h} \p_yu_\Phi\|^2_{L^2}\bigr)\\
=&-\bigl(\D_k^\h\left(u\p_xu\right)_\Phi | \D_k^{\rm h}\bigl(\p_tu+\f12u\bigr)_\Phi\bigr)_{L^2}-\bigl(\D_k^\h\left(v\p_yu\right)_\Phi | \D_k^{\rm h}\bigl(\p_t u +\f12u\bigr)_\Phi\bigr)_{L^2}.
\end{split}
\eeq
where $\frak{F}_0(t)$ is determined by
\begin{align*}
\frak{F}_0(t)\eqdefa &\f12\Bigl(\|\D_k^{\rm h}(\p_tu)_\Phi(t)\|_{L^2}^2+\|\D_k^{\rm h}\p_yu_\Phi(t)\|_{L^2}^2+\bigl(\D_k^{\rm h}(\p_tu)_\Phi(t) | \D_k^\h u_\Phi(t)\bigr)_{L^2}\\
&+\lam\dot{\tht}(t)\||D_x|^{\f14}\D_k^{\rm h}u_\Phi(t)\|_{L^2}^2\Bigr)+\f14\|\D_k^{\rm h}u_\Phi(t)\|_{L^2}^2
+\lam^2\dot{\tht}^2(t)\||D_x|^{\f12}\D_k^{\rm h}u_\Phi(t)\|_{L^2}^2.
\end{align*}

Observing that
\begin{align*}
\f13\|\D_k^{\rm h}(\p_tu)_\Phi(t)\|_{L^2}^2+\f12\bigl(\D_k^{\rm h}(\p_tu)_\Phi(t) |& \D_k^\h u_\Phi(t)\bigr)_{L^2}+\f3{16}\|\D_k^{\rm h}u_\Phi(t)\|_{L^2}^2\\
&=\f13\bigl\|\D_k^{\rm h}(\p_tu)_\Phi(t)+\f34\D_k^\h u_\Phi(t)\bigr\|_{L^2}^2\geq 0,
\end{align*}
we have
\beq\label{S4eq5e}
\begin{split}
\f16&\|\D_k^{\rm h}(\p_tu)_\Phi(t)\|_{L^2}^2+\f12\|\D_k^{\rm h}\p_yu_\Phi(t)\|_{L^2}^2+\f1{16}\|\D_k^{\rm h}u_\Phi(t)\|_{L^2}^2+\f\lam2\dot{\tht}(t)\||D_x|^{\f14}\D_k^{\rm h}u_\Phi(t)\|_{L^2}^2\\
&+\lam^2\dot{\tht}^2(t)\||D_x|^{\f12}\D_k^{\rm h}u_\Phi(t)\|_{L^2}^2
\leq \frak{F}_0(t)\leq \f34\|\D_k^{\rm h}(\p_tu)_\Phi(t)\|_{L^2}^2+\f12\|\D_k^{\rm h}\p_yu_\Phi(t)\|_{L^2}^2\\
&+\f12\|\D_k^{\rm h}u_\Phi(t)\|_{L^2}^2+\f\lam2\dot{\tht}(t)\||D_x|^{\f14}\D_k^{\rm h}u_\Phi(t)\|_{L^2}^2
+\lam^2\dot{\tht}^2(t)\||D_x|^{\f12}\D_k^{\rm h}u_\Phi(t)\|_{L^2}^2.
\end{split}\eeq

While due to $\left(u_{\Phi}, v_{\Phi}\right)|_{y=0}=\left(u_\Phi, v_\Phi\right)|_{y=1}=0,$ we get, by applying Poincar\'e inequality, that
\beq \label{S4eq5f}
\|\D_k^\h u_\Phi\|_{L^2}^2\leq {\cK}\|\p_y\D_k^\h u_\Phi\|_{L^2}^2.
\eeq
Moreover, by virtue of \eqref{S4eq-1}, one has $\ddot\tht(t)=-\f{\fk}2\dot\tht(t),$
then
by  \eqref{S4eq5f}, if we take $\fk\eqdefa \min\left(\f1{6},\f1{4(1+\cK)}\right),$ we find
\beq \label{S4eq5g}
\begin{split}
\lam\bigl(\f{\dot{\tht}(t)}4-\f{\ddot{\tht}(t)}2\bigr)\||D_x|^{\f14}&\D_k^{\rm
h} u_\Phi(t)\|^2_{L^2}+\lam^2(\f{\dot\tht(t)^2}{2}-2\dot{\tht}(t)\ddot\tht(t))\||D_x|^{\f12}\D_k^{\rm h}u_\Phi(t)\|_{L^2}^2\\
&+\f14\bigl(\|\D_k^{\rm h}(\p_tu)_\Phi\|_{L^2}^2+\|\D_k^{\rm
h} \p_yu_\Phi\|^2_{L^2}\bigr)-2\fk\frak{F}_0(t)\geq 0.
\end{split}
\eeq
Then by multiplying \eqref{S4eq5d} by $e^{2\frak{K}t}$  and  integrating the resulting inequality over $[0,t],$ we achieve
\beq \label{S4eq6rt}
\begin{split}
\frak{G}_k(u)(t)
\leq & \f3{4}\bigl\|e^{a|D_x|^{\f12}}\D_k^{\rm h}(\p_tu)(0)\bigr\|_{L^2}^2+\f12\bigl\|e^{a|D_x|^{\f12}}\D_k^{\rm h}\p_yu_0\bigr\|_{L^2}^2
+\f1{2}\bigl\|e^{a|D_x|^{\f12}}\D_k^{\rm h}u_0\bigr\|_{L^2}^2\\
&+\f\lam2\de^{\f12}2^{\f{k}2}\bigl\|e^{a|D_x|^{\f12}}\D_k^{\rm h}u_0\bigr\|_{L^2}^2+{\lam^2}\de2^k\bigl\|e^{a|D_x|^{\f12}}\D_k^{\rm h}u_0\bigr\|_{L^2}^2\\
&+\int_0^t\bigl|\bigl(\ektp\D_k^\h\left(u\p_xu\right)_\Phi | \ektp\D_k^\h \bigl(\p_tu+\f12u\bigr)_\Phi\bigr)_{L^2}\bigr|\,dt'\\
&\qquad\qquad\qquad\quad+\int_0^t\bigl|\bigl(\ektp\D_k^\h\left(v\p_yu\right)_\Phi | \ektp\D_k^\h \bigl(\p_tu+\f12u\bigr)_\Phi\bigr)_{L^2}\bigr|\,dt',
\end{split}
\eeq
with $\frak{G}_{k}(u)(t)$ being determined by
\beq \label{S4eq6a}
\begin{split}
\frak{G}_k&(u)(t)\eqdefa \f16\|\ektp\D_k^{\rm h}(\p_tu)_\Phi\|_{L^\infty_t(L^2)}^2+\f12\|\ektp\D_k^{\rm h}\p_yu_\Phi\|_{L^\infty_t(L^2)}^2+\f1{16}\|\ektp\D_k^{\rm h}u_\Phi\|_{L^\infty_t(L^2)}^2\\
&+\f\lam2\dot{\tht}(t)\|\ektp|D_x|^{\f14}\D_k^{\rm h}u_\Phi(t)\|_{L^2}^2+{\lam^2}\dot{\tht}^2(t)\|\ektp|D_x|^{\f12}\D_k^{\rm h}u_\Phi(t)\|_{L^2}^2\\
&+\lam \int_0^t\dot{\tht}(t')\Bigl(\f14\|\ektp|D_x|^{\f14}\D_k^{\rm
h}u_\Phi\|^2_{L^2}+\bigl\|\ektp|D_x|^{\f14}\D_k^{\rm
h}(\p_tu_\Phi, \p_yu_\Phi)(t')\bigr\|^2_{L^2}\Bigr)\,dt'\\
&+\f{{\lam^2}}{2}\int_0^t\dot{\tht}^2(t')\|\ektp|D_x|^{\f12}\D_k^{\rm h}u_\Phi(t')\|_{L^2}^2\,dt'+\lam^3
\int_0^t\dot{\tht}^3(t')\|\ektp|D_x|^{\f34}\D_k^{\rm
h}u_\Phi(t')\|^2_{L^2}\,dt'\\
&+\f14\bigl(\|\ektp\D_k^{\rm h}(\p_tu)_\Phi\|_{L^2_t(L^2)}^2+\|\D_k^{\rm
h} \p_yu_\Phi\|^2_{L^2_t(L^2)}\bigr).
\end{split}
\eeq

In what follows, we shall always assume that $t<T^\star_0$ with
$T^\star_0$ being determined by
\beq\label{eq4.3} T^\star_0\eqdefa
\sup\bigl\{\ t>0,\ \  \|\p_y u_\Phi(t)\|_{\cB^{\f12}}\leq \de e^{-\frak{K} t} \  \bigr\}. \eeq

The estimates of the last two terms in \eqref{S4eq6rt} rely on the following two lemmas:

\begin{lem}\label{lem3.1}
{\sl For any $s\in \left]-\f14,\f34\right]$ and $t\leq T^\star_0,$ there holds
\beq \label{S4eq7wq}
\begin{split}
\int_0^t\bigl|\bigl(\ektp\D_k^\h(u\p_x \fa)_\Phi\ &|\ \ektp\D_k^\h \fb_\Phi\bigr)_{L^2}\bigr|\,dt'\\
&\lesssim d_k^2
2^{-2ks}\|\ektp \fa_\Phi\|_{\wt{L}^2_{t,\dot{\tht}^3(t)}(\cB^{s+\frac34})}\|\ektp \fb_\Phi\|_{\wt{L}^2_{t,\dot{\tht}(t)}(\cB^{s+\frac14})},
\end{split}
\eeq where the time-weighted Chemin-Lerner type norms $\|\cdot\|_{\wt{L}^2_{t,\dot{\tht}(t)}(\cB^{s+\frac14})}$ and
$\|\cdot\|_{\wt{L}^2_{t,\dot{\tht}^3(t)}(\cB^{s+\frac34})}$ will be given by Definition \ref{def1.1}.
In particular, when $\fa=u$ or $\fa=\p_yu,$ \eqref{S4eq7wq} holds for any $s>-\f14.$}
\end{lem}

\begin{lem}\label{lem3.2}
{\sl Let us assume that $\|\fa_\Phi(t)\|_{\cB^{\f12}}\leq \de e^{-\fk t}$ for $t\leq T^\star_\fa.$ Then  for any $s\in \left]-\f14,\f34\right]$ and $t\leq T^\star_\fa,$ there holds
\beq \label{S4eq9g}\begin{split}
\int_0^t\bigl|\bigl(\ektp\D_k^\h(v\fa)_\Phi\ &|\ \ektp\D_k^\h \fb_\Phi\bigr)_{L^2}\bigr|\,dt'\\
&\lesssim d_k^2
2^{-2ks}\|\ektp u_\Phi\|_{\wt{L}^2_{t,\dot{\tht}^3(t)}(\cB^{s+\frac34})}\|\ektp \fb_\Phi\|_{\wt{L}^2_{t,\dot{\tht}(t)}(\cB^{s+\frac14})}.
\end{split}
\eeq
Furthermore for $t\leq T^\star_{1,\fa}\eqdefa \sup \left\{ \ t\leq T^\star_\fa, \ \|u_\Phi(t)\|_{\cB^{1}}\leq \de^{\f12} e^{-\f{\fk}2 t}\ \right\}$ and
$s>-\f14,$ one has
\beq \label{S4eq9gq}\begin{split}
\int_0^t&\bigl|\bigl(\ektp\D_k^\h(v\fa)_\Phi\ |\ \ektp\D_k^\h \fb_\Phi\bigr)_{L^2}\bigr|\,dt'\\
&\lesssim d_k^2
2^{-2ks}\bigl(\|\ektp u_\Phi\|_{\wt{L}^2_{t,\dot{\tht}^3(t)}(\cB^{s+\frac34})}+\|\ektp \fa_\Phi\|_{\wt{L}^2_{t,\dot{\tht}(t)}(\cB^{s+\frac14})}\bigr)
\|\ektp \fb_\Phi\|_{\wt{L}^2_{t,\dot{\tht}(t)}(\cB^{s+\frac14})}.
\end{split}
\eeq
}
\end{lem}

We admit the above lemmas for the time being and continue our proof of Theorem \ref{th1.2}. Indeed
we first  deduce from Lemma \ref{lem3.1} that for $t\leq T^\star_0$ and $s\in \left]-\f14,\f34\right]$
\beq \label{S4eq9a}
\begin{split}
\int_0^t\bigl|\bigl(\ektp\D_k^\h(u\p_x u)_\Phi\ &|\ \ektp\D_k^\h u_\Phi\bigr)_{L^2}\bigr|\,dt'\\
&\lesssim d_k^2
2^{-2ks}\|\ektp u_\Phi\|_{\wt{L}^2_{t,\dot{\tht}^3(t)}(\cB^{s+\frac34})}\|\ektp u_\Phi\|_{\wt{L}^2_{t,\dot{\tht}(t)}(\cB^{s+\frac14})},
\end{split}
\eeq
and
\beq \label{S4eq9b}
\begin{split}
\int_0^t&\bigl|\bigl(\ektp\D_k^\h\left(u\p_xu\right)_\Phi | \ektp\D_k^\h (\p_tu)_\Phi\bigr)_{L^2}\bigr|\,dt'\\
&=\int_0^t\bigl|\bigl(\ektp\D_k^\h\left(u\p_xu\right)_\Phi | \ektp\D_k^\h \bigl(\p_tu_\Phi+\lam\dot{\tht}(t){|D_x|}^{\f12}u_\Phi\bigr)\bigr)_{L^2}\bigr|\,dt'\\
&\lesssim d_k^2
2^{-2ks}\|\ektp u_\Phi\|_{\wt{L}^2_{t,\dot{\tht}^3(t)}(\cB^{s+\frac34})}\bigl(\|\ektp \p_tu_\Phi\|_{\wt{L}^2_{t,\dot{\tht}(t)}(\cB^{s+\frac14})}
+\lam \|\ektp u_\Phi\|_{\wt{L}^2_{t,\dot{\tht}^3(t)}(\cB^{s+\frac34})}\bigr).
\end{split}
\eeq

Whereas it follows from Lemma \ref{lem3.2} that for $t\leq T^\star_0$ and $s\in \left]-\f14,\f34\right]$
\beq \label{S4eq9c}
\begin{split}
\int_0^t\bigl|\bigl(\ektp\D_k^\h(v\p_y u)_\Phi\ &|\ \ektp\D_k^\h u_\Phi\bigr)_{L^2}\bigr|\,dt'\\
&\lesssim d_k^2 2^{-2ks}\|\ektp u_\Phi\|_{\wt{L}^2_{t,\dot{\tht}^3(t)}(\cB^{s+\frac34})}\|\ektp u_\Phi\|_{\wt{L}^2_{t,\dot{\tht}(t)}(\cB^{s+\frac14})},
\end{split}
\eeq
and
\beq \label{S4eq9dqw}
\begin{split}
\int_0^t&\bigl|\bigl(\ektp\D_k^\h\left(v\p_yu\right)_\Phi | \ektp\D_k^\h (\p_tu)_\Phi\bigr)_{L^2}\bigr|\,dt'\\
&\lesssim d_k^2
2^{-2ks}\|\ektp u_\Phi\|_{\wt{L}^2_{t,\dot{\tht}^3(t)}(\cB^{s+\frac34})}\bigl(\|\ektp \p_tu_\Phi\|_{\wt{L}^2_{t,\dot{\tht}(t)}(\cB^{s+\frac14})}
+\lam \|\ektp u_\Phi\|_{\wt{L}^2_{t,\dot{\tht}^3(t)}(\cB^{s+\frac34})}\bigr).
\end{split}\eeq

Without loss of generality, we may assume that $\lam\geq 1.$ Then
by
inserting the above estimates into \eqref{S4eq6rt}, we find that for $t\leq T^\star_0$ and $s\in \left]-\f14,\f34\right]$
\beno
\begin{split}
\frak{G}_k(u)(t)
\leq & \f3{4}\bigl\|e^{a|D_x|^{\f12}}\D_k^{\rm h}u_1\bigr\|_{L^2}^2+\f12\bigl\|e^{a|D_x|^{\f12}}\D_k^{\rm h}\p_yu_0\bigr\|_{L^2}^2
+\f12\bigl\|e^{a|D_x|^{\f12}}\D_k^{\rm h}u_0\bigr\|_{L^2}^2\\
&+\f\lam2\de^{\f12}2^{\f{k}2}\bigl\|e^{a|D_x|^{\f12}}\D_k^{\rm h}u_0\bigr\|_{L^2}^2+{\lam^2}\de2^k\bigl\|e^{a|D_x|^{\f12}}\D_k^{\rm h}u_0\bigr\|_{L^2}^2\\
&+ Cd_k^2
2^{-2ks}\Bigl(\bigl\|\ektp (u_\Phi,\p_tu_\Phi)\bigr\|_{\wt{L}^2_{t,\dot{\tht}(t)}(\cB^{s+\frac14})}^2+\lam\|\ektp u_\Phi\|_{\wt{L}^2_{t,\dot{\tht}^3(t)}(\cB^{s+\frac34})}^2
\Bigr).
\end{split}
\eeno
Then thanks to \eqref{S4eq-1}, by multiplying the above inequality by $2^{2ks}$ and then
taking square root of the resulting inequality, and finally by summing up the resulting ones over $\Z,$ we obtain
\beq \label{AAAQ}
\begin{split}
&{\rm E}_{s,\lam}(u)(t) \leq C\Bigl(\bigl\|e^{a|D_x|^{\f12}}(u_0,u_1,\p_yu_0)\bigr\|_{\cB^{s}}+\sqrt{a\fk}\bigl\|e^{a|D_x|^{\f12}}u_0\bigr\|_{\cB^{s+\f14}}\\
&\quad+
 {a\fk}\bigl\|e^{a|D_x|^{\f12}}u_0\bigr\|_{\cB^{s+\f12}}+\bigl\|\ektp \bigl(u_\Phi, \p_tu_\Phi\bigr)\bigr\|_{\wt{L}^2_{t,\dot\tht(t)}(\cB^{s+\f14})}
+\lam^{\f12}\|\ektp u_\Phi\|_{\wt{L}^2_{t,\dot\tht^3(t)}(\cB^{s+\f34})}\Bigr),
\end{split}
\eeq for $t\leq T^\star_0$ and $s\in \left]-\f14,\f34\right],$
where we used \eqref{S4eq-1} so that $ \lam\de^{\f12}=\f{a\fk}4,$ and
\beq\label{Esu}
\begin{split}
{\rm E}_{s,\lam}&(u)(t)\eqdefa
 \bigl\|\ektp \bigl(u_\Phi, \p_yu_\Phi, (\p_tu)_\Phi\bigr)\bigr\|_{\wt{L}^\infty_t(\cB^{s})}+\sqrt{a\fk}\|e^{\f34\fk t'} u_\Phi\|_{\wt{L}^\infty_{t}(\cB^{s+\f14})}\\
&+a\fk\|e^{\f12\fk t'} u_\Phi\|_{\wt{L}^\infty_{t}(\cB^{s+\f12})}+\sqrt{\lam}\bigl\|\ektp \bigl(u_\Phi, \p_tu_\Phi, \p_yu_\Phi\bigr)\bigr\|_{\wt{L}^2_{t,\dot\tht(t)}(\cB^{s+\f14})}\\
&+{\lam}\|\ektp u_\Phi\|_{\wt{L}^2_{t,\dot\tht^2(t)}(\cB^{s+\f12})}+{\lam}^{\f32}\|\ektp u_\Phi\|_{\wt{L}^2_{t,\dot\tht^3(t)}(\cB^{s+\f34})}
+\bigl\|\ektp\bigl((\p_t u)_\Phi, \p_y u_\Phi\bigr)\bigr\|_{\wt{L}^2_t(\cB^{s})}.
\end{split}
\eeq

Taking $\lam =4C^2$ in \eqref{AAAQ} leads to
\beq \label{Esuu}
{\rm E}_{s,\f{\lam}4}(t)\leq C\frak{H}^0_{s}(u_0,u_1),
\eeq
for   $t\leq T^\star_0$ and $s\in \left]-\f14,\f34\right],$ and $\frak{H}^0_{s}(u_0,u_1)$ being determined by \eqref{S1eq3}.
In particular, we deduce from \eqref{S1eq3} and \eqref{Esuu}  that for   $t\leq T^\star_0$
\beq \label{Esuuv}
{\rm E}_{\f12,\f{\lam}4}(t)\leq C\frak{H}^0_{\f12}(u_0,u_1)\leq Cc_0,
\eeq
which implies
\beq\label{S4eq19}
\|\p_yu_\Phi(t)\|_{\cB^{\f12}}\leq Cc_0 e^{-\frak{K}t}\leq \f\de2e^{-\frak{K}t} \quad \forall\ t\leq T^\star_0,
\eeq
if we take $c_0$ in \eqref{S1eq1} to be so small that
$ Cc_0\leq \f\de2.$ Then
we deduce by a continuous argument that $T^\star_0$ determined by \eqref{eq4.3} equals $+\infty$ and \eqref{S4eq10} holds.
This completes the proof of part (1) of Theorem \ref{th1.2}.

\no{\bf Part (2).} We also deduce from  \eqref{Esu} and \eqref{Esuuv} that
\beno
a\fk \|e^{\f{\fk}2t'}u_\Phi\|_{L^\infty(\R^+;\cB^1)}\leq Cc_0.
\eeno
By taking $c_0$ in \eqref{S1eq1} to be so small that $\f{Cc_0}{a\fk}\leq \de^{\f12},$ we deduce that
\beq \label{S4eq19ui}
\|u_\Phi(t)\|_{\cB^1}\leq \de^{\f12} e^{-\f{\fk}2t}\quad \forall\ t\in\R^+.
\eeq
Thanks to \eqref{S4eq19ui}, we get, by applying \eqref{S4eq9gq}, that for any $s>-\f14$
\begin{align*}
\int_0^t&\bigl|\bigl(\ektp\D_k^\h\left(v\p_yu\right)_\Phi | \ektp\D_k^\h \Bigl(\p_tu+\f12u)_\Phi\bigr)_{L^2}\bigr|\,dt'\lesssim d_k^2
2^{-2ks}\bigl(\|\ektp u_\Phi\|_{\wt{L}^2_{t,\dot{\tht}^3(t)}(\cB^{s+\frac34})}\\
&\qquad+\|\ektp \p_y u_\Phi\|_{\wt{L}^2_{t,\dot{\tht}^3(t)}(\cB^{s+\frac14})}\Bigr)\Bigl(\|\ektp \p_tu_\Phi\|_{\wt{L}^2_{t,\dot{\tht}(t)}(\cB^{s+\frac14})}
+\lam \|\ektp u_\Phi\|_{\wt{L}^2_{t,\dot{\tht}^3(t)}(\cB^{s+\frac34})}\Bigr).
\end{align*}
Then along the same line to the proof \eqref{AAAQ}, we deduce that for any $s>-\f14$
\begin{align*}
&{\rm E}_{s,\lam}(t) \leq C\Bigl(\frak{H}^0_s(u_0,u_1)+\bigl\|\ektp \bigl(u_\Phi, \p_yu_\Phi, \p_tu_\Phi\bigr)\bigr\|_{\wt{L}^2_{t,\dot\tht(t)}(\cB^{s+\f14})}
+\lam\|\ektp u_\Phi\|_{\wt{L}^2_{t,\dot\tht^3(t)}(\cB^{s+\f34})}\Bigr),
\end{align*}
Taking $\lam=4C^2$ in the above inequality gives rise to
\beno
{\rm E}_{s,\f\lam4}(t) \leq C\frak{H}^0_s(u_0,u_1)\quad \mbox{for} \ s>-\f14,
\eeno
which implies \eqref{S1eq12}.
This finishes the proof of part (2) of Theorem \ref{th1.2}.
\end{proof}

Now let us present the proof of Lemmas \ref{lem3.1} and \ref{lem3.2}. Indeed, we observe that it amounts to prove these lemmas for $\frak{K}=0.$ Without loss of generality, we may assume that $\widehat{\fa}\ge 0$ and $\widehat{\fb}\ge0$ (and similar assumption for the proof of the product law in the rest of this paper, one may check \cite{CGP} for detail).

\begin{proof}[Proof of Lemma \ref{lem3.1}] We first get, by applying Bony's decomposition \eqref{Bony}  in the horizontal variable to $u\p_x\fa$, that
\beno
u\p_x\fa=T^\h_{u}\p_x\fa+T^\h_{\p_x\fa}u+R^h(u,\p_x\fa).
\eeno
Accordingly, we shall handle the following three terms:\smallskip

\no $\bullet$ \underline{Estimate of
$\int_0^t\bigl|\bigl(\D_k^{\rm h}(T^\h_{u}\p_x\fa)_\Phi\ |\ \D_k^{\rm
h}\fb_\Phi\bigr)_{L^2}\bigr|\,dt'$}

Considering the support properties to the Fourier transform of the terms in $T^\h_{u}\p_x\fa,$ we find
\beno
\begin{split}
\int_0^t\bigl|\bigl(\D_k^{\rm h}(T^\h_{u}\p_x\fa)_\Phi\ |&\ \D_k^{\rm
h}\fb_\Phi\bigr)_{L^2}\bigr|\,dt'\\
\lesssim & \sum_{|k'-k|\leq 4}\int_0^t\|S_{k'-1}^\h u_\Phi(t')\|_{L^\infty}
\|\D_{k'}^\h\p_x\fa_\Phi(t')\|_{L^2}\|\D_k^\h \fb_\Phi(t')\|_{L^2}\,dt'.
\end{split}
\eeno
However, due to $u|_{y=0}=0,$ we write $u(t,x,y)=\int_0^y\p_yu(t,x,y')\,dy',$ so that we deduce   from Lemma \ref{lem:Bern} and Poincar\'e inequality that
\beq\label{S3eq8}
\begin{split}
\|\D_k^\h u_\Phi(t)\|_{L^\infty}\lesssim & 2^{\frac{k}2}\|\D_k^\h u_\Phi(t)\|_{L^2_\h (L^\infty_{\rm v})}\\
\lesssim& 2^{\frac{k}2}\|\D_k^\h u_\Phi(t)\|_{L^2}^{\f12}\|\D_k^\h \p_y u_\Phi(t)\|_{L^2}^{\f12}\\
\lesssim & 2^{\frac{k}2} \|\D_k^\h \p_y u_\Phi(t)\|_{L^2}\lesssim d_j(t)\|\p_y u_\Phi(t)\|_{\cB^{\f12}},
\end{split}
\eeq
which implies
\beno
\|S_{k'-1}^\h u_\Phi(t)\|_{L^\infty}\lesssim \|\p_yu_\Phi(t)\|_{\cB^{\f12}}.
\eeno
This together with \eqref{S4eq-1} and \eqref{eq4.3} ensures that for $t\leq T^\star_0,$
\beno
\begin{split}
\int_0^t\bigl|\bigl(\D_k^{\rm h}(T^\h_{u}\p_x\fa)_\Phi\ |&\ \D_k^{\rm
h}\fb_\Phi\bigr)_{L^2}\bigr|\,dt'\\
\lesssim & \sum_{|k'-k|\leq 4}2^{k'}\int_0^t\|\p_yu_\Phi(t')\|_{\cB^{\f12}}
\|\D_{k'}^\h \fa_\Phi(t')\|_{L^2}\|\D_k^\h \fb_\Phi(t')\|_{L^2}\,dt'\\
\lesssim & \sum_{|k'-k|\leq 4}2^{k'}\int_0^t \dot{\tht}^2(t')
\|\D_{k'}^\h \fa_\Phi(t')\|_{L^2}\|\D_k^\h \fb_\Phi(t')\|_{L^2}\,dt'.
\end{split}
\eeno
Then we get, by applying H\"older inequality and using Definition \ref{def1.1}, that
\beq \label{S3eq8wr}
\begin{split}
\int_0^t&\bigl|\bigl(\D_k^{\rm h}(T^\h_{u}\p_x\fa)_\Phi\ |\ \D_k^{\rm
h}\fb_\Phi\bigr)_{L^2}\bigr|\,dt'\\
\lesssim & \sum_{|k'-k|\leq 4}2^{k'}\Bigl(\int_0^t\dot{\tht}^3(t')\|\D_{k'}^\h \fa_\Phi(t')\|_{L^2}^2\,dt'\Bigr)^{\f12}\Bigl(\int_0^t\dot{\tht}(t') \|\D_k^\h \fb_\Phi(t')\|_{L^2}^2\,dt'\Bigr)^{\f12}\\
\lesssim & d_k2^{-2ks}\|\fa_\Phi\|_{\wt{L}_{t,\dot{\tht}^3(t)}(\cB^{s+\frac34})}\|\fb_\Phi\|_{\wt{L}_{t,\dot{\tht}(t)}(\cB^{s+\frac14})}
\Bigl(\sum_{|k'-k|\leq 4}d_{k'}2^{(k-k')\left(s-\f14\right)}\Bigr)\\
\lesssim & d_k^22^{-2ks}\|\fa_\Phi\|_{\wt{L}_{t,\dot{\tht}^3(t)}(\cB^{s+\frac34})}\|\fb_\Phi\|_{\wt{L}_{t,\dot{\tht}(t)}(\cB^{s+\frac14})}.
\end{split}
\eeq

\no $\bullet$ \underline{Estimate of
$\int_0^t\bigl|\bigl(\D_k^{\rm h}(T^\h_{\p_x\fa}u)_\Phi\ |\ \D_k^{\rm
h}\fb_\Phi\bigr)_{L^2}\bigr|\,dt'$}

Again considering the support properties to the Fourier transform of the terms in $T^\h_{\p_x\fa}u$ and thanks to  \eqref{S3eq8},
we deduce that for $t\leq T^\star_0$
\beq\label{S3eq8ap}
\begin{split}
\int_0^t\bigl|\bigl(&\D_k^{\rm h}(T^\h_{\p_x
\fa}u)_\Phi\ |\ \D_k^{\rm
h}\fb_\Phi\bigr)_{L^2}\bigr|\,dt'\\
\lesssim & \sum_{|k'-k|\leq 4}\int_0^t\|S_{k'-1}^\h \p_x\fa_\Phi(t')\|_{L^\infty_\h(L^2_{\rm v})}
\|\D_{k'}^\h u_\Phi(t')\|_{L^2_\h(L^\infty_{\rm v})}\|\D_k^\h \fb_\Phi(t')\|_{L^2}\,dt'\\
\lesssim & \sum_{|k'-k|\leq 4}2^{-\frac{k'}2}\int_0^td_{k'}(t)\|S_{k'-1}^\h\p_x\fa_\Phi(t')\|_{L^\infty_\h(L^2_{\rm v})}\|\p_yu_\Phi(t')\|_{\cB^{\f12}}
\|\D_k^\h \fb_\Phi(t')\|_{L^2}\,dt'\\
\lesssim & \sum_{|k'-k|\leq 4}2^{-\frac{k'}2}\Bigl(\int_0^t\|S_{k'-1}^\h\p_x\fa_\Phi(t')\|_{L^\infty_\h(L^2_{\rm v})}^2\dot{\tht}^3(t')\,dt'\Bigr)^{\frac12}\\
&\qquad\qquad\qquad\qquad\qquad\qquad\qquad\times \Bigl(\int_0^t\|\D_k^\h \fb_\Phi(t')\|_{L^2}^2\dot{\tht}(t')\,dt'\Bigr)^{\frac12}.
\end{split} \eeq
Yet we observe from Definition \ref{def1.1} and $s< \f34$ that
\beq\label{S3eq8ah}
\begin{split}
\Bigl(\int_0^t\|S_{k'-1}^\h\p_x\fa_\Phi(t')\|_{L^\infty_\h(L^2_{\rm v})}^2\dot{\tht}^3(t')\,dt'\Bigr)^{\frac12}
\lesssim
&\sum_{\ell\leq k'-2}2^{\frac{3\ell}2}\Bigl(\int_0^t\|\D_\ell^\h \fa_\Phi(t')\|_{L^2}^2\dot{\tht}^3(t')\,dt'\Bigr)^{\frac12}\\
\lesssim
&\sum_{\ell\leq k'-2}d_\ell 2^{\ell\left(\f34-s\right)}\|\fa_\Phi\|_{\wt{L}^2_{t,\dot{\tht}^3(t)}(\cB^{s+\frac34})}\\
\lesssim & d_{k'}2^{k'\left(\f34-s\right)}\|\fa_\Phi\|_{\wt{L}^2_{t,\dot{\tht}^3(t)}(\cB^{s+\frac34})}.
\end{split}
\eeq
As a result, for $s< \f34$, it comes out
\beq \label{S3eq8au}
\begin{split}
\int_0^t\bigl|\bigl(\D_k^{\rm h}(T^\h_{\p_x\fa}u)_\Phi\ |\ \D_k^{\rm
h}\fb_\Phi\bigr)_{L^2}\bigr|\,dt' \lesssim & d_k^22^{-2ks}\|\fa_\Phi\|_{\wt{L}_{t,\dot{\tht}^3(t)}(\cB^{s+\frac34})}\|\fb_\Phi\|_{\wt{L}_{t,\dot{\tht}(t)}(\cB^{s+\frac14})}.
\end{split}
\eeq

On the other hand, we get, by a similar derivation of \eqref{S3eq8ap} that
\begin{align*}
\int_0^t\bigl|\bigl(&\D_k^{\rm h}(T^\h_{\p_x
\fa}u)_\Phi\ |\ \D_k^{\rm
h}\fb_\Phi\bigr)_{L^2}\bigr|\,dt'\\
\lesssim & \sum_{|k'-k|\leq 4}2^{-\frac{k'}2}\int_0^td_{k'}(t)\|\p_x\fa_\Phi(t')\|_{L^\infty_\h(L^2_{\rm v})}\|\p_yu_\Phi(t')\|_{\cB^{\f12}}
\|\D_k^\h \fb_\Phi(t')\|_{L^2}\,dt'\\
\lesssim & \sum_{|k'-k|\leq 4}2^{-\frac{k'}2}d_{k'}\Bigl(\int_0^t\|\p_x\fa_\Phi(t')\|_{L^\infty_\h(L^2_{\rm v})}^2\dot{\tht}^3(t')\,dt'\Bigr)^{\frac12}\Bigl(\int_0^t\|\D_k^\h \fb_\Phi(t')\|_{L^2}^2\dot{\tht}(t')\,dt'\Bigr)^{\frac12}.
\end{align*}
While it follows from the derivation of \eqref{S3eq8ah} that
\begin{align*}
\Bigl(\int_0^t\|\p_x\fa_\Phi(t')\|_{L^\infty_\h(L^2_{\rm v})}^2\dot{\tht}^3(t')\,dt'\Bigr)^{\frac12}
\lesssim \|\fa_\Phi\|_{\wt{L}^2_{t,\dot{\tht}^3(t)}(\cB^{\frac32})}.
\end{align*}
Therefore \eqref{S3eq8au} holds for $s=\f34.$\\

\no $\bullet$ \underline{Estimate of
$\int_0^t\bigl|\bigl(\D_k^{\rm h}(R^\h(u,\p_x\fa))_\Phi\ |\ \D_k^{\rm
h}\fb_\Phi\bigr)_{L^2}\bigr|\,dt'$}\vspace{0.2cm}

Again considering the support properties to the Fourier transform of the terms in $R^\h(u,\p_x\fa),$ we get, by applying Lemma
\ref{lem:Bern} and \eqref{S3eq8}, that for $t\leq T^\star_0$
\beno
\begin{split}
\int_0^t\bigl|&\bigl(\D_k^{\rm h}(R^\h(u,\p_x\fa))_\Phi\ |\ \D_k^{\rm
h}\fb_\Phi\bigr)_{L^2}\bigr|\,dt'\\
\lesssim &2^{\f{k}2}\sum_{k'\geq k-3}\int_0^t\|\wt{\D}_{k'}^\h u_\Phi(t')\|_{L^2_\h(L^\infty_{\rm v})}\|{\D}_{k'}^\h \p_x \fa_\Phi(t')\|_{L^2}\|\D_k^\h \fb_\Phi(t')\|_{L^2}\,dt'\\
\lesssim & 2^{\f{k}2}\sum_{k'\geq k-3}2^{\f{k'}2}\int_0^t\dot{\tht}^2(t')\|{\D}_{k'}^\h  \fa_\Phi(t')\|_{L^2}\|\D_k^\h \fb_\Phi(t')\|_{L^2}\,dt'.
\end{split}
\eeno
Due to $s>-\f14,$ by applying H\"older inequality and using Definition \ref{def1.1}, we obtain
\beq \label{S3eq8hk}
\begin{split}
\int_0^t\bigl|&\bigl(\D_k^{\rm h}(R^\h(u,\p_x\fa))_\Phi\ |\ \D_k^{\rm
h}\fb_\Phi\bigr)_{L^2}\bigr|\,dt'\\
\lesssim & 2^{\f{k}2}\sum_{k'\geq k-3}2^{\f{k'}2}\Bigl(\int_0^t\|{\D}_{k'}^\h  \fa_\Phi(t')\|_{L^2}^2\dot{\tht}^3(t')\,dt'\Bigr)^{\f12}
 \Bigl(\int_0^t\|\D_k^\h \fb_\Phi(t')\|_{L^2}^2\dot{\tht}(t')\,dt'\Bigr)^{\f12}\\
\lesssim &d_k2^{-2{k}s}\|\fa_\Phi\|_{\wt{L}_{t,\dot{\tht}^3(t)}(\cB^{s+\frac34})}\|\fb_\Phi\|_{\wt{L}_{t,\dot{\tht}(t)}(\cB^{s+\frac14})} \Bigl(\sum_{k'\geq k-3}d_{k'}2^{(k-k')\left(s+\f14\right)}\Bigr)\\
\lesssim &d_k^22^{-2ks}\|\fa_\Phi\|_{\wt{L}_{t,\dot{\tht}^3(t)}(\cB^{s+\frac34})}\|\fb_\Phi\|_{\wt{L}_{t,\dot{\tht}(t)}(\cB^{s+\frac14})}.
\end{split}
\eeq

By summing up the estimates \eqref{S3eq8wr}, \eqref{S3eq8au} and \eqref{S3eq8hk}, we conclude the proof of \eqref{S4eq7wq} for $s\in \left]-\f14,\f34\right].$

On the other hand, we observe from \eqref{S3eq8} that \beno
\|S_{k'-1}^\h \p_x u_\Phi(t)\|_{L^\infty}\lesssim 2^k\|\p_yu_\Phi(t)\|_{\cB^{\f12}},
\eeno
so that we get, by a similar derivation of \eqref{S3eq8ap} that for any $s\in\R,$
\begin{align*}
\int_0^t\bigl|\bigl(&\D_k^{\rm h}(T^\h_{\p_x
u}u)_\Phi\ |\ \D_k^{\rm
h}\fb_\Phi\bigr)_{L^2}\bigr|\,dt'\\
\lesssim & \sum_{|k'-k|\leq 4}\int_0^t\|S_{k'-1}^\h \p_xu_\Phi(t')\|_{L^\infty_\h}
\|\D_{k'}^\h u_\Phi(t')\|_{L^2}\|\D_k^\h \fb_\Phi(t')\|_{L^2}\,dt'\\
\lesssim & \sum_{|k'-k|\leq 4}2^{k'}\int_0^t\|\p_yu_\Phi(t')\|_{\cB^{\f12}}\|\D_{k'}^\h u_\Phi(t')\|_{L^2}
\|\D_k^\h \fb_\Phi(t')\|_{L^2}\,dt'\\
\lesssim & \sum_{|k'-k|\leq 4}2^{k'}\Bigl(\int_0^t\|\D_{k'}^\h u_\Phi(t')\|_{L^2}^2\dot{\tht}^3(t')\,dt'\Bigr)^{\frac12}
\Bigl(\int_0^t\|\D_k^\h \fb_\Phi(t')\|_{L^2}^2\dot{\tht}(t')\,dt'\Bigr)^{\frac12}\\
 \lesssim &d_k^22^{-2ks}\|u_\Phi\|_{\wt{L}_{t,\dot{\tht}^3(t)}(\cB^{s+\frac34})}\|\fb_\Phi\|_{\wt{L}_{t,\dot{\tht}(t)}(\cB^{s+\frac14})},
\end{align*}
which together with \eqref{S3eq8wr} and \eqref{S3eq8hk} ensures that \eqref{S4eq7wq} holds for any $s>-\f14$ in case
$\fa=u.$

Finally when $\fa=\p_yu,$ we get, by applying \eqref{S3eq8}, that for any $s\in\R$
\begin{align*}
\int_0^t\bigl|\bigl(&\D_k^{\rm h}(T^\h_{\p_x\p_y
u}u)_\Phi\ |\ \D_k^{\rm
h}\fb_\Phi\bigr)_{L^2}\bigr|\,dt'\\
\lesssim & \sum_{|k'-k|\leq 4}\int_0^t\|S_{k'-1}^\h \p_x\p_yu_\Phi(t')\|_{L^\infty_\h(L^2_{\rm v})}
\|\D_{k'}^\h u_\Phi(t')\|_{L^2_\h(L^\infty_{\rm v})}\|\D_k^\h \fb_\Phi(t')\|_{L^2}\,dt'\\
\lesssim & \sum_{|k'-k|\leq 4}2^{k'}\int_0^t\|\p_yu_\Phi(t')\|_{\cB^{\f12}}\|\D_{k'}^\h \p_yu_\Phi(t')\|_{L^2}
\|\D_k^\h \fb_\Phi(t')\|_{L^2}\,dt'\\
\lesssim & \sum_{|k'-k|\leq 4}2^{k'}\Bigl(\int_0^t\|\D_{k'}^\h \p_yu_\Phi(t')\|_{L^2}^2\dot{\tht}^3(t')\,dt'\Bigr)^{\frac12}
\Bigl(\int_0^t\|\D_k^\h \fb_\Phi(t')\|_{L^2}^2\dot{\tht}(t')\,dt'\Bigr)^{\frac12}\\
 \lesssim &d_k^22^{-2ks}\|\p_yu_\Phi\|_{\wt{L}_{t,\dot{\tht}^3(t)}(\cB^{s+\frac34})}\|\fb_\Phi\|_{\wt{L}_{t,\dot{\tht}(t)}(\cB^{s+\frac14})},
\end{align*}
which together with \eqref{S3eq8wr} and \eqref{S3eq8hk} ensures that \eqref{S4eq7wq} holds for any $s>-\f14$ in case
$\fa=\p_yu.$
This finishes the proof of Lemma \ref{lem3.1}.
\end{proof}

\begin{proof}[Proof of Lemma \ref{lem3.2}] Once again,  by applying Bony's decomposition \eqref{Bony}  for
 the horizontal variable to  $v\fa$,  we obtain
\beno
v\fa=T^\h_{v}\fa+T^\h_{\fa}v+R^h(v,\fa).
\eeno
Accordingly, we shall handle the following three terms:\smallskip

\no $\bullet$ \underline{Estimate of
$\int_0^t\bigl(\D_k^{\rm h}(T^\h_{v}\fa)_\Phi\ |\ \D_k^{\rm
h}\fb_\Phi\bigr)_{L^2}\,dt'$}

We first observe from \eqref{S4eq-1} and the assumption: $\|\fa_\Phi(t)\|_{\cB^{\f12}}\leq \de e^{-\fk t}$ for $t\leq T^\star_\fa.\,$ that
\beno
\begin{split}
\int_0^t\bigl|\bigl(\D_k^{\rm h}(T^\h_{v}\fa)_\Phi\ |&\ \D_k^{\rm
h}\fb_\Phi\bigr)_{L^2}\bigr|\,dt'\\
\lesssim & \sum_{|k'-k|\leq 4}\int_0^t\|S_{k'-1}^\h v_\Phi(t')\|_{L^\infty}
\|\D_{k'}^\h\fa_\Phi(t)\|_{L^2}\|\D_k^\h \fb_\Phi(t')\|_{L^2}\,dt'\\
\lesssim & \sum_{|k'-k|\leq 4} 2^{-\f{k'}2}\int_0^t\|S_{k'-1}^\h v_\Phi(t')\|_{L^\infty}\|\fa_\Phi(t')\|_{\cB^{\f12}}
\|\D_k^\h \fb_\Phi(t')\|_{L^2}\,dt'\\
\lesssim & \sum_{|k'-k|\leq 4} 2^{-\f{k'}2}\int_0^t\dot{\tht}^2(t')\|S_{k'-1}^\h v_\Phi(t')\|_{L^\infty}
\|\D_k^\h \fb_\Phi(t')\|_{L^2}\,dt'.
\end{split}
\eeno
Due to $\p_xu+\p_yv=0$ and \eqref{S1eq2},  we write $v(t,x,y)=-\int_0^y\p_x u(t,x,y')\,dy'.$ Then we deduce
from Lemma \ref{lem:Bern} that
\beq \label{S3eq11}
\begin{split}
\|\D_k^\h v_\Phi(t)\|_{L^\infty}\leq & \int_0^1\|\D_k^\h\p_xu_\Phi(t,\cdot,y')\|_{L^\infty_\h}\,dy'\\
\lesssim &2^{\f{3k}2}\int_0^1\|\D_k^\h u_\Phi(t,\cdot,y')\|_{L^2_h}\,dy'\lesssim 2^{\f{3k}2}\|\D_k^\h u_\Phi(t)\|_{L^2},
\end{split}
\eeq
from which and $s<\f34$, we infer
\beq \label{S3eq12}
\begin{split}
\Bigl(\int_0^t\|S_{k'-1}^\h v_\Phi(t')\|_{L^\infty}^2\dot{\tht}^3(t')\,dt'\Bigr)^{\f12}
\leq &\sum_{\ell\leq k'-2}2^{\f{3\ell}2}\Bigl(\int_0^t\|\D_\ell^\h u_\Phi(t)\|_{L^2}^2\dot{\tht}^3(t')\,dt'\Bigr)^{\f12}\\
\lesssim &\sum_{\ell\leq k'-2}d_\ell 2^{\ell\left(\f34-s\right)}\|u_\Phi\|_{\wt{L}^2_{t,\dot{\tht}^3(t)}(\cB^{s+\f34})}\\
\lesssim &d_{k'}2^{{k'}\left(\f34-s\right)}\|u_\Phi\|_{\wt{L}^2_{t,\dot{\tht}^3(t)}(\cB^{s+\frac34})}.
\end{split}
\eeq
Consequently, by virtue of Definition \ref{def1.1}, we obtain
\beno
\begin{split}
\int_0^t&\bigl|\bigl(\D_k^{\rm h}(T^\h_{v}\fa)_\Phi\ |\ \D_k^{\rm
h}\fb_\Phi\bigr)_{L^2}\bigr|\,dt'\\
\lesssim & \sum_{|k'-k|\leq 4}2^{-\f{k'}2}\Bigl(\int_0^t\|S_{k'-1}^\h v_\Phi(t')\|_{L^\infty}^2\dot{\tht}^3(t')\,dt'\Bigr)^{\f12}\Bigl(\int_0^t\|\D_k^\h \fb_\Phi(t')\|_{L^2}^2\dot{\tht}(t')\,dt'\Bigr)^{\f12}\\
\lesssim &d_k^22^{-2ks}\|u_\Phi\|_{\wt{L}^2_{t,\dot{\tht}^3(t)}(\cB^{s+\f34})}\|\fb_\Phi\|_{\wt{L}^2_{t,\dot{\tht}^3(t)}(\cB^{s+\f14})}.
\end{split}
\eeno
The case for $s=\f34$ can be proved along the same line to \eqref{S3eq8au}.

\no $\bullet$ \underline{Estimate of
$\int_0^t\bigl|\bigl(\D_k^{\rm h}(T^\h_{\fa}v)_\Phi\ |\ \D_k^{\rm
h}\fb_\Phi\bigr)_{L^2}\bigr|\,dt'$}

Notice that
\beno
\begin{split}
\int_0^t\bigl|\bigl(&\D_k^{\rm h}(T^\h_{\fa}v)_\Phi\ |\ \D_k^{\rm h}\fb_\Phi\bigr)_{L^2}\bigr|\,dt'\\
\lesssim & \sum_{|k'-k|\leq 4}\int_0^t\|S_{k'-1}^\h \fa_\Phi(t')\|_{L^\infty_\h(L^2_{\rm v})}
\|\D_{k'}^\h v_\Phi(t)\|_{L^2_\h(L^\infty_{\rm v})}\|\D_k^\h \fb_\Phi(t')\|_{L^2}\,dt',
\end{split}
\eeno
which together with \eqref{S4eq-1} and \eqref{S3eq11} ensures that for $t\leq T^\star_\fa$
\beno
\begin{split}
\int_0^t\bigl|\bigl(&\D_k^{\rm h}(T^\h_{\fa}v)_\Phi\ |\ \D_k^{\rm
h}\fb_\Phi\bigr)_{L^2}\bigr|\,dt'\\
\lesssim & \sum_{|k'-k|\leq 4}2^{k'}\int_0^t\|\fa_\Phi(t')\|_{\cB^{\f12}}\|\D_{k'}^\h u_\Phi(t')\|_{L^2}
\|\D_k^\h \fb_\Phi(t')\|_{L^2}\,dt'\\
\lesssim & \sum_{|k'-k|\leq 4}2^{k'}\Bigl(\int_0^t\|\D^\h_{k'}u_\Phi(t')\|_{L^2}^2\dot{\tht}^3(t')\,dt'\Bigr)^{\frac12}
\Bigl(\int_0^t\|\D_k^\h \fb_\Phi(t')\|_{L^2}^2\dot{\tht}(t')\,dt'\Bigr)^{\frac12}.
\end{split}
\eeno
Then thanks to Definition \ref{def1.1}, for any $s\in\R,$ we achieve
\beq \label{S3eq12op}
\begin{split}
\int_0^t\bigl|\bigl(\D_k^{\rm h}(T^\h_{\fa}v)_\Phi\ |\ \D_k^{\rm
h}\fb_\Phi\bigr)_{L^2}\bigr|\,dt' \lesssim & d_k^22^{-2ks}\|u_\Phi\|_{\wt{L}^2_{t,\dot{\tht}^3(t)}(\cB^{s+\f34})}\|\fb_\Phi\|_{\wt{L}^2_{t,\dot{\tht}(t)}(\cB^{s+\f14})}.
\end{split}
\eeq

\no $\bullet$ \underline{Estimate of
$\int_0^t\bigl|\bigl(\D_k^{\rm h}(R^\h(v,\fa))_\Phi\ |\ \D_k^{\rm
h}\fb_\Phi\bigr)_{L^2}\bigr|\,dt'$}\vspace{0.2cm}

We get, by applying lemma
\ref{lem:Bern} and \eqref{S3eq11}, that
\beno
\begin{split}
\int_0^t\bigl|&\bigl(\D_k^{\rm h}(R^\h(v,\fa))_\Phi\ |\ \D_k^{\rm
h}\fb_\Phi\bigr)_{L^2}\bigr|\,dt'\\
\lesssim &2^{\f{k}2}\sum_{k'\geq k-3}\int_0^t\|{\D}_{k'}^\h v_\Phi(t')\|_{L^2_\h(L^\infty_{\rm v})}\|\wt{\D}_{k'}^\h \fa_\Phi(t')\|_{L^2}\|\D_k^\h \fb_\Phi(t')\|_{L^2}\,dt'\\
\lesssim & 2^{\f{k}2}\sum_{k'\geq k-3}2^{\f{k'}2}\int_0^t\|{\D}_{k'}^\h u_\Phi(t')\|_{L^2}\|\fa_\Phi(t')\|_{\cB^{\f12}}\|\D_k^\h \fb_\Phi(t')\|_{L^2}\,dt'\\
\lesssim & 2^{\f{k}2}\sum_{k'\geq k-3}2^{\f{k'}2}\Bigl(\int_0^t\|{\D}_{k'}^\h  u_\Phi(t')\|_{L^2}^2\dot{\tht}^3(t')\,dt'\Bigr)^{\f12}
 \Bigl(\int_0^t\|\D_k^\h \fb_\Phi(t')\|_{L^2}^2\dot{\tht}(t')\,dt'\Bigr)^{\f12},
\end{split}
\eeno
which together with Definition \ref{def1.1} and $s>-\f14$ ensures that
\beq\label{S3eq12we}
\begin{split}
\int_0^t&\bigl|\bigl(\D_k^{\rm h}(R^\h(v,\fa)_\Phi\ |\ \D_k^{\rm
h}\fb_\Phi\bigr)_{L^2}\bigr|\,dt'\\
\lesssim &d_k2^{-{2ks}}\|u_\Phi\|_{\wt{L}^2_{t,\dot{\tht}^3(t)}(\cB^{s+\f34})}\|\fb_\Phi\|_{\wt{L}^2_{t,\dot{\tht}(t)}(\cB^{s+\f14})} \Bigl(\sum_{k'\geq k-3}d_{k'}2^{(k-k')\left(\f14+s\right)}\Bigr)\\
\lesssim &d_k^22^{-2ks}\|u_\Phi\|_{\wt{L}^2_{t,\dot{\tht}^3(t)}(\cB^{s+\f34})}\|\fb_\Phi\|_{\wt{L}^2_{t,\dot{\tht}(t)}(\cB^{s+\f14})} .
\end{split}
\eeq

By summing up the above estimates, we  arrive at \eqref{S4eq9g} for $s\in\left]-\f14,\f34\right].$

On the other hand, we deduce from \eqref{S3eq11} that for $ t\leq T^\star_{1,\fa}$
\beno
\|S_{k'-1}^\h v_\Phi(t')\|_{L^\infty}\lesssim 2^{\f{k'}2}\|u_\Phi(t')\|_{\cB^1}\lesssim 2^{\f{k'}2}\dot\tht(t'),
\eeno
from which, we deduce that
\begin{align*}
\int_0^t\bigl|\bigl(\D_k^{\rm h}(T^\h_{v}\fa)_\Phi\ |&\ \D_k^{\rm
h}\fb_\Phi\bigr)_{L^2}\bigr|\,dt'\\
\lesssim & \sum_{|k'-k|\leq 4}\int_0^t\|S_{k'-1}^\h v_\Phi(t')\|_{L^\infty}
\|\D_{k'}^\h\fa_\Phi(t')\|_{L^2}\|\D_k^\h \fb_\Phi(t')\|_{L^2}\,dt'\\
\lesssim & \sum_{|k'-k|\leq 4} 2^{\f{k'}2}\Bigl(\int_0^t\dot\tht(t')\|\D_{k'}^\h\fa_\Phi(t')\|_{L^2}^2\,dt'\Bigr)^{\f12}
\Bigl(\int_0^t\dot\tht(t')\|\D_k^\h \fb_\Phi(t')\|_{L^2}\,dt'\Bigr)^{\f12}\\
\lesssim & d_k^22^{-2ks}\|\fa_\Phi\|_{\wt{L}^2_{t,\dot{\tht}(t)}(\cB^{s+\frac14})}
\| \fb_\Phi\|_{\wt{L}^2_{t,\dot{\tht}(t)}(\cB^{s+\frac14})},
\end{align*}
which together with \eqref{S3eq12op} and \eqref{S3eq12we} ensures \eqref{S4eq9gq}. This completes the proof of Lemma \ref{lem3.2}.
\end{proof}

\setcounter{equation}{0}


\section{Propagation of regularities for $\p_y u_\Phi$}\label{Sect3}

The purpose of this section is to present the proof of part (3) of Theorem \ref{th1.2}.
Due to the boundary conditions for $v:$ $v|_{y=0}=0=v|_{y=1},$ we get, by integrating $\pa_xu+\pa_yv=0$ over $[0,1]$
with respect to $y$ variable,
that
\beno
\p_x\int_0^1u(t,x,y)\,dy=0 \Rightarrow  \int_0^1u(t,x,y)\,dy\eqdefa f(t).
\eeno
Let us determine the pressure $p$ via
\beq\label{S7eq2}
\pa_xp=\pa_yu(t,x,1)-\pa_yu(t,x,0)-\f12\pa_x\int_0^1u^2(t,x,y)\,dy.
\eeq
So that by integrating the equation $\p_t^2u+\pa_tu+u\pa_xu+v\pa_yu-\pa_y^2u+\pa_xp=0$ for $y\in [0,1],$ we obtain
\beno
f''(t)+f(t)=0,
\eeno
which together with the compatibility conditions that $\int_0^1u_0(x,y)\,dy=0=\int_0^1u_1(x,y)\,dy,$ that is,
$f(0)=f'(0)=0,$ ensures that $f(t)=0.$

\begin{proof}[Proof of \eqref{S1eq14}]
In what follows, we shall always denote ${\om}\eqdefa \p_yu$ and assume that $t<T^\star_1$ with
$T^\star_1$ being determined by
\beq\label{S7eq1} T^\star_1\eqdefa
\sup\bigl\{\ t>0,\ \  \|\p_y\om_\Phi(t)\|_{\cB^{\f12}}\leq \de e^{-\frak{K}t} \  \bigr\}. \eeq

Thanks to $\p_xu+\p_yv=0,$ we get, by applying $\p_y$ to \eqref{S1eq1}, that
\beno
\p_t^2 {\om}+\p_t{\om}+u\p_x{\om}+v\p_y\om-\p_y^2{\om}+\p_y\p_x p =0.
\eeno
In view of \eqref{S4eq1}, we get, by applying the Fourier multiplier $e^{\Phi(t,D)}$ to the above equation, that
\beno
\p_t(\p_t {\om})_\Phi+\lambda\dot\tht|D_x|^{\f 12}(\p_t{\om})_\Phi+(\p_t{\om})_\Phi+(u\p_x{\om})_\Phi+(v\p_y {\om})_\Phi-\p_y^2{\om}_\Phi+\p_y\p_x p_\Phi =0.
\eeno
By taking $L^2$ inner product of the above equation  by $(\p_t{\om})_\Phi+\f 12 {\om}_\Phi$ and using the fact that $-\p_y {\om} +\p_x p$ vanishes
 on the boundary of $\cS,$ we get, by using a similar derivation of \eqref{S4eq6rt}, that
\beq \label{S7eq3}
\begin{split}
\frak{G}_k({\om})(t)
\leq & \f34\bigl\|e^{a|D_x|^{\f12}}\D_k^{\rm h}(\p_t{\om})(0)\bigr\|_{L^2}^2+\f12\bigl\|e^{a|D_x|^{\f12}}\D_k^{\rm h}\p_y{\om}_0\bigr\|_{L^2}^2
+\f1{2}\bigl\|e^{a|D_x|^{\f12}}\D_k^{\rm h}\om_0\bigr\|_{L^2}^2\\
&+\f\lam2\de^{\f12}2^{\f{k}2}\bigl\|e^{a|D_x|^{\f12}}\D_k^{\rm h}\om_0\bigr\|_{L^2}^2+\lam^2\de2^k\bigl\|e^{a|D_x|^{\f12}}\D_k^{\rm h}\om_0\bigr\|_{L^2}^2\\
&+\int_0^t\bigl|\bigl(\ektp\D_k^\h\left(u\p_x\om\right)_\Phi | \ektp\D_k^\h \bigl(\p_t\om+\f12\om\bigr)_\Phi\bigr)_{L^2}\bigr|\,dt'\\
&+\int_0^t\bigl|\bigl(\ektp\D_k^\h\left(v\p_y \om\right)_\Phi | \ektp\D_k^\h \bigl(\p_t\om+\f12\om\bigr)_\Phi\bigr)_{L^2}\bigr|\,dt'+I_k(t),
\end{split}
\eeq
where the functional $\frak{G}_k(\om)(t)$ is determined by \eqref{S4eq6a}
and
\beq \label{S7eq4}
I_k(t)\eqdefa \int_0^t\bigl(\ekt\D_k^h \p_xp_\Phi | \ekt\D_k^h\bigl(\p_t\p_y\om_\Phi+\f12\p_y \om_\Phi\bigr)\bigr)_{L^2}\,dt'. \eeq

It follows from Lemma \ref{lem3.1} that for any $\frak{s}>-\f14$
\begin{align*}
\int_0^t\bigl|\bigl(\ektp\D_k^\h\left(u\p_x\om\right)_\Phi &| \ektp\D_k^\h \bigl(\p_t\om+\f12\om)_\Phi\bigr)_{L^2}\bigr|\,dt'
\lesssim d_k^2
2^{-2k\fs}\|\ektp \om_\Phi\|_{\wt{L}^2_{t,\dot{\tht}^3(t)}(\cB^{\fs+\frac34})}\\
&\qquad\times
\bigl(\|\ektp(\om_\Phi, \p_t\om_\Phi)\|_{\wt{L}^2_{t,\dot{\tht}(t)}(\cB^{\fs+\frac14})}
+\lam \|\ektp \om_\Phi\|_{\wt{L}^2_{t,\dot{\tht}^3(t)}(\cB^{\fs+\frac34})}\bigr).
\end{align*}

Whereas it follows from \eqref{S4eq19ui}, Lemma \ref{lem3.2} and \eqref{S7eq1} that for  $t\leq T^\star_1$ and for any $\fs>-\f14$
\begin{align*}
\int_0^t\bigl|\bigl(&\ektp\D_k^\h\left(v\p_y\om\right)_\Phi | \ektp\D_k^\h \bigl(\p_t\om+\f12\om\bigr)_\Phi\bigr)_{L^2}\bigr|\,dt'
\lesssim d_k^2
2^{-2k\fs}\bigl(\|\ektp u_\Phi\|_{\wt{L}^2_{t,\dot{\tht}^3(t)}(\cB^{\fs+\frac34})}\\
&+\|\ektp \p_y\om_\Phi\|_{\wt{L}^2_{t,\dot{\tht}(t)}(\cB^{\fs+\frac14})}\bigr)
\bigl(\|\ektp (\om_\Phi,\p_t\om_\Phi)\|_{\wt{L}^2_{t,\dot{\tht}(t)}(\cB^{\fs+\frac14})}
+\lam \|\ektp \om_\Phi\|_{\wt{L}^2_{t,\dot{\tht}^3(t)}(\cB^{\fs+\frac34})}\bigr).
\end{align*}

We claim that
\beq \label{S7eq5}
\begin{split}
|I_k(t)|\leq  C&d_k^22^{-2k\fs}\Bigl(\frak{N}_\fs+\|e^{\f{\fk}2t'}u_\Phi\|_{\wt{L}^\infty_t(\cB^{\fs+\f34})}^{2}
+\|\ektp (\p_tu)_\Phi\|_{\wt{L}^2_{t,\dot{\tht}(t)}(\cB^{\fs+\f12})}^2\\
&+\|\ektp \om_\Phi\|_{\wt{L}^\infty_t(\cB^{\fs})}\|\ektp \p_y\om_\Phi\|_{\wt{L}^\infty_t(\cB^{\fs})}
+\|\ektp \om_\Phi\|_{\wt{L}^2_t(\cB^{\fs})}\|\ektp \p_y\om_\Phi\|_{\wt{L}^2_t(\cB^{\fs})}
\\
&+\f{a\fk}{4C}
\|e^{\f{\fk}2t'}\om_\Phi\|_{\wt{L}^\infty_t(\cB^{\fs+\f12})}^2+\|\ektp (\om,\p_y\om)_\Phi\|_{\wt{L}^2_{t,\dot{\tht}(t)}(\cB^{\fs+\f14})}^2
\\
&+\lam\bigl(\|\ektp (u,\om)_\Phi\|_{\wt{L}^2_{t,\dot{\tht}^3(t)}(\cB^{\fs+\f34})}^2+\|\ektp u_\Phi\|_{\wt{L}^2_{t,\dot{\tht}^3(t)}(\cB^{\fs+1})}^2\bigr)
\Bigr),
\end{split}
\eeq
where
\beq \label{S7eq5g}\begin{split}
\frak{N}_\fs\eqdefa &\bigl\|e^{\de |D_x|^{\f12}} \om_0\bigr\|_{\cB^{\fs}}\bigl\|e^{\de |D_x|^{\f12}} \p_y\om_0\bigr\|_{\cB^{\fs}}\\
&+\bigl\|e^{\de|D_x|^{\f12}}u_0\bigr\|_{\cB^{\f12}}^{\f12}\bigl\|e^{\de|D_x|^{\f12}}\om_0\bigr\|_{\cB^{\f12}}^{\f12}
\bigl\|e^{\de|D_x|^{\f12}}u_0\bigr\|_{\cB^{1+\fs}}
\bigl\|e^{\de|D_x|^{\f12}}\om_0\bigr\|_{\cB^{\fs}}^{\f12}\bigl\|e^{\de|D_x|^{\f12}}\p_y\om_0\bigr\|_{\cB^{\fs}}^{\f12}.
\end{split}
\eeq
The proof of \eqref{S7eq5} involves technicality, which we postpone in the Appendix \ref{sectb}.

Let us denote
\beq \label{S7eq10}
\begin{split}
\frak{L}_s\eqdefa & \bigl\|e^{\de |D_x|^{\f12}} (\om_0,\p_yu_1,\p_y\om_0)\bigr\|_{\cB^{\fs}}^2+a\fk\bigl\|e^{\de|D_x|^{\f12}}\om_0\bigr\|_{\cB^{\fs+\f14}}^{2}
\\
&+a^2\fk^2\bigl\|e^{\de|D_x|^{\f12}}\om_0\bigr\|_{\cB^{\fs+\f12}}^{2}
+\bigl\|e^{\de|D_x|^{\f12}}u_0\bigr\|_{\cB^{\f12}}\bigl\|e^{\de|D_x|^{\f12}}\om_0\bigr\|_{\cB^{\f12}}\bigl\|e^{\de|D_x|^{\f12}}u_0\bigr\|_{\cB^{\fs+1}}^2.
\end{split}
\eeq
Then
by
substituting the above estimates into \eqref{S7eq3}, we find
\beno
\begin{split}
\frak{G}_k(\om)(t)
&\leq  Cd_k^2
2^{-2k\fs}\Bigl(\frak{L}_\fs+\|e^{\f{\fk}2t'}u_\Phi\|_{\wt{L}^\infty_t(\cB^{\fs+\f34})}^{2}
+\|\ektp (\p_tu)_\Phi\|_{\wt{L}^2_{t,\dot{\tht}(t)}(\cB^{\fs+\f12})}^2\\
&+\|\ektp \om_\Phi\|_{\wt{L}^\infty_t(\cB^{\fs})}\|\ektp \p_y\om_\Phi\|_{\wt{L}^\infty_t(\cB^{\fs})}
+\|\ektp \om_\Phi\|_{\wt{L}^2_t(\cB^{\fs})}\|\ektp \p_y\om_\Phi\|_{\wt{L}^2_t(\cB^{\fs})}
\\
&+\f{a\fk}{4C}
\|e^{\f{\fk}2t'}\om_\Phi\|_{\wt{L}^\infty_t(\cB^{\fs+\f12})}^2+\bigl\|\ektp\bigl(\om,\p_y\om_\Phi,\p_t\om_\Phi\bigr)\bigr\|_{\wt{L}^2_{t,\dot{\tht}(t)}(\cB^{\fs+\f14})}^2
\\
&
+\lam\bigl(\|\ektp u_\Phi\|_{\wt{L}^2_{t,\dot{\tht}^3(t)}(\cB^{\fs+\f34})}^2
+\|\ektp u_\Phi\|_{\wt{L}^2_{t,\dot{\tht}^3(t)}(\cB^{\fs+1})}^2\bigr)
+\lam^2\|\ektp \om_\Phi\|_{\wt{L}^2_{t,\dot{\tht}^3(t)}(\cB^{\fs+\f34})}^2
\Bigr).
\end{split}
\eeno
Then thanks to \eqref{S4eq-1} and \eqref{Esu}, by multiplying the above inequality by $2^{2k\fs}$ and then
taking square root of the resulting inequality, finally  by summing up the inequalities for $k$ over $\Z,$ we achieve
\begin{align*}
{\rm E}_{\fs,\lam}(\om)(t)
 \leq &\frac12\bigl(\|\ektp \p_y\om_\Phi\|_{\wt{L}^\infty_t(\cB^{\fs})}+a\fk\|e^{\f12\fk t'} \om_\Phi\|_{\wt{L}^\infty_{t}(\cB^{\fs+\f12})}
 +\|\ektp\p_y \om_\Phi\|_{\wt{L}^2_t(\cB^{\fs})}\bigr)\\
 &+  C\Bigl(\frak{L}_\fs
+\|e^{\f{\fk}2t'}u_\Phi\|_{\wt{L}^\infty_t(\cB^{\fs+\f34})}
+\|\ektp (\p_tu)_\Phi\|_{\wt{L}^2_{t,\dot{\tht}(t)}(\cB^{\fs+\f12})}\\
&+\|\ektp (w,\p_yw,\p_t\om)_\Phi\|_{\wt{L}^2_{t,\dot{\tht}(t)}(\cB^{\fs+\f14})}^2
+\sqrt{\lam}\bigl(\|\ektp u_\Phi\|_{\wt{L}^2_{t,\dot{\tht}^3(t)}(\cB^{\fs+\f34})}\\
&+\|\ektp u_\Phi\|_{\wt{L}^2_{t,\dot{\tht}^3(t)}(\cB^{\fs+1})}\bigr)+\lam\|\ektp \om_\Phi\|_{\wt{L}^2_{t,\dot{\tht}^3(t)}(\cB^{\fs+\f34})}
\Bigr),
\end{align*}
where we used \eqref{S4eq-1} so that $ \lam\de^{\f12}=\f{a\fk}4.$ By taking $\lam=4C^2,$ we obtain
\begin{align*}
{\rm E}_{\fs,\f\lam4}(\om)(t)
 \leq
  C\Bigl(&\frak{L}_\fs+
\|e^{\f{\fk}2t'}u_\Phi\|_{\wt{L}^\infty_t(\cB^{\fs+\f34})}
+\|\ektp \p_tu_\Phi\|_{\wt{L}^2_{t,\dot{\tht}(t)}(\cB^{\fs+\f12})}\\
&+\sqrt{\lam}\|\ektp u_\Phi\|_{\wt{L}^2_{t,\dot{\tht}^3(t)}(\cB^{\fs+\f34})}+\sqrt{\lam}\|\ektp u_\Phi\|_{\wt{L}^2_{t,\dot{\tht}^3(t)}(\cB^{\fs+1})}
\Bigr),
\end{align*}
which together with \eqref{Esuu} ensures that for $t\leq T^\star_1,$
\beq \label{S7eq11}
{\rm E}_{\fs,\f\lam4}(\om)(t)
 \leq
  C\bigl(\frak{L}_\fs+\frak{H}_\fs(u_0,u_1)+\frak{H}_{\fs+\f14}(u_0,u_1)\bigr).
\eeq
In particular, we deduce from \eqref{S1eq13} and \eqref{S7eq11} that
\begin{align*}
\|\ekt\p_y\om_\Phi(t)\|_{\cB^{\f12}}\leq {\rm E}_{\f12,\f\lam4}(\om)(t)\leq C\bigl(\frak{L}_{\f12}+\frak{H}_{\f12}(u_0,u_1)+\frak{H}_{\f34}(u_0,u_1)\bigr)\leq Cc_1,
\end{align*}
so that for $t\leq T^\star_1,$ there holds
\beno
\|\p_y\om_\Phi(t)\|_{\cB^{\f12}}\leq Cc_1e^{-\fk t}\leq \f\de2e^{-\fk t},
\eeno
for $c_1$ sufficiently small. Then a continuous argument shows that $T^\star_1$ defined by
\eqref{S7eq1} equals $\infty.$ Moreover, \eqref{S1eq14} holds for  ${\rm E}_{\fs}(\p_y u)(t).$

Let us turn to the estimate of $\|\ektp(\p_t^2u)_\Phi\|_{\wt{L}_t^2(\cB^{\fs})}.$
  Indeed,
by  applying $\D_k^\h$ to  \eqref{S4eq2} and taking $L^2$ inner product of the resulting equation
with $\D_k^\h\left(\p^2_tu\right)_\Phi,$  we obtain
\beno
\|\D_k^\h\left(\p^2_tu\right)_\Phi\|_{L^2}\leq \|\D_k^\h\pa_y^2 u_\Phi \|_{L^2}+\|\D_k^\h(\p_t u)_\Phi  \|_{L^2}
+\|\D_k^\h(u\pa_x u)_\Phi  \|_{L^2}+\|\D_k^\h(v\pa_yu)_\Phi  \|_{L^2}. \eeno
By multiplying the above inequality by $2^{2k\fs}e^{2\frak{K}t}$ and then integrating the resulting inequality over $[0,t]$,
and finally summing up $k$ over $\Z,$ we achieve
\beq\label{S7eq30}
\begin{split}
\|\ekt\left(\p^2_tu\right)_\Phi\|_{\wt{L}_t^2(\cB^{\fs})}\leq \bigl\|e^{\frak{K}t}\bigl(\p_t u,\pa_y^2 u\bigr)_\Phi\|_{\wt{L}^2_t(\cB^{\fs})}
+\bigl\|e^{\frak{K}t}\bigl(u\pa_x u,v\pa_yu\bigr)_\Phi  \bigr\|_{\wt{L}^2_t(\cB^{\fs})}.
\end{split}
\eeq
Yet it follows from the proof of Lemmas \ref{lem3.1} and \ref{lem3.2} that
\begin{align*}
\bigl\|e^{\frak{K}t}\bigl(u\pa_x u,v\pa_yu\bigr)_\Phi  \bigr\|_{\wt{L}^2_t(\cB^{\fs})}
\lesssim \bigl\|\ektp u_\Phi\bigr\|_{\wt{L}^2_{t,\dot\tht^3(t)}(\cB^{\fs+1})}.
\end{align*}
By substituting the above estimates, \eqref{Esuuv} and \eqref{S7eq11} into  \eqref{S7eq30}, we obtain the estimate \eqref{S1eq14} for $\|\ekt\left(\p^2_tu\right)_\Phi\|_{\wt{L}_t^2(\cB^{\fs})}.$
This completes the proof of \eqref{S1eq14} and thus Theorem \ref{th1.2}.
\end{proof}

\renewcommand{\theequation}{\thesection.\arabic{equation}}
\setcounter{equation}{0}
\section{Global well-posedness of the  system \eqref{S1eq1.8}}\label{Sect4}

In this section, we prove the global well-posedness of the scaled anisotropic hyperbolic
 Navier-Stokes system \eqref{S1eq1.8} with small Gevrey class 2 data and establish uniform estimates
 for such solutions, namely, we are going to present the proof of Theorem \ref{thm1.1}.

\begin{proof}[Proof of Theorem \ref{thm1.1}]

In the rest of this section, we shall prove that under the assumption of \eqref{S1eq8}, there holds the {\it a priori} estimate
 \eqref{S3eq18} for smooth enough solutions of \eqref{S1eq1.8}, and neglect the regularization procedure. For simplicity,
 we shall neglect the supscript $\e.$ Then in view of \eqref{S1eq1.8} and \eqref{S4eq1}, we observe that
$(u_\Phi, v_\Phi)$ verifies
\beq\label{S3eq4}
\left\{
\begin{array}{ll}
\p_t(\p_tu)_{\Phi}+\lam\dot{\tht}|D_x|^{\f12}(\p_tu)_{\Phi}+(\p_tu)_{\Phi}+\left(u\p_xu\right)_{\Phi}+\left(v\p_y u\right)_{\Phi}\\
\qquad-\e^2\p_x^2 u_{\Phi}
-\p_y^2u_{\Phi}+\p_x p_{\Phi}=0, \\
\e^2\left(\p_t(\p_tv)_{\Phi}+\lam\dot{\tht}|D_x|^{\f12}(\p_tv)_{\Phi}+(\p_tv)_{\Phi}+\left(u\p_xv\right)_{\Phi}+\left(v\p_y v\right)_{\Phi}\right.\\
\qquad\left.-\e^2\p_x^2 v_{\Phi}
-\p_y^2v_{\Phi}\right)+\p_y p_{\Phi}=0,\\
\p_x u_{\Phi}+\p_yv_{\Phi} =0\quad\mbox{for}  \quad (t,x,y)\in\R_+\times\cS,\\
\left(u_{\Phi}, v_{\Phi}\right)|_{y=0}=\left(u_\Phi, v_\Phi\right)|_{y=1}=0.
\end{array}
\right.
 \eeq
By applying the dyadic operator
$\D_k^{\rm h}$ to \eqref{S3eq4} and then taking the $L^2$ inner
product of the resulting equations with $\left(\D_k^{\rm
h}(\p_tu)_{\Phi}, \D_k^\h (\p_tv)_{\Phi}\right),$  we find \beq \label{S3eq5}
\begin{split}
\f12\f{d}{dt}&\bigl\|\D_k^{\rm h}(\p_tu, \e \p_tv)_{\Phi}(t)\bigr\|_{L^2}^2+\lam\dot{\tht}\bigl\||D_x|^{\f14}\D_k^{\rm
h}(\p_tu,\e \p_tv)_{\Phi}\bigr\|_{L^2}^2+\bigl\|\D_k^{\rm
h}(\p_tu,\e \p_tv)_{\Phi}\bigr\|_{L^2}^2\\
&-\e^2\bigl(\D_k^\h \p_x^2(u, \e v)_{\Phi} |  \D_k^\h (\p_tu, \e \p_tv)_{\Phi} \bigr)_{L^2}-\bigl(\D_k^\h \p_y^2(u, \e v)_{\Phi} |  \D_k^\h (\p_tu, \e \p_tv)_{\Phi} \bigr)_{L^2}\\
=&-\bigl(\D_k^\h\left(u\p_xu\right)_{\Phi} | \D_k^\h (\p_tu)_{\Phi}\bigr)_{L^2}-\bigl(\D_k^\h\left(v\p_yu\right)_{\Phi} | \D_k^\h (\p_tu)_{\Phi}\bigr)_{L^2}\\
&-\e^2\bigl(\D_k^\h\left(u\p_xv\right)_{\Phi} | \D_k^\h (\p_tv)_{\Phi}\bigr)_{L^2}-\e^2\bigl(\D_k^\h\left(v\p_yv\right)_{\Phi} | \D_k^\h (\p_tv)_{\Phi}\bigr)_{L^2},
\end{split}
\eeq
where we used the fact that $\p_x u_{\Phi}+\p_y v_{\Phi} =0,$  so that
\beno
\bigl(\D_k^{\rm h}\na p_{\Phi}\ |\ \D_k^{\rm
h}(\p_tu, \p_tv)_{\Phi}\bigr)_{L^2}=0.\eeno
Here and in all that follows, we always denote $$\left(\vec{\fa}, \vec{\fb}\right)_{L^2}\eqdefa \int_{\cS}\left(\fa_1(x)\fb_1(x)+\fa_2(x)\fb_2(x)\right)\,dx
$$ for $\vec{\fa}(x)=(\fa_1(x), \fa_2(x))$ and  $\vec{\fb}(x)=(\fb_1(x), \fb_2(x)).$

We first get, by a similar derivations as those equalities  above \eqref{S4eq5a}, that
\begin{align*}
\lam\dot{\tht}\||D_x|^{\f14}\D_k^{\rm
h}(\p_tu,\e\p_tv)_\Phi\|^2_{L^2}=&\lam\dot{\tht}\||D_x|^{\f14}\D_k^{\rm
h}(\p_tu_\Phi,\e\p_tv_\Phi)\|^2_{L^2}+\lam^3\dot{\tht}^3\||D_x|^{\f34}\D_k^{\rm
h}(u,\e v)_\Phi\|^2_{L^2}\\
&+\lam^2\dot{\tht}^2\f{d}{dt}\||D_x|^{\f12}\D_k^{\rm h}(u,\e v)_\Phi(t)\|_{L^2}^2,
\end{align*}
and
\begin{align*}
-\bigl(\D_k^\h\p_x^2(u,\e^2v)_\Phi | \D_k^{\rm h}(\p_tu,\p_tv)_\Phi\bigr)_{L^2}
=&\f{1}2\f{d}{dt}\|\D_k^{\rm h}\p_x(u,\e v)_\Phi(t)\|_{L^2}^2+\lam\dot{\tht}\||D_x|^{\f14}\D_k^{\rm
h} \p_x(u,\e v)_\Phi\|^2_{L^2},\\
-\bigl(\D_k^\h\p_y^2(u,\e^2v)_\Phi | \D_k^{\rm h}(\p_tu,\p_tv)_\Phi\bigr)_{L^2}
=&\f12\f{d}{dt}\|\D_k^{\rm h}\p_y(u,\e v)_\Phi(t)\|_{L^2}^2+\lam\dot{\tht}\||D_x|^{\f14}\D_k^{\rm
h} \p_y(u,\e v)_\Phi\|^2_{L^2}.
\end{align*}
By inserting the above equalities into \eqref{S3eq5}, we obtain
\beq\label{S3eq1}
\begin{split}
\f12&\f{d}{dt}\Bigl(\|\D_k^{\rm h}(\p_tu,\e\p_tv)_\Phi(t)\|_{L^2}^2+2\lam^2\dot{\tht}^2(t)\||D_x|^{\f12}\D_k^{\rm h}(u,\e v)_\Phi(t)\|_{L^2}^2\\
&+{\e^2}\|\D_k^{\rm h}\p_x(u,\e v)_\Phi(t)\|_{L^2}^2+\|\D_k^{\rm h}\p_y(u,\e v)_\Phi(t)\|_{L^2}^2\Bigr)\\
&+\|\D_k^{\rm h}(\p_tu,\e\p_tv)_\Phi\|_{L^2}^2+\lam\dot{\tht}\Bigl(\||D_x|^{\f14}\D_k^{\rm
h}\p_t(u,\e v)_\Phi\|^2_{L^2}\\
&+\e^2\||D_x|^{\f14}\D_k^{\rm
h} \p_x(u,\e v)_\Phi\|^2_{L^2}+\||D_x|^{\f14}\D_k^{\rm
h} \p_y(u,\e v)_\Phi\|^2_{L^2}\Bigr)\\
&-2\lam^2\dot{\tht}\ddot{\tht}\||D_x|^{\f12}\D_k^{\rm h}(u,\e v)_\Phi\|_{L^2}^2+\lam^3\dot{\tht}^3\||D_x|^{\f34}\D_k^{\rm
h}(u,\e v)_\Phi\|^2_{L^2}\\
=&-\bigl(\D_k^\h\left(u\p_xu\right)_{\Phi} | \D_k^\h (\p_tu)_{\Phi}\bigr)_{L^2}-\bigl(\D_k^\h\left(v\p_yu\right)_{\Phi} | \D_k^\h (\p_tu)_{\Phi}\bigr)_{L^2}\\
&-\e^2\bigl(\D_k^\h\left(u\p_xv\right)_{\Phi} | \D_k^\h (\p_tv)_{\Phi}\bigr)_{L^2}-\e^2\bigl(\D_k^\h\left(v\p_yv\right)_{\Phi} | \D_k^\h (\p_tv)_{\Phi}\bigr)_{L^2}.
\end{split}
\eeq

On the other hand, we get,
by taking $L^2$ inner product of \eqref{S3eq4} with $\D_k^\h (u,v)_\Phi,$ that
\beq \label{S3eq2}
\begin{split}
\bigl(&\D_k^{\rm h}\p_t(\p_tu,\e^2\p_tv)_\Phi | \D_k^\h (u,v)_\Phi\bigr)_{L^2}+\lam\dot{\tht}\bigl(|D_x|^{\f12}\D_k^{\rm
h}(\p_tu,\e^2\p_tv)_\Phi  | \D_k^\h (u,v)_\Phi\bigr)_{L^2}\\
&+\bigl(\D_k^{\rm h}(\p_tu,\e^2\p_tv)_\Phi  | \D_k^\h (u,v)_\Phi\bigr)_{L^2}-\e^2\bigl(\D_k^\h\p_x^2(u,\e^2v)_\Phi | \D_k^{\rm h} (u,v)_\Phi\bigr)_{L^2}\\
&-\bigl(\D_k^\h\p_y^2(u,\e^2v)_\Phi | \D_k^{\rm h} (u,v)_\Phi\bigr)_{L^2}+\bigl(\D_k^\h \na p_\Phi  | \D_k^{\rm h} (u,v)_\Phi\bigr)_{L^2}\\
=&-\bigl(\D_k^\h\left(u\p_xu\right)_\Phi | \D_k^{\rm h} u_\Phi\bigr)_{L^2}-\bigl(\D_k^\h\left(v\p_yu\right)_\Phi | \D_k^{\rm h} u_\Phi\bigr)_{L^2}\\
&-\e^2\bigl(\D_k^\h\left(u\p_xv\right)_{\Phi} | \D_k^\h v_{\Phi}\bigr)_{L^2}-\e^2\bigl(\D_k^\h\left(v\p_yv\right)_{\Phi} | \D_k^\h v_{\Phi}\bigr)_{L^2}.
\end{split}
\eeq
By using integration by parts and \eqref{S4eq5p}, we find
\begin{align*}
\bigl(&\D_k^{\rm h}\p_t(\p_tu,\e^2\p_tv)_\Phi | \D_k^\h (u,v)_\Phi\bigr)_{L^2}\\
&=\f{d}{dt}\bigl(\D_k^{\rm h}(\p_tu,\e \p_tv)_\Phi | \D_k^\h (u,\e v)_\Phi\bigr)_{L^2}-\|\D_k^\h(\p_tu,\e\p_tv)_\Phi\|_{L^2}^2\\
&+\f\lam2\dot{\tht}
\f{d}{dt}\||D_x|^{\f14}\D_k^{\rm h}(u,\e v)_\Phi(t)\|_{L^2}^2+\lam^2\dot{\tht}^2\||D_x|^{\f12}\D_k^{\rm
h}(u,\e v)_\Phi\|^2_{L^2},
\end{align*}
and
\begin{align*}
\bigl(\D_k^{\rm h}(\p_tu,\e^2\p_tv)_\Phi  | \D_k^\h (u,v)_\Phi\bigr)_{L^2}=\f12\f{d}{dt}\|\D_k^{\rm h}(u,\e v)_\Phi(t)\|_{L^2}^2
+\lam\dot{\tht}\||D_x|^{\f14}\D_k^{\rm
h}(u,\e v)_\Phi\|^2_{L^2}.
\end{align*}
By substituting the above inequalities into \eqref{S3eq2}, we write
\beq
\label{S3eq3}
\begin{split}
\f{d}{dt}&\Bigl(\bigl(\D_k^{\rm h}(\p_tu,\e\p_tv)_\Phi(t) | \D_k^\h (u,\e v)_\Phi(t)\bigr)_{L^2}
+\lam\dot{\tht}(t)\||D_x|^{\f14}\D_k^{\rm h}(u,\e v)_\Phi(t)\|_{L^2}^2\\
&+\f12\|\D_k^{\rm h}(u, \e v)_\Phi(t)\|_{L^2}^2\Bigr)-\|\D_k^\h(\p_tu,\e\p_t v)_\Phi\|_{L^2}^2+\lam\bigl(\dot{\tht}-\ddot{\tht}\bigr)\||D_x|^{\f14}\D_k^{\rm
h}(u,\e v)_\Phi\|^2_{L^2}\\
&+\e^2\|\D_k^{\rm
h} \p_x(u,\e v)_\Phi\|^2_{L^2}+\|\D_k^{\rm
h} \p_y(u,\e v)_\Phi\|^2_{L^2}+2\lam^2\dot{\tht}^2\||D_x|^{\f12}\D_k^{\rm h}(u,\e v)_\Phi\|_{L^2}^2\\
=&-\bigl(\D_k^\h\left(u\p_xu\right)_\Phi | \D_k^{\rm h} u_\Phi\bigr)_{L^2}-\bigl(\D_k^\h\left(v\p_yu\right)_\Phi | \D_k^{\rm h} u_\Phi\bigr)_{L^2}\\
&-\e^2\bigl(\D_k^\h\left(u\p_xv\right)_{\Phi} | \D_k^\h v_{\Phi}\bigr)_{L^2}-\e^2\bigl(\D_k^\h\left(v\p_yv\right)_{\Phi} | \D_k^\h v_{\Phi}\bigr)_{L^2}.
\end{split}
\eeq
Let us denote $\frak{F}_1(t)$  by
\begin{align*}
\frak{F}_1(t)\eqdefa &\f12\Bigl(\|\D_k^{\rm h}(\p_tu,\e\p_tv)_\Phi(t)\|_{L^2}^2+{\e^2}\|\D_k^{\rm h}\p_x(u,\e v)_\Phi(t)\|_{L^2}^2+\|\D_k^{\rm h}\p_y(u,\e v)_\Phi(t)\|_{L^2}^2\\
&+\bigl(\D_k^{\rm h}(\p_tu,\e\p_tv)_\Phi(t) | \D_k^\h (u,\e v)_\Phi(t)\bigr)_{L^2}+\lam\dot{\tht}(t)\||D_x|^{\f14}\D_k^{\rm h}(u,\e v)_\Phi(t)\|_{L^2}^2\Bigr)\\
&+\f14\|\D_k^{\rm h}(u,\e v)_\Phi(t)\|_{L^2}^2
+\lam^2\dot{\tht}^2(t)\||D_x|^{\f12}\D_k^{\rm h}(u,\e v)_\Phi(t)\|_{L^2}^2.
\end{align*}

Then we get, by summing up \eqref{S3eq1} with $\f12\times$\eqref{S3eq3}, that
\beq\label{S3eq3a}
\begin{split}
\f{d}{dt}&\frak{F}_1(t)+\f\lam2\bigl(\dot{\tht}-\ddot{\tht}\bigr)\||D_x|^{\f14}\D_k^{\rm
h}(u,\e v)_\Phi\|^2_{L^2}+\lam\dot{\tht}\Bigl(\||D_x|^{\f14}\D_k^{\rm
h}\p_t(u,\e v)_\Phi\|^2_{L^2}\\
&+\e^2\||D_x|^{\f14}\D_k^{\rm
h} \p_x(u,\e v)_\Phi\|^2_{L^2}+\||D_x|^{\f14}\D_k^{\rm
h} \p_y(u,\e v)_\Phi\|^2_{L^2}\Bigr)\\
&+\lam^2\bigl(\dot{\tht}^2-2\dot{\tht}\ddot{\tht}\bigr)\||D_x|^{\f12}\D_k^{\rm h}(u,\e v)_\Phi\|_{L^2}^2+\lam^3\dot{\tht}^3\||D_x|^{\f34}\D_k^{\rm
h}(u,\e v)_\Phi\|^2_{L^2}\\
&+\f12\Bigl(\|\D_k^{\rm h}(\p_tu,\e \p_tv)_\Phi\|_{L^2}^2+{\e^2}\|\D_k^{\rm h}\p_x(u,\e v)_\Phi\|_{L^2}^2+\|\D_k^{\rm
h} \p_y(u,\e v)_\Phi\|^2_{L^2}\Bigr)\\
=&-\bigl(\D_k^\h\left(u\p_xu\right)_{\Phi} | \D_k^\h (\p_tu+\f12u)_{\Phi}\bigr)_{L^2}-\bigl(\D_k^\h\left(v\p_yu\right)_{\Phi} | \D_k^\h (\p_tu+\f12u)_{\Phi}\bigr)_{L^2}\\
&-\e^2\bigl(\D_k^\h\left(u\p_xv\right)_{\Phi} | \D_k^\h (\p_tv+\f12v)_{\Phi}\bigr)_{L^2}-\e^2\bigl(\D_k^\h\left(v\p_yv\right)_{\Phi} | \D_k^\h (\p_tv+\f12v)_{\Phi}\bigr)_{L^2}.
\end{split}
\eeq

While due to $\left(u_{\Phi}, v_{\Phi}\right)|_{y=0}=\left(u_{\Phi}, v_{\Phi}\right)|_{y=1}=0,$ by applying Poincar\'e inequality, we have
\beq \label{S3eq3b}
\|\D_k(u_{\Phi},\e v_{\Phi})\|_{L^2}^2\leq {\cK}\bigl\|\p_y\D_k(u_{\Phi}, \e v_{\Phi})\bigr\|_{L^2}^2.
\eeq
Let $\fk\eqdefa\min\left(\f1{6},\f1{4(1+\cK)}\right).$
Then, by virtue of \eqref{S4eq5e}, we get, by  using Lemma \ref{lem:Bern} and by multiplying \eqref{S3eq3a} by $e^{2\frak{K}t}$ and then  integrating the resulting inequality over $[0,t],$ that
\beq \label{S3eq6qw}
\begin{split}
{\rm G}_k(u,v)(t)
\leq  &\f34\bigl\|e^{a|D_x|^{\f12}}\D_k^{\rm h}(u_1,\e v_1)(0)\bigr\|_{L^2}^2+\f{\e^2}2\bigl\|e^{a|D_x|^{\f12}}\D_k^{\rm h}\p_x(u_0,\e v_0)\bigr\|_{L^2}^2\\
&+\f12\bigl\|e^{a|D_x|^{\f12}}\D_k^{\rm h}\p_y(u_0,\e v_0)\bigr\|_{L^2}^2+\f12\bigl\|e^{a|D_x|^{\f12}}\D_k^{\rm h}(u_0,\e v_0)\bigr\|_{L^2}^2\\
&+\f\lam2\de^{\f12}2^{\f{k}2}\bigl\|e^{a|D_x|^{\f12}}\D_k^{\rm h}(u_0,\e v_0)\bigr\|_{L^2}^2+\lam^2\de2^k\bigl\|e^{a|D_x|^{\f12}}\D_k^{\rm h}(u_0,\e v_0)\bigr\|_{L^2}^2\\
&+\int_0^t\Bigl(\bigl|\bigl(\ektp\D_k^\h\left(u\p_xu+v\p_yu\right)_\Phi | \ektp\D_k^\h (\p_tu+\f12u)_\Phi\bigr)_{L^2}\bigr|\\
&\qquad\quad+\e^2\bigl|\bigl(\ektp\D_k^\h\left(u\p_xv+v\p_yv\right)_\Phi | \ektp\D_k^\h (\p_tv+\f12v)_\Phi\bigr)_{L^2}\bigr|\Bigr)\,dt'.
\end{split}
\eeq
where ${\rm G}_k(u,v)(t)$ is defined  by
\beq \label{S3eq7}
\begin{split}
{\rm G}_k&(u,v)(t)\eqdefa \f16\|\ektp\D_k^{\rm h}(\p_tu,\e\p_tv)_\Phi\|_{L^\infty_t(L^2)}^2+\f{\e^2}2\|\ektp\D_k^{\rm h}\p_x(u,\e v)_\Phi\|_{L^\infty_t(L^2)}^2\\
&+\f12\|\ektp\D_k^{\rm h}\p_y(u,\e v)_\Phi\|_{L^\infty_t(L^2)}^2+\f1{16}\|\ektp\D_k^{\rm h}(u,\e v)_\Phi\|_{L^\infty_t(L^2)}^2\\
&+\f\lam2\dot{\tht}(t)\|\ektp|D_x|^{\f14}\D_k^{\rm h}(u,\e v)_\Phi(t)\|_{L^2}^2+{\lam^2}\dot{\tht}^2(t)\|\ektp|D_x|^{\f12}\D_k^{\rm h}(u,\e v)_\Phi(t)\|_{L^2}^2\\
&+\f{\lam}4\int_0^t\dot{\tht}(t')\bigl\|\ektp|D_x|^{\f14}\D_k^{\rm
h}\bigl((u,\e v)_\Phi,\p_t(u,\e v)_\Phi,
\p_y(u,\e v)_\Phi\bigr)(t')\|^2_{L^2}\bigr)\,dt'\\
&+\f{\lam^2}2\int_0^t\dot{\tht}^2(t')\|\ektp|D_x|^{\f12}\D_k^{\rm h}(u,\e v)_\Phi(t')\|_{L^2}^2\,dt'\\
&+\lam^3\int_0^t\dot{\tht}^3(t')\|\ektp|D_x|^{\f34}\D_k^{\rm
h}(u,\e v)_\Phi(t')\|^2_{L^2}\,dt'\\
&+\f14\bigl\|\ektp\D_k^{\rm h}\bigl((\p_tu,\e\p_tv)_\Phi,
\e\p_x(u,\e v)_\Phi,\p_y(u,\e v)_\Phi\bigr)\bigr\|^2_{L^2_t(L^2)}.
\end{split}
\eeq

In what follows, we shall always assume that $t<T_2^\star$ with
$T^\star_2$ being determined by
\beq\label{S3eq20} T_2^\star\eqdefa
\sup\bigl\{\ t>0,\ \  \|(\p_y u,\e\p_xu)_\Phi(t)\|_{\cB^{\f12}}\leq \de e^{-\frak{K}t} \  \bigr\}. \eeq

Then it follows from Lemma \ref{lem3.1} that  $t\leq T_2^\star$
\beno
\e^2\int_0^t\bigl|\bigl(\ektp\D_k^\h(u\p_x v)_\Phi\ |\ \ektp\D_k^\h v_\Phi\bigr)_{L^2}\bigr|\,dt'\lesssim d_k^2
2^{-k}\|\ektp \e v_\Phi\|_{\wt{L}^2_{t,\dot{\tht}^3(t)}(\cB^{\frac54})}\|\ektp \e v_\Phi\|_{\wt{L}^2_{t,\dot{\tht}(t)}(\cB^{\frac34})},
\eeno
and
\begin{align*}
\e^2\int_0^t&\bigl|\bigl(\ektp\D_k^\h\left(u\p_xv\right)_\Phi | \ektp\D_k^\h (\p_tv)_\Phi\bigr)_{L^2}\bigr|\,dt'\\
&\lesssim d_k^2
2^{-k}\|\ektp \e v_\Phi\|_{\wt{L}^2_{t,\dot{\tht}^3(t)}(\cB^{\frac54})}\bigl(\|\ektp \e\p_tv_\Phi\|_{\wt{L}^2_{t,\dot{\tht}(t)}(\cB^{\frac34})}
+\lam \|\ektp \e v_\Phi\|_{\wt{L}^2_{t,\dot{\tht}^3(t)}(\cB^{\frac54})}\bigr).
\end{align*}

The estimate of the remaining terms in \eqref{S3eq6qw} relies on the following lemma, which we admit for the time being.

\begin{lem}\label{lem3.3}
{\sl For $t\leq T_2^\star,$ there holds
\beq \label{S3eq10}
\begin{split}
\e\int_0^t\bigl|\bigl(\ektp\D_k^\h(v\p_y v)_{\Phi}\ &|\ \ektp\D_k^\h \fb_{\Phi}\bigr)_{L^2}\bigr|\,dt'\\
&\lesssim d_k^2
2^{-k}\bigl\|\ektp (u, \e v)_{\Phi}\bigr\|_{\wt{L}^2_{t,\dot{\tht}^3(t)}(\cB^{\f54})}\|\ektp \fb_{\Phi}\|_{\wt{L}^2_{t,\dot{\tht}(t)}(\cB^{\f34})}.
\end{split}
\eeq}
\end{lem}

We deduce  from Lemma \ref{lem3.3} that for  $t\leq T_2^\star$
\begin{align*}
\e^2&\int_0^t\bigl|\bigl(\ektp\D_k^\h\left(v\p_yv\right)_\Phi | \ektp\D_k^\h (\p_tv+\f12v)_\Phi\bigr)_{L^2}\bigr|\,dt'\\
&\lesssim d_k^2
2^{-k}\|\ektp (u,\e v)_\Phi\|_{\wt{L}^2_{t,\dot{\tht}^3(t)}(\cB^{\frac54})}\bigl(\|\ektp \e(v_\Phi,\p_tv_\Phi)\|_{\wt{L}^2_{t,\dot{\tht}(t)}(\cB^{\frac34})}
+\lam \|\ektp \e v_\Phi\|_{\wt{L}^2_{t,\dot{\tht}^3(t)}(\cB^{\frac54})}\bigr).
\end{align*}

Without loss of generality, we may assume that $\lam\geq 1,$ then
by
inserting the above estimates and (\ref{S4eq9a}-\ref{S4eq9dqw}) into \eqref{S3eq6qw}, we find
\beno
\begin{split}
{\rm G}_k(u,v)(t)
\leq &  \f34\bigl\|e^{a|D_x|^{\f12}}\D_k^{\rm h}(u_1,\e v_1)(0)\bigr\|_{L^2}^2+\f{\e^2}2\bigl\|e^{a|D_x|^{\f12}}\D_k^{\rm h}\p_x(u_0,\e v_0)\bigr\|_{L^2}^2\\
&+\f12\bigl\|e^{a|D_x|^{\f12}}\D_k^{\rm h}\p_y(u_0,\e v_0)\bigr\|_{L^2}^2+\f12\bigl\|e^{a|D_x|^{\f12}}\D_k^{\rm h}(u_0,\e v_0)\bigr\|_{L^2}^2\\
&+\f\lam2\de^{\f12}2^{\f{k}2}\bigl\|e^{a|D_x|^{\f12}}\D_k^{\rm h}(u_0,\e v_0)\bigr\|_{L^2}^2+\lam^2\de2^k\bigl\|e^{a|D_x|^{\f12}}\D_k^{\rm h}(u_0,\e v_0)\bigr\|_{L^2}^2\\
&+ Cd_k^2
2^{-k}\Bigl(\bigl\|\ektp \bigl((u,\e v)_\Phi, \p_tu_\Phi,\e \p_tv_\Phi\bigr)\bigr\|_{\wt{L}^2_{t,\dot{\tht}(t)}(\cB^{\frac34})}^2+\lam\|\ektp (u,\e v)_\Phi\|_{\wt{L}^2_{t,\dot{\tht}^3(t)}(\cB^{\frac54})}^2
\Bigr).
\end{split}
\eeno
By multiplying the above inequality by $2^k$ and then
taking square root of the resulting inequality,  and   summing up the inequalities  over $\Z,$ we find that for $t\leq T_2^\ast$
\begin{align*}
\frak{E}^\e_{\fk,\lam}(u,v)
 &\leq  C\Bigl(\bigl\|e^{a|D_x|^{\f12}}\bigl(u_0,\e v_0,\e\p_x(u_0,\e v_0),\p_y(u_0,\e v_0), u_1,\e v_1\bigr)\bigr\|_{\cB^{\f12}}\\
&\quad+\sqrt{a\fk}\bigl\|e^{a|D_x|^{\f12}}(u_0,\e v_0)\bigr\|_{\cB^{\f34}}+a\fk\bigl\|e^{a|D_x|^{\f12}}(u_0,\e v_0)\bigr\|_{\cB^{1}}\\
&\quad+\bigl\|\ektp \bigl((u,\e v)_\Phi, \p_tu_\Phi,\e \p_tv_\Phi\bigr)\bigr\|_{\wt{L}^2_{t,\dot\tht(t)}(\cB^{\f34})}
+\lam^{\f12}\|\ektp (u,\e v)_\Phi\|_{\wt{L}^2_{t,\dot\tht^3(t)}(\cB^{\f54})}\Bigr),
\end{align*}
where we used \eqref{S4eq-1} so that $\lam\de^{\f12}=\f{a\fk}4,$
and
\beq\label{Senerg2}
\begin{split}
\frak{E}^\e_{\fk,\lam}&(u,v)\eqdefa \bigl\|\ektp \bigl(u,\e v,\e\p_x(u,\e v), \p_y(u,\e v), (\p_tu,\e\p_tv)\bigr)_\Phi\bigr\|_{\wt{L}^\infty_t(\cB^{\f12})}\\
+&\sqrt{a\fk}\|e^{\f34\fk t'} (u,\e v)_\Phi\|_{\wt{L}^\infty_{t}(\cB^{\f34})}+\bigl\|\ektp \bigl(\e\p_x(u,\e v),\p_y (u,\e v), \p_t( u, \e v)\bigr)_\Phi\bigr\|_{\wt{L}^2_t(\cB^{\f12})}\\
+&a\fk\|e^{\f12\fk t'} (u,\e v)_\Phi\|_{\wt{L}^\infty_{t}(\cB^{1})}+\sqrt{\lam}\bigl\|\ektp \bigl((u,\e v,  \p_y(u,\e v))_\Phi, (\p_tu_\Phi,\e \p_tv_\Phi)\bigr)\bigr\|_{\wt{L}^2_{t,\dot\tht(t)}(\cB^{\f34})}\\
+&{\lam}\|\ektp (u,\e v)_\Phi\|_{\wt{L}^2_{t,\dot\tht^2(t)}(\cB^{1})}+{\lam}^{\f32}\|\ektp (u,\e v)_\Phi\|_{\wt{L}^2_{t,\dot\tht^3(t)}(\cB^{\f54})}
\end{split}
\eeq
Taking $\lam =\max\left(C^2,1\right)$ in the above inequality leads to \eqref{S1eq16}
for $t\leq T^\star_2.$

In particular, we deduce from \eqref{S1eq8} and \eqref{S1eq16}   that
\beq\label{S4eq19}
\begin{split}
\|(\e\p_x u, \p_yu)_\Phi(t)\|_{\cB^{\f12}}\leq &C\left(1+\sqrt{a\fk}+a\fk\right)c_2 e^{-\fk t}\\
\leq &\f\de2 e^{-\fk t} \quad \mbox{for any}\ t\leq T^\star_2,
\end{split}
\eeq
 if we take $c_1$ in \eqref{S1eq8} to be so small that
$C\left(1+\sqrt{a\fk}+a\fk\right)c_2\leq  \f\de2.$  Then
we deduce by a continuous argument that $T^\star_2$ determined by \eqref{S3eq20} equals $+\infty$ and \eqref{S1eq16} holds for any $t>0.$
This completes the proof of Theorem \ref{thm1.1}.
\end{proof}

\begin{proof}[Proof of Lemma \ref{lem3.3}] Once again we shall prove \eqref{S3eq10} for $\fk=0.$ By applying Bony's decomposition \eqref{Bony}  to $v\p_yv$, we write
\beno
v\p_yv=T^\h_{v}\p_yv+T^\h_{\p_yv}v+R^h(v,\p_yv).
\eeno
Let us handle the following three terms:\smallskip

\no $\bullet$ \underline{Estimate of
$\e\int_0^t\bigl|\bigl(\D_k^{\rm h}(T^\h_{v}\p_yv)_{\Phi}\ |\ \D_k^{\rm
h}\fb_{\Phi}\bigr)_{L^2}\bigr|\,dt'$}

Due to $\p_yv=-\p_xu,$ for $t\leq T_2^\star,$ one deduces from \eqref{S4eq-1} and \eqref{S3eq20} that
\beno
\begin{split}
\e\int_0^t\bigl|\bigl(&\D_k^{\rm h}(T^\h_{v}\p_yv)_{\Phi}\ |\ \D_k^{\rm
h}\fb_{\Phi}\bigr)_{L^2}\bigr|\,dt'\\
\lesssim & \e\sum_{|k'-k|\leq 4}\int_0^t\|S_{k'-1}^\h v_{\Phi}(t')\|_{L^\infty}
\|\D_{k'}^\h\p_yv_{\Phi}(t)\|_{L^2}\|\D_k^\h \fb_{\Phi}(t')\|_{L^2}\,dt'\\
\lesssim & \sum_{|k'-k|\leq 4}2^{-\f{k'}2}\int_0^t\|S_{k'-1}^\h v_{\Phi}(t')\|_{L^\infty}\e\|\p_xu_{\Phi}(t')\|_{\cB^{\f12}}
\|\D_k^\h \fb_{\Phi}(t')\|_{L^2}\,dt'\\
\lesssim & \sum_{|k'-k|\leq 4}2^{-\f{k'}2}\Bigl(\int_0^t\|S_{k'-1}^\h v_{\Phi}(t')\|_{L^\infty}^2\dot{\tht}^3(t')\,dt'\Bigr)^{\f12}\Bigl(\int_0^t\|\D_k^\h \fb_{\Phi}(t')\|_{L^2}^2\dot{\tht}(t')\,dt'\Bigr)^{\f12},
\end{split}
\eeno
from which, \eqref{S3eq12} and Definition \ref{def1.1}, we infer
\beno
\e\int_0^t\bigl|\bigl(\D_k^{\rm h}(T^\h_{v}\p_yv)_{\Phi}\ |\ \D_k^{\rm
h}\fb_{\Phi}\bigr)_{L^2}\bigr|\,dt'\leq d_k^22^{-k}\|u_{\Phi}\|_{\wt{L}^2_{t,\dot{\tht}^3(t)}(\cB^{\f54})}\|\fb_{\Phi}\|_{\wt{L}^2_{t,\dot{\tht}(t)}(\cB^{\f34})}.
\eeno

\no $\bullet$ \underline{Estimate of
$\int_0^t\bigl|\bigl(\D_k^{\rm h}(T^\h_{\p_yv}v)_{\Phi}\ |\ \D_k^{\rm
h}\fb_{\Phi}\bigr)_{L^2}\bigr|\,dt'$}

Observing that
\beno
\begin{split}
\int_0^t\bigl|\bigl(&\D_k^{\rm h}(T^\h_{\p_yv}v)_{\Phi}\ |\ \D_k^{\rm h}\fb_{\Phi}\bigr)_{L^2}\bigr|\,dt'\\
\lesssim & \sum_{|k'-k|\leq 4}\int_0^t\|S_{k'-1}^\h \p_x u_{\Phi}(t')\|_{L^\infty}
\|\D_{k'}^\h v_{\Phi}(t)\|_{L^2}\|\D_k^\h \fb_{\Phi}(t')\|_{L^2}\,dt',
\end{split}
\eeno
from which, Lemma \ref{lem:Bern}, \eqref{S3eq8} and \eqref{S3eq20}, we deduce that
\beno
\begin{split}
\int_0^t\bigl|\bigl(&\D_k^{\rm h}(T^\h_{\p_yv}v)_{\Phi}\ |\ \D_k^{\rm h}\fb_{\Phi}\bigr)_{L^2}\bigr|\,dt'\\
\lesssim & \sum_{|k'-k|\leq 4}2^{k'}\int_0^t\|S_{k'-1}^\h u_{\Phi}(t')\|_{L^\infty}
\|\D_{k'}^\h v_{\Phi}(t)\|_{L^2}\|\D_k^\h \fb_{\Phi}(t')\|_{L^2}\,dt'\\
\lesssim & \sum_{|k'-k|\leq 4}2^{k'}\int_0^t\|\p_y u_{\Phi}(t')\|_{\cB^{\f12}}
\|\D_{k'}^\h v_{\Phi}(t)\|_{L^2}\|\D_k^\h \fb_{\Phi}(t')\|_{L^2}\,dt'\\
\lesssim & \sum_{|k'-k|\leq 4}2^{k'}\Bigl(\int_0^t
\|\D_{k'}^\h v_{\Phi}(t)\|_{L^2}^2\dot{\tht}^3(t')\,dt'\Bigr)^{\f12}\Bigl(\int_0^t\|\D_k^\h \fb_{\Phi}(t')\|_{L^2}^2\dot{\tht}(t')\,dt'\Bigr)^{\f12}
\end{split}
\eeno
Then thanks to Definition \ref{def1.1}, we achieve
\beno
\begin{split}
\int_0^t\bigl|\bigl(\D_k^{\rm h}(T^\h_{\p_yv}v)_{\Phi}\ |\ \D_k^{\rm
h}\fb_{\Phi}\bigr)_{L^2}\bigr|\,dt' \lesssim & d_k^22^{-k}\|v_{\Phi}\|_{\wt{L}^2_{t,\dot{\tht}^3(t)}(\cB^{\f54})}\|\fb_{\Phi}\|_{\wt{L}^2_{t,\dot{\tht}(t)}(\cB^{\f34})}.
\end{split}
\eeno

\no $\bullet$ \underline{Estimate of
$\int_0^t\bigl|\bigl(\D_k^{\rm h}(R^\h(v,\p_yv))_{\Phi}\ |\ \D_k^{\rm
h}\fb_{\Phi}\bigr)_{L^2}\bigr|\,dt'$}\vspace{0.2cm}

Due to $\p_xu+\p_yv=0,$ we get, by applying lemma
\ref{lem:Bern} and \eqref{S3eq8}, that
\beno
\begin{split}
\int_0^t\bigl|&\bigl(\D_k^{\rm h}(R^\h(v,\p_yv))_{\Phi}\ |\ \D_k^{\rm
h}\fb_{\Phi}\bigr)_{L^2}\bigr|\,dt'\\
\lesssim &2^{\f{k}2}\sum_{k'\geq k-3}\int_0^t\|{\D}_{k'}^\h v_{\Phi}(t')\|_{L^2}\|\wt{\D}_{k'}^\h \p_x u_{\Phi}(t')\|_{L^2_\h(L^\infty_{\rm v})}\|\D_k^\h \fb_{\Phi}(t')\|_{L^2}\,dt'\\
\lesssim & 2^{\f{k}2}\sum_{k'\geq k-3}2^{\f{k'}2}\int_0^t\|{\D}_{k'}^\h v_{\Phi}(t')\|_{L^2}\|\p_yu_{\Phi}(t')\|_{\cB^{\f12}}\|\D_k^\h \fb_{\Phi}(t')\|_{L^2}\,dt'\\
\lesssim & 2^{\f{k}2}\sum_{k'\geq k-3}2^{\f{k'}2}\Bigl(\int_0^t\|{\D}_{k'}^\h  v_{\Phi}(t')\|_{L^2}^2\dot{\tht}^3(t')\,dt'\Bigr)^{\f12}
 \Bigl(\int_0^t\|\D_k^\h v_{\Phi}(t')\|_{L^2}^2\dot{\tht}(t')
 \,dt'\Bigr)^{\f12},
\end{split}
\eeno
which together with Definition \ref{def1.1} ensures that
\beno
\begin{split}
\int_0^t\bigl|\bigl(\D_k^{\rm h}(R^\h(v,\p_yu))_{\Phi}\ |\ \D_k^{\rm
h}\fb_{\Phi}\bigr)_{L^2}\bigr|\,dt'
\lesssim d_k^22^{-k}\|v_{\Phi}\|_{\wt{L}^2_{t,\dot{\tht}^3(t)}(\cB^{\f54})}\|\fb_{\Phi}\|_{\wt{L}^2_{t,\dot{\tht}^3(t)}(\cB^{\f34})}.
\end{split}
\eeno

By summing up the above estimates, we obtain \eqref{S3eq10}. This concludes the proof of Lemma \ref{lem3.3}.
\end{proof}

\setcounter{equation}{0}

\section{The Convergence to the hyperbolic Prandtl system} \label{Sect5}

In this section, we justify the limit from the scaled anisotropic hyperbolic Navier-Stokes system \eqref{S1eq1.8}
 to the hydrostatic hyperbolic Navier-Stokes system \eqref{S1eq1} in a 2-D striped domain. To this end, we introduce
\beno
w^1_\e\eqdefa u^\e-u,\quad w^2_\e\eqdefa v^\e-v,\quad q_\e\eqdefa p^\e-p,
\eeno
where $(u^\e,v^\e,p^\e)$ (resp. $(u,v,p)$) are solutions to the system \eqref{S1eq1.8} (resp. \eqref{S1eq1}) obtained in Theorem
\ref{thm1.1} (resp. Theorem \ref{th1.2}).
Then $(w^1_\e,w^2_\e,q_\e)$ verifies
\begin{equation}\label{S5eq1}
 \quad\left\{\begin{array}{l}
\displaystyle \p_t^2 w^1_\e+\p_tw^1_\e-\e^2\p_x^2w^1_\e-\p_y^2w^1_\e+\p_xq_\e=R^1_\e\ \quad \mbox{in} \ \ \cS\times ]0,\infty[,\\
\displaystyle \e^2\left(\p_t^2 w^2_\e+\p_tw^2_\e-\e^2\p_x^2w^2_\e-\p_y^2w^2_\e\right)+\p_yq_\e=R^2_\e,\\
\displaystyle \p_xw^1_\e+\p_yw^2_\e=0,\\
\displaystyle \left(w^1_\e,w^2_\e\right)|_{y=0}=\left(w^1_\e,w^2_\e\right)|_{y=1}=0,\\
\displaystyle \left(w^1_\e, w^2_\e\right)|_{t=0}=\left(u_0^\e-u_0, v_0^\e-v_0\right),\quad \left(\p_tw^1_\e, \p_tw^2_\e\right)|_{t=0}=\left(u_1^\e-u_1, v_1^\e-v_1\right).
\end{array}\right.
\end{equation}
where $v_i, i=0,1,$ is determined from $u_i$ via $\p_xu_i+\p_yv_i=0$ and $v_i|_{y=0}=v_i|_{y=1}=0,$ and
\beq\label{S5eq2}
\begin{split}
&R^1_\e=\e^2\p_x^2u-\big(u^\e\p_xu^\e-u\p_xu\big)-\big(v^\e\p_yu^\e-v\p_yu\big),\\
&R^2_\e=-\e^2\big(\p_t^2v+\p_tv-\e^2\p_x^2v-\p_y^2v+u^\e\p_xv^\e+v^\e\p_yv^\e\big).
\end{split} \eeq

We now present the proof of Theorem \ref{thm3}.
In what follows, we shall neglect the subscript $\e$ in $(w^1_\e,w^2_\e).$

\begin{proof}[Proof of Theorem \ref{thm3}]
In view of \eqref{S5eq1}, we get, by using a similar derivation of \eqref{S3eq6qw}, that
\beq \label{S5eq3}
\begin{split}
{\rm G}_k(&w^1,w^2)(t)
\leq  \f34\bigl\|e^{a|D_x|^{\f12}}\D_k^{\rm h}(w_0^1,\e w_0^2)\bigr\|_{L^2}^2+\f{\e^2}2\bigl\|e^{a|D_x|^{\f12}}\D_k^{\rm h}\p_x(w^1_0,\e w^2_0)\bigr\|_{L^2}^2\\
&+\f12\bigl\|e^{a|D_x|^{\f12}}\D_k^{\rm h}\p_y(w^1_0,\e w^2_0)\bigr\|_{L^2}^2+\f1{2}\bigl\|e^{a|D_x|^{\f12}}\D_k^{\rm h}(w^1_0,\e w^2_0)\bigr\|_{L^2}^2\\
&+\f\lam2\de^{\f12}2^{\f{k}2}\bigl\|e^{a|D_x|^{\f12}}\D_k^{\rm h}(w^1_0,\e w^2_0)\bigr\|_{L^2}^2+\f{\lam^2}{2}\de2^k\bigl\|e^{a|D_x|^{\f12}}\D_k^{\rm h}(w^1_0,\e w^2_0)\bigr\|_{L^2}^2\\
&+\int_0^t\Bigl(\bigl|\bigl(\D_k^\h R^1_\Phi | \D_k^\h \bigl(\p_t w^1+\frac12 w^1\bigr)_\Phi\bigr)_{L^2}\bigr|+\bigl|\bigl(\D_k^\h R^2_\Phi | \D_k^\h \bigl(\p_t w^2+\frac 12 w^2\bigr)_\Phi\bigr)_{L^2}\bigr|\Bigr)\,dt'.
\end{split}
\eeq
where  the functional ${\rm G}_k(w^1,w^2)(t)$ is defined  by \eqref{S3eq7}.

We now claim that
\beq \label{S5eq7}
\begin{split}
\int_0^t&\bigl|\bigl(\D_k^\h R^1_\Phi | \D_k^\h \bigl(\p_tw^1+\f12w^1\bigr)_\Phi\bigr)_{L^2}\bigr|\,dt'\\
&\lesssim  d_k^22^{-k}
\Bigl(\e^2\|u_{\Phi}\|_{\wt{L}^2_{t}(\cB^{\f52})}\|(w^1,\p_tw^1)_{\Phi}\|_{\wt{L}^2_{t}(\cB^{\f12})}
+\bigl(\|u_\Th\|^{\f12}_{L^\infty_t(\cB^{\f32})}
\|\p_yw^1_\Th\|_{\wt{L}^2_t(\cB^{\f12})}
\\
&\qquad\qquad\ \ +\|w^1_\Th\|_{\wt{L}^2_{t,\dot\tht^3(t)}(\cB^{\f54})}
\bigr)\bigl(\|(w^1_\Phi,\p_tw^1_\Th)\|_{\wt{L}^2_{t,\dot\tht(t)}(\cB^{\f34})}+\lam \| w^1_\Th\|_{\wt{L}^2_{t,\dot{\tht}^3(t)}(\cB^{\f54})}\bigr)
\Bigr),
\end{split}
\eeq
and
\beq \label{S5eq8}
\begin{split}
\int_0^t\bigl|&\bigl(\D_k^\h R^2_\Phi | \D_k^\h \bigl(\p_tw^1+\f12w^1\bigr)_\Phi\bigr)_{L^2}\bigr|\,dt'\\
 \lesssim & d_k^22^{-k}
 \Bigl(\e^2\bigl(\bigl\|(\p_t^2u, \p_tu, \p_yu_\Th)_\Th\bigr\|_{\wt{L}^2_{t}(\cB^{\f32})}+ \e^2\|u_\Th\|_{\wt{L}^2_{t}(\cB^{\f72})}\bigr)\|(\p_yw^2,\p_tw^2)_\Th\|_{\wt{L}^2_{t}(\cB^{\f12})}
 \\
 &+\e\bigl(\e\|u^\e_\Th\|_{L^\infty_t(\cB^{\f12})}^{\f12}\| u_\Th\|_{\wt{L}_t^2(\cB^{\f94})}+ \bigl\|(u, \e v, w^1_\Th,\e w^2_\Th)\bigr\|_{\wt{L}^2_{t,\dot\tht^3(t)}(\cB^{\f54})}\\
 &\quad
 +\e\|u_\Th\|^{\f12}_{L^\infty_t(\cB^{2})}
\|\p_yw^2_\Th\|_{\wt{L}^2_t(\cB^{\f12})}
\bigr)
 \bigl(\|(w^2_\Phi,\p_tw^2_\Th)\|_{\wt{L}^2_{t,\dot\tht(t)}(\cB^{\f34})}+\lam \| w^2_\Th\|_{\wt{L}^2_{t,\dot{\tht}^3(t)}(\cB^{\f54})}\bigr)\Bigr),
 \end{split}
\eeq
the proof of which will be postponed after we finish the proof of Theorem \ref{thm3}.

On the other hand, due to $\frak{H}_3^0(u_0,u_1)<\infty,$ we deduce from \eqref{S4eq10} and \eqref{S1eq12} that there exists some
positive constant $M$ so that
\beq\label{S5eqop}
{\rm E}_{\f12}(u)(t)+{\rm E}_3(u)(t)\leq M,
\eeq
for the energy functional ${\rm E}_s(u)$ being determined by \eqref{S1eq10}.

Moreover, it follows from \eqref{S1eq14} and Theorem \ref{thm1.1} that
\beq\label{S5eqopp}
{\rm E}_{\f32}(\p_yu)(t)+\|\ektp(\p_t^2u)_\Phi\|_{\wt{L}^2_t(\cB^{\f32})}+{\rm E}_\fk^1(u^\e,v^\e)\leq M.
\eeq

In particular, we deduce from \eqref{S5eqop} that
\beq \label{S5eqou}
\begin{split}
\|u_\Phi\|_{\wt{L}^2_t(\cB^{\f72})}=\sum_{j\in\Z}2^{\f72k}\|\D_k^\h u_\Phi\|_{L^2_t(L^2)}\lesssim& \sum_{j\in\Z}2^{\f72k}\|e^{\f{\fk}2t}\D_k^\h u_\Phi\|_{L^\infty_t(L^2)}\\
\lesssim& \|e^{\f{\fk}2t} u_\Phi\|_{\wt{L}^\infty_t(\cB^{\f72})}\leq C{\rm E}_3(u)(t)
\leq CM.
\end{split}
\eeq
Similarly, we have
\beq \label{S5eqol}
\|u_\Phi\|_{\wt{L}^2_t(\cB^{\f94})}+\|u_\Phi\|_{\wt{L}^2_t(\cB^{\f52})}\leq CM.
\eeq

Let $\frak{E}^\e_{0,\lam}(w^1,w^2)$ be given by \eqref{Senerg2}. Thanks to (\ref{S5eqop}-\ref{S5eqol}), by substituting \eqref{S5eq7} and \eqref{S5eq8} into \eqref{S5eq3}, and  taking square root of the resulting inequality,
and then multiplying it by $2^{\f{k}2}$ and summing up the inequalities for $k$ over $\Z,$ we obtain
\begin{align*}
\frak{E}^\e_{0,\lam}(&w^1,w^2) \leq  C\Bigl({\frak H}_0^1(w^1_0,w^2_0)+(M\e)^{\f12}\bigl(1+\bigl\|\bigl(w^1_\Phi,\p_tw^1_\Th, \e w^2_\Phi, \e \p_y w^2_\Phi, \e\p_tw^2_\Phi\bigr)\bigr\|_{\wt{L}^2_{t}(\cB^{\f12})}^{\f12}\bigr)\\
&+\bigl(\|w^1_\Th\|_{\wt{L}^2_{t,\dot\tht^3(t)}(\cB^{\f54})}
+M^{\f12}
\|\p_yw^1_\Th\|_{\wt{L}^2_t(\cB^{\f12})}
\bigr)^{\f12}\bigl(\|(w^1_\Phi,\p_tw^1_\Th)\|_{\wt{L}^2_{t,\dot\tht(t)}(\cB^{\f34})}+\lam \| w^1_\Th\|_{\wt{L}^2_{t,\dot{\tht}^3(t)}(\cB^{\f54})}\bigr)^{\f12}\\
&+\e^{\f12}\bigl(M\e+\bigl\|(w^1_\Th,\e w^2_\Th)\bigr\|_{\wt{L}^2_{t,\dot\tht^3(t)}(\cB^{\f54})}
 +\e M^{\f12}
\|\p_yw^2_\Th\|_{\wt{L}^2_t(\cB^{\f12})}
\bigr)^{\f12}\\
&\qquad\qquad\qquad\qquad\qquad\qquad\qquad\qquad\times \bigl(\|(w^2_\Phi,\p_tw^2_\Th)\|_{\wt{L}^2_{t,\dot\tht(t)}(\cB^{\f34})}+\lam \| w^2_\Th\|_{\wt{L}^2_{t,\dot{\tht}^3(t)}(\cB^{\f54})}\bigr)^{\f12}\Bigr).
\end{align*}
Applying Young's inequality leads to
\beq\label{S5eq30}
\begin{split}
\frak{E}^\e_{0,\lam}&(w^1,w^2) \leq \f12\bigl\|\bigl(w^1_\Phi, \p_yw^1_\Phi,\p_tw^1_\Th, \e w^2_\Phi, \e \p_yw^2_\Phi, \e\p_tw^2_\Phi\bigr)\bigr\|_{\wt{L}^2_{t}(\cB^{\f12})}+C\Bigl({\frak H}_0^1(w^1_0,w^2_0)\\
&\qquad\ \ +M\e
+\bigl\|\bigl(w^1_\Phi,\p_tw^1_\Th, \e w^2_\Phi, \e\p_tw^2_\Th)\|_{\wt{L}^2_{t,\dot\tht(t)}(\cB^{\f34})}
+\lam \bigl\|(w^1, \e w^2)_\Th\|_{\wt{L}^2_{t,\dot{\tht}^3(t)}(\cB^{\f54})}\Bigr).
\end{split}
\eeq
By taking $\lam=4C^2$ in \eqref{S5eq30}, we achieve
\beno
\frak{E}^\e_{0,\f{\lam}4}(w^1,w^2) \leq C\bigl({\frak H}_0^1(w^1_0,w^2_0)+M\e\bigr),
\eeno
which leads to \eqref{S1eq15}. This completes the proof of Theorem \ref{thm3}.
\end{proof}

Let us now present the proof of \eqref{S5eq7} and \eqref{S5eq8}.

\begin{proof}[Proof of \eqref{S5eq7}]
According \eqref{S5eq2}, we write
\beno
R^1_\e=\e^2\p_x^2u-\big(u^\e\p_xw^1+w^1\p_xu\big)-\big(v^\e\p_yw^1+w^2\p_yu\big).
\eeno

We first observe that
\beq \label{S5eq9}
\begin{split}
\e^2\int_0^t\big|\big(\Delta_k^\h\p_x^2u_{\Phi}|\Delta_k^\h (w^1,\p_tw^1)_{\Phi}\big)_{L^2}\big|\,dt'\lesssim d_k^22^{-k}\e^2\| u_{\Phi}\|_{\wt{L}^2_{t}(\cB^{\f52})}\|(w^1,\p_tw^1)_{\Phi}\|_{\wt{L}^2_{t}(\cB^{\f12})}.
\end{split}
\eeq

\noindent$\bullet$\underline{
The estimate of $\int_0^t\big|\big(\Delta_k^\h(u^\e\p_xw^1+w^1\p_xu)_{\Phi} | \D_k^\h \bigl(\p_t w^1+\frac12 w^1\bigr)_\Phi
\big)_{L^2}\big|\,dt'$.}

It follows from \eqref{S1eq16} and Lemma \ref{lem3.1} that
\beq\label{S5eq10}
\begin{split}
\int_0^t&\big|\big(\Delta_k^\h(u^\e\p_xw^1)_{\Phi}|\Delta_k^\h \bigl(\p_tw^1+\f12w^1\bigr)_{\Phi}\big)_{L^2}\big|\,dt'\\
&\lesssim d_k^22^{-k} \| w^1_\Th\|_{\wt{L}^2_{t,\dot{\tht}^3(t)}(\cB^{\f54})}\bigl(\|(w^1_\Th,\p_tw^1_\Phi)\|_{\wt{L}^2_{t,\dot{\tht}(t)}(\cB^{\f34})}
+\lam \| w^1_\Th\|_{\wt{L}^2_{t,\dot{\tht}^3(t)}(\cB^{\f54})}\bigr).
\end{split}
\eeq
While we get, by applying Bony's decomposition \eqref{Bony}  to $w^1\p_xu$, that
\beno
w^1\p_xu=T^\h_{w^1}\p_xu+T^\h_{\p_xu}w^1+R^h(w^1,\p_xu).
\eeno
Notice that
\beno
\|\D_{k'}^\h\p_x u_\Th(t')\|_{L^2_\h(L^\infty_{\rm v})}\lesssim d_{k'}(t)2^{-\f{k'}4}\|u_\Th(t')\|_{\cB^{2}}^{\f12}\|\p_y u_\Th(t')\|_{\cB^{\f12}}^{\f12},
\eeno
we deduce from \eqref{S4eq10}  and \eqref{S4eq-1} that
\beno
\begin{split}
\int_0^t&\bigl|\bigl(\D_k^{\rm h}(T^\h_{w^1}\p_xu)_\Th\ |\ \D_k^{\rm
h}\bigl(\p_tw+\f12w\bigr)_\Th^1\bigr)_{L^2}\bigr|\,dt'\\
\lesssim & \sum_{|k'-k|\leq 4}\int_0^t\|S_{k'-1}^\h w^1_\Th(t')\|_{L^\infty_\h(L^2_{\rm v})}
\|\D_{k'}^\h\p_xu_\Th(t')\|_{L^2_\h(L^\infty_{\rm v})}\|\D_k^\h (w^1,\p_tw^1)_\Th(t')\|_{L^2}\,dt'\\
\lesssim & \sum_{|k'-k|\leq 4}d_{k'}2^{-\f{k'}4}\|u_\Th\|_{L^\infty_t(\cB^{2})}^{\f12}\| w^1_\Th\|_{L^2_t(L^\infty_\h(L^2_{\rm v}))}\Bigl(\int_0^t
\dot\tht(t')\|\D_k^\h (w^1,\p_tw^1)_\Th(t')\|_{L^2}^2\,dt'\Bigr)^{\f12}\\
\lesssim & d_k^22^{-k}\|u_\Th\|^{\f12}_{L^\infty_t(\cB^{2})}
\|\p_yw^1_\Th\|_{\wt{L}^2_t(\cB^{\f12})}\bigl(\|(w^1_\Th,\p_tw^1_\Th)\|_{\wt{L}^2_{t,\dot\tht(t)}(\cB^{\f34})}+\lam \|w^1_\Th\|_{\wt{L}^2_{t,\dot\tht^3(t)}(\cB^{\f54})}\bigr).
\end{split}
\eeno

While we get, by applying Lema \ref{lem:Bern}, that
\beno
\begin{split}
\int_0^t\bigl|&\bigl(\D_k^{\rm h}(R^\h(w^1,\p_xu)_\Th\ |\ \D_k^{\rm
h}w^1_\Th\bigr)_{L^2}\bigr|\,dt'\\
\lesssim &2^{\f{k}2}\sum_{k'\geq k-3}\int_0^t\|{\D}_{k'}^\h w^1_\Th(t')\|_{L^2}\|\wt{\D}_{k'}^\h\p_x u_\Th(t')\|_{L^2_\h(L^\infty_{\rm v})}\|\D_k^\h w^1_\Th(t')\|_{L^2}\,dt'\\
\lesssim & 2^{\f{k}2}\sum_{k'\geq k-3}2^{\f{k'}2}\Bigl(\int_0^t\|\D_{k'}^\h w_\Th^1(t')\|_{L^2}^2\dot\tht^3(t')\,dt'\Bigr)^{\frac12}\\
&\qquad\qquad\qquad\qquad\qquad\times \Bigl(\int_0^t\|\D_k^\h (w^1,\p_tw^1)_\Th(t')\|_{L^2}^2\dot\tht(t')\,dt'\Bigr)^{\frac12}\\
\lesssim & d_k^22^{-k}\|w^1_\Th\|_{\wt{L}^2_{t,\dot\tht^3(t)}(\cB^{\f54})}\bigl(\|(w^1_\Phi,\p_tw^1_\Th)\|_{\wt{L}^2_{t,\dot\tht(t)}(\cB^{\f34})}+\lam \| w^1_\Th\|_{\wt{L}^2_{t,\dot{\tht}^3(t)}(\cB^{\f54})}\bigr).
\end{split}
\eeno
The same estimate holds for $\int_0^t\bigl|\bigl(\D_k^{\rm h}(T^\h_{\p_xu}w^1)_\Th\ |\ \D_k^{\rm
h}\bigl(\p_tw^1+\f12w^1\bigr)_\Th\bigr)_{L^2}\bigr|\,dt'.$

As a result, it comes out
\beq \label{S5eq11}
\begin{split}
\int_0^t\bigl|\bigl(\D_k^{\rm h}&({w^1}\p_xu)_\Th\ |\ \D_k^{\rm
h}\bigl(w_\Th^1+\f12(\p_tw^1)_\Phi\bigr)_{L^2}\bigr|\,dt'\\
\lesssim & d_k^22^{-k} \bigl(\|u_\Th\|^{\f12}_{L^\infty_t(\cB^{2})}
\|\p_yw^1_\Th\|_{\wt{L}^2_t(\cB^{\f12})}+
\|w^1_\Th\|_{\wt{L}^2_{t,\dot\tht^3(t)}(\cB^{\f54})}\bigr)\\
&\qquad\qquad\times\bigl(\|(w^1_\Phi,\p_tw^1_\Th)\|_{\wt{L}^2_{t,\dot\tht(t)}(\cB^{\f34})}+\lam \| w^1_\Th\|_{\wt{L}^2_{t,\dot{\tht}^3(t)}(\cB^{\f54})}\bigr).
\end{split}
\eeq

\noindent$\bullet$\underline{
The estimate of $\int_0^t\big|\big(\Delta_k^\h(v^\e\p_yw^1)_\Th | \Delta_k^\h \bigl(\p_tw^1+\f12w^1\bigr)_\Th\big)_{L^2}\big|\,dt'$.}

We write
\beno
v^\e\p_y w^1=w^2\p_yw^1+v\p_yw^1.
\eeno
We first deduce from \eqref{S4eq10}, \eqref{S1eq16} and Lemma \ref{lem3.2} that
\beq\label{S5eq11a}
\begin{split}
\int_0^t\big|\big(&\Delta_k^\h(w^2\p_yw^1)_\Th | \Delta_k^\h (w^1,\p_tw^1)_\Th\big)_{L^2}\big|\,dt'\\
\lesssim & d_k^22^{-k}\|w_\Th^1\|_{\wt{L}^2_{t,\dot{\tht}^3(t)}(\cB^{\f54})}\bigl(\|(w^1_\Phi,\p_tw^1_\Th)\|_{\wt{L}^2_{t,\dot\tht(t)}(\cB^{\f34})}+\lam \| w^1_\Th\|_{\wt{L}^2_{t,\dot{\tht}^3(t)}(\cB^{\f54})}\bigr).
\end{split}
\eeq

Whereas by applying Bony's decomposition \eqref{Bony} to $v\p_yw^1$, we find
\beno
v\p_y w^1=T^\h_{v}\p_y w^1+T^\h_{\p_yw^1}v+R^h(v,\p_y w^1).
\eeno
It follows from \eqref{S3eq11} that
\beno
\begin{split}
\|S_{k'-1}^\h v_\Th(t')\|_{L^\infty}\lesssim& \sum_{\ell\leq k'-2}2^{\f{3\ell}2}\|\D_{\ell}^\h u_\Th(t')\|_{L^2}^{\f12}
\|\D_{\ell}^\h \p_y u_\Th(t')\|_{L^2}^{\f12}\\
\lesssim &
d_{k'}(t)2^{\f{k'}4}\|u_\Th(t')\|_{\cB^{2}}^{\f12}\|\p_y u_\Th(t')\|_{\cB^{\f12}}^{\f12},
\end{split}
\eeno
from which, we infer
\beno
\begin{split}
\int_0^t&\bigl|\bigl(\D_k^{\rm h}(T^\h_{v}\p_y w^1)_\Th\ |\ \D_k^{\rm
h}\bigl(\p_tw^1+\f12w^1\bigr)_\Th\bigr)_{L^2}\bigr|\,dt'\\
\lesssim & \sum_{|k'-k|\leq 4}\int_0^t\|S_{k'-1}^\h v_\Th(t')\|_{L^\infty}
\|\D_{k'}^\h\p_y w^1_\Th(t')\|_{L^2}\|\D_k^\h (w^1,\p_tw^1)_\Th(t')\|_{L^2}\,dt'\\
\lesssim & \sum_{|k'-k|\leq 4}2^{\f{k'}4}\|u_\Th\|_{L^\infty_t(\cB^{2})}^{\f12}\|\D_{k'}^\h\p_y w^1_\Th(t')\|_{L^2_t(L^2)}
\Bigl(\int_0^t
\dot\tht(t')\|\D_k^\h (w^1,\p_tw^1)_\Th(t')\|_{L^2}^2\,dt'\Bigr)^{\f12}\\
\lesssim & d_k^22^{-k}\|u_\Th\|^{\f12}_{L^\infty_t(\cB^{2})}
\|\p_yw^1_\Th\|_{\wt{L}^2_t(\cB^{\f12})}\bigl(\|(w^1_\Th,\p_tw^1_\Th)\|_{\wt{L}^2_{t,\dot\tht(t)}(\cB^{\f34})}
+\lam\|w^1_\Th\|_{\wt{L}^2_{t,\dot\tht^3(t)}(\cB^{\f54})}\bigr).
\end{split}
\eeno
Whereas thanks to \eqref{S3eq11},
we get, by applying Lemma \ref{lem:Bern}, that
\beno
\begin{split}
\int_0^t\bigl|&\bigl(\D_k^{\rm h}(R^\h(v, \p_yw^1))_\Th\ |\ \D_k^{\rm
h}\bigl(\p_tw^1+\f12w^1\bigr)_\Th\bigr)_{L^2}\bigr|\,dt'\\
\lesssim &2^{\f{k}2}\sum_{k'\geq k-3}\int_0^t\|{\D}_{k'}^\h v_\Th(t')\|_{L^\infty_\h(L^2_{\rm v})}\|\wt{\D}_{k'}^\h\p_y w^1_\Th(t')\|_{L^2}
\|\D_k^\h (w^1,\p_tw^1)_\Th(t')\|_{L^2}\,dt'\\
\lesssim & 2^{\f{k}2}\sum_{k'\geq k-3}2^{-\f{k'}4}\|u_\Th\|^{\f12}_{L^\infty_t(\cB^{2})}
\|\wt{\D}_{k'}^\h\p_y w^1_\Th\|_{L^2_t(L^2)} \Bigl(\int_0^t\|\D_k^\h (w^1,\p_tw^1)_\Th(t')\|_{L^2}^2\dot\tht(t')\,dt'\Bigr)^{\frac12}\\
\lesssim &d_k^22^{-k}\|u_\Th\|^{\f12}_{L^\infty_t(\cB^{2})}
\|\p_yw^1_\Th\|_{\wt{L}^2_t(\cB^{\f12})}\bigl(\|(w^1_\Th,\p_tw^1_\Th)\|_{\wt{L}^2_{t,\dot\tht(t)}(\cB^{\f34})}
+\lam\|w^1_\Th\|_{\wt{L}^2_{t,\dot\tht^3(t)}(\cB^{\f54})}\bigr).
\end{split}
\eeno
The same estimate holds for $
\int_0^t\bigl|\bigl(\D_k^{\rm h}(T^\h_{\p_yw^1}v)_\Th\ |\ \D_k^{\rm
h}\bigl(\p_tw^1+\f12w^1\bigr)_\Th\bigr)_{L^2}\bigr|\,dt'.$

As a consequence, we arrive at
\beq \label{S5eq12}
\begin{split}
\int_0^t&\bigl|\bigl(\D_k^{\rm h}({v}\p_yw^1)_\Th\ |\ \D_k^{\rm
h}\bigl(\p_tw^1+\f12w^1\bigr)_\Th\bigr)_{L^2}\bigr|\,dt'\\
&\lesssim d_k^22^{-k}\|u_\Th\|^{\f12}_{L^\infty_t(\cB^{2})}
\|\p_yw^1_\Th\|_{\wt{L}^2_t(\cB^{\f12})}\bigl(\|(w^1_\Th,\p_tw^1_\Th)\|_{\wt{L}^2_{t,\dot\tht(t)}(\cB^{\f34})}
+\lam\|w^1_\Th\|_{\wt{L}^2_{t,\dot\tht^3(t)}(\cB^{\f54})}\bigr).
\end{split}
\eeq

\noindent$\bullet$\underline{
The estimate of $\int_0^t\big|\big(\Delta_k^\h(w^2\p_yu)_\Th | \Delta_k^\h(\p_tw^1+\f12 w^1)_\Th\big)_{L^2}\big|\,dt'$.}

By applying Bony's decomposition \eqref{Bony}  to $w^2\p_yu$, we write
\beno
w^2\p_y u=T^\h_{w^2}\p_yu+T^\h_{\p_yu}w^2+R^h(w^2,\p_yu).
\eeno
In view of  \eqref{S3eq12}, we have
\beno
\begin{split}
\Bigl(\int_0^t&\|S_{k'-1}^\h w^2_\Th(t')\|_{L^\infty}^2\dot\tht^3(t')\,dt'\Bigr)^{\f12}
\lesssim d_{k'} 2^{\f{k'}4}\|w^1_\Th\|_{\wt{L}^2_{t,\dot{\tht}^3(t)}(\cB^{\f54})},
\end{split}
\eeno
so that we get, by applying H\"older's inequality, that
\beno
\begin{split}
\int_0^t&\bigl|\bigl(\D_k^{\rm h}(T^\h_{w^2}\p_y u)_\Th\ |\ \D_k^{\rm
h}(\p_tw^1,w^1)_\Th\bigr)_{L^2}\bigr|\,dt'\\
\lesssim & \sum_{|k'-k|\leq 4}2^{-\f{k'}2}\int_0^t\|S_{k'-1}^\h w^2_\Th(t')\|_{L^\infty}\|\p_yu_\Th(t')\|_{\cB^{\f12}}\|\D_k^\h (\p_tw^1,w^1)_\Th(t')\|_{L^2}\,dt'\\
\lesssim & \sum_{|k'-k|\leq 4}2^{-\f{k'}2}\Bigl(\int_0^t\|S_{k'-1}^\h w^2_\Th(t')\|_{L^\infty}^2\dot\tht^3(t')\,dt'\Bigr)^{\f12}\\
&\qquad\qquad\qquad\times\Bigl(\int_0^t\|\D_k^\h (\p_tw^1,w^1)_\Th(t')\|_{L^2}^2\dot\tht(t')\,dt'\Bigr)^{\f12}\\
\lesssim & d_k^22^{-k}\|w^1_\Th\|_{\wt{L}^2_{t,\dot\tht^3(t)}(\cB^{\f54})}\bigl(\|(w^1_\Phi,\p_tw^1_\Th)\|_{\wt{L}^2_{t,\dot\tht(t)}(\cB^{\f34})}+\lam \| w^1_\Th\|_{\wt{L}^2_{t,\dot{\tht}^3(t)}(\cB^{\f54})}\bigr).
\end{split}
\eeno
While thanks to \eqref{S3eq11}, we get, by applying Lemma \ref{lem:Bern}, that
\beno
\begin{split}
\int_0^t\bigl|&\bigl(\D_k^{\rm h}(R^\h(w^2, \p_y u))_\Th\ |\ \D_k^{\rm
h}(w^1,\p_tw^1)_\Th\bigr)_{L^2}\bigr|\,dt'\\
\lesssim &2^{\f{k}2}\sum_{k'\geq k-3}\int_0^t\|{\D}_{k'}^\h w^2_\Th(t')\|_{L^2_\h(L^\infty_{\rm v})}
\|\wt{\D}_{k'}^\h\p_y u_\Th(t')\|_{L^2}\|\D_k^\h w^1_\Th(t')\|_{L^2}\,dt'\\
\lesssim & 2^{\f{k}2}\sum_{k'\geq k-3}2^{\f{k'}2}\int_0^t\|{\D}_{k'}^\h w^1_\Th(t')\|_{L^2}\|\p_yu_\Th(t')\|_{\cB^{\f12}}\|\D_k^\h w^1_\Th(t')\|_{L^2}\,dt'\\
\lesssim &d_k^22^{-k}\|w^1_\Th\|_{\wt{L}^2_{t,\dot\tht^3(t)}(\cB^{\f54})}\bigl(\|(w^1_\Phi,\p_tw^1_\Th)\|_{\wt{L}^2_{t,\dot\tht(t)}(\cB^{\f34})}+\lam \| w^1_\Th\|_{\wt{L}^2_{t,\dot{\tht}^3(t)}(\cB^{\f54})}\bigr).
\end{split}
\eeno
The same estimate holds for $
\int_0^t\bigl|\bigl(\D_k^{\rm h}(T^\h_{\p_yu}w^2)_\Th\ |\ \D_k^{\rm
h}(w^1,\p_tw^1)_\Th\bigr)_{L^2}\bigr|\,dt'.$

Hence we obtain
\beq \label{S5eq13}
\begin{split}
\int_0^t\bigl|\bigl(&\D_k^{\rm h}(w^2\p_yu)_\Th\ |\ \D_k^{\rm
h}(\p_tw^1+\f12w^1)_\Th\bigr)_{L^2}\bigr|\,dt'\\
\lesssim & d_k^22^{-k}\|w^1_\Th\|_{\wt{L}^2_{t,\dot\tht^3(t)}(\cB^{\f54})}\bigl(\|(w^1_\Phi,\p_tw^1_\Th)\|_{\wt{L}^2_{t,\dot\tht(t)}(\cB^{\f34})}+\lam \| w^1_\Th\|_{\wt{L}^2_{t,\dot{\tht}^3(t)}(\cB^{\f54})}\bigr).
\end{split}
\eeq

By summing up (\ref{S5eq9}-\ref{S5eq13}), we conclude the proof of \eqref{S5eq7}.
\end{proof}

\begin{proof}[Proof of \eqref{S5eq8}] We first observe from $\p_xu+\p_yv=0$ and Poincar\'e inequality that
\beq \label{S5eq14}
\begin{split}
\int_0^t\big|\big(\Delta_k^\h (\p_t^2v)_{\Th}|\Delta_k^\h (\p_tw^2,w^2)_\Th\big)_{L^2}\big|dt'\lesssim& d_k^22^{-k}\|(\p_t^2u)_\Th\|_{\wt{L}^2_{t}(\cB^{\f32})}\|(\p_yw^2,\p_tw^2)_\Th\|_{\wt{L}^2_{t}(\cB^{\f12})},\\
\int_0^t\big|\big(\Delta_k^\h (\p_tv)_{\Th}|\Delta_k^\h (w^2,\p_tw^2)_\Th\big)_{L^2}\big|dt'\lesssim&d_k^22^{-k}\|(\p_tu)_\Th\|_{\wt{L}^2_{t}(\cB^{\f32})}
\|(\p_yw^2,\p_tw^2)_\Th\|_{\wt{L}^2_{t}(\cB^{\f12})},\\
\int_0^t\big|\big(\Delta_k^\h (\p_y^2v)_{\Th}|\Delta_k^\h (w^2,\p_tw^2)_\Th\big)_{L^2}\big|dt'\lesssim& d_k^22^{-k}\|\p_yu_\Th\|_{\wt{L}^2_{t}(\cB^{\f32})}\|(\p_yw^2,\p_tw^2)_\Th\|_{\wt{L}^2_{t}(\cB^{\f12})},\\
\int_0^t\big|\big(\Delta_k^\h (\p_x^2v)_{\Th}|\Delta_k^\h (w^2,\p_tw^2)_\Th\big)_{L^2}\big|dt'\lesssim& d_k^22^{-k}\|u_\Th\|_{\wt{L}^2_{t}(\cB^{\f72})}\|(\p_yw^2,\p_tw^2)_\Th\|_{\wt{L}^2_{t}(\cB^{\f12})}.
\end{split}
\eeq

\noindent$\bullet$\underline{
The estimate of $\int_0^t\big|\big(\Delta_k^\h(u^\e\p_xv^\e)_\Th | \Delta_k^\h (\p_tw^2+\f12w^2)_\Th\big)_{L^2}\big|\,dt'$.}

We write $$ u^\e\p_xv^\e=u^\e\p_xw^2+u^\e\p_xv.$$
It follows from \eqref{S3eq18}, \eqref{S4eq-1} and Lemma \ref{lem3.1} that
\beq\label{S5eq17}
\begin{split}
\int_0^t\big|\big(&\Delta_k^\h(u^\e\p_x w^2)_{\Th}|\Delta_k^\h (w^2,\p_tw^2)_\Th\big)_{L^2}\big|\,dt'\\
&\lesssim d_k^22^{-k} \|w^2_\Th\|_{\wt{L}^2_{t,\dot{\tht}^3(t)}(\cB^{\f54})}\bigl(\|(w^2_\Phi,\p_tw^2_\Th)\|_{\wt{L}^2_{t,\dot\tht(t)}(\cB^{\f34})}+\lam \| w^2_\Th\|_{\wt{L}^2_{t,\dot{\tht}^3(t)}(\cB^{\f54})}\bigr).
\end{split}
\eeq
By applying Bony's decomposition in the horizontal variable to $u^\e\p_x v$  gives
\beno u^\e\p_x v=T^\h_{u^\e}\p_x v+T^\h_{\p_x v}{u^\e}+R^\h({u^\e},\p_xv).
\eeno
Due to $$ \|S_{k'-1}^\h u^\e_\Th(t')\|_{L^\infty}\lesssim \|u^\e_\Th(t')\|_{\cB^{\f12}}^{\f12}\|\p_y u^\e_\Th(t')\|_{\cB^{\f12}}^{\f12},$$
and \eqref{S3eq11}, we have
\beno
\begin{split}
\int_0^t&\bigl|\bigl(\D_k^{\rm h}(T^\h_{u^\e}\p_x v)_\Th\ |\ \D_k^{\rm
h}(w,\p_tw^2)_\Th\bigr)_{L^2}\bigr|\,dt'\\
\lesssim & \sum_{|k'-k|\leq 4}\int_0^t\|S_{k'-1}^\h u^\e_\Th(t')\|_{L^\infty}
\|\D_{k'}^\h\p_x v_\Th(t')\|_{L^2}\|\D_k^\h (w^2,\p_tw^2)_\Th(t')\|_{L^2}\,dt'\\
\lesssim & \sum_{|k'-k|\leq 4}2^{2{k'}}\|u^\e_\Th\|_{L^\infty_t(\cB^{\f12})}^{\f12}\|\D_k^\h u_\Th\|_{L^2_t(L^2)}
\Bigl(\int_0^t\dot\tht(t')\|\D_k^\h (w^2,\p_tw^2)_\Th(t')\|_{L^2}^2\,dt'\Bigr)^{\f12}\\
\lesssim & d_k^22^{-k}\|u^\e_\Th\|_{L^\infty_t(\cB^{\f12})}^{\f12}\| u_\Th\|_{\wt{L}_t^2(\cB^{\f94})}\bigl(\|(w^2_\Phi,\p_tw^2_\Th)\|_{\wt{L}^2_{t,\dot\tht(t)}(\cB^{\f34})}+\lam \| w^2_\Th\|_{\wt{L}^2_{t,\dot{\tht}^3(t)}(\cB^{\f54})}\bigr).
\end{split}
\eeno
While again thanks to \eqref{S3eq11}, we find
\beno
\|S_{k'-1}^\h\p_xv_\Th\|_{L^2_t(L^\infty)}\lesssim d_{k'}2^{\f{k'}4}\|u_\Th\|_{\wt{L}^2_t(\cB^{\f94})},
\eeno
from which, we infer
\beno
\begin{split}
\int_0^t&\bigl|\bigl(\D_k^{\rm h}(T^\h_{\p_xv} u^\e)_\Th\ |\ \D_k^{\rm
h}(w^2,\p_tw^2)_\Th\bigr)_{L^2}\bigr|\,dt'\\
\lesssim & \sum_{|k'-k|\leq 4}\int_0^t\|S_{k'-1}^\h \p_xv_\Th(t')\|_{L^\infty}
\|\D_{k'}^\h u^\e_\Th(t')\|_{L^2}\|\D_k^\h (w^2,\p_tw^2)_\Th(t')\|_{L^2}\,dt'\\
\lesssim & \sum_{|k'-k|\leq 4}d_{k'}2^{-\f{k'}4}\|u_\Th\|_{\wt{L}^2_t(\cB^{\f94})}\|u^\e\|_{L^\infty_t(\cB^{\f12})}^{\f12}\Bigl(\int_0^t
\dot\tht(t')\|\D_k^\h (w^2,\p_tw^2)_\Th(t')\|_{L^2}^2\,dt'\Bigr)^{\f12}\\
\lesssim & d_k^22^{-k}\|u^\e_\Th\|_{L^\infty_t(\cB^{\f12})}^{\f12}\| u_\Th\|_{\wt{L}_t^2(\cB^{\f94})}\bigl(\|(w^2_\Phi,\p_tw^2_\Th)\|_{\wt{L}^2_{t,\dot\tht(t)}(\cB^{\f34})}+\lam \| w^2_\Th\|_{\wt{L}^2_{t,\dot{\tht}^3(t)}(\cB^{\f54})}\bigr).
\end{split}
\eeno
Along the same line, we obtain
\beno
\begin{split}
\int_0^t\bigl|&\bigl(\D_k^{\rm h}(R^\h(u^\e, \p_xv))_\Th\ |\ \D_k^{\rm
h}(w^2,\p_tw^2)_\Th\bigr)_{L^2}\bigr|\,dt'\\
\lesssim &2^{\f{k}2}\sum_{k'\geq k-3}\int_0^t\|{\D}_{k'}^\h u^\e_\Th(t')\|_{L^2_\h(L^\infty_{\rm v})}
\|\wt{\D}_{k'}^\h\p_xv_\Th(t')\|_{L^2}\|\D_k^\h w^1_\Th(t')\|_{L^2}\,dt'\\
\lesssim & 2^{\f{k}2}\sum_{k'\geq k-3}2^{\f{3k'}2}\|u^\e_\Th\|_{L^\infty_t(\cB^{\f12})}^{\f12}\|\D_k^\h u_\Th\|_{L^2_t(L^2)}
\Bigl(\int_0^t\dot\tht(t')\|\D_k^\h w^2_\Th(t')\|_{L^2}^2\,dt'\Bigr)^{\f12}\\
\lesssim & d_k^22^{-k}\|u^\e_\Th\|_{L^\infty_t(\cB^{\f12})}^{\f12}\| u_\Th\|_{\wt{L}_t^2(\cB^{\f94})}\bigl(\|(w^2_\Phi,\p_tw^2_\Th)\|_{\wt{L}^2_{t,\dot\tht(t)}(\cB^{\f34})}+\lam \| w^2_\Th\|_{\wt{L}^2_{t,\dot{\tht}^3(t)}(\cB^{\f54})}\bigr).
\end{split}
\eeno
This gives rise to
\beq \label{S5eq18}
\begin{split}
\int_0^t\bigl|\bigl(&\D_k^{\rm h}(u^\e\p_xv)_\Th\ |\ \D_k^{\rm
h}(\p_tw^2+\f12w^2)_\Th\bigr)_{L^2}\bigr|\,dt'\\
&\lesssim d_k^22^{-k}\|u^\e_\Th\|_{L^\infty_t(\cB^{\f12})}^{\f12}\| u_\Th\|_{\wt{L}_t^2(\cB^{\f94})}\bigl(\|(w^2_\Phi,\p_tw^2_\Th)\|_{\wt{L}^2_{t,\dot\tht(t)}(\cB^{\f34})}+\lam \| w^2_\Th\|_{\wt{L}^2_{t,\dot{\tht}^3(t)}(\cB^{\f54})}\bigr).
\end{split}
\eeq

\noindent$\bullet$\underline{
The estimate of $\int_0^t\big|\big(\Delta_k^\h(v^\e\p_yv^\e)_\Th | \Delta_k^\h (\p_tw^2+\f12w^2)_\Th\big)_{L^2}\big|\,dt'$.}

We first note that
\beno
v^\e\p_y v^\e=v\p_yw^2+w^2\p_yw^2+v\p_yv+w^2\p_yv.\eeno
We first deduce Lemma \ref{lem3.3} that
\begin{align*}
\e\int_0^t\bigl|\bigl(&\D_k^{\rm h}(w^2\p_yw^2)_\Th\ |\ \D_k^{\rm
h}(w^2,\p_tw^2)_\Th\bigr)_{L^2}\bigr|\,dt'\\
&\lesssim d_k^22^{-k}
\bigl\|(w^1_\Th,\e w^2_\Th)\bigr\|_{\wt{L}^2_{t,\dot\tht^3(t)}(\cB^{\f54})}\bigl(\|(w^2_\Phi,\p_tw^2_\Th)\|_{\wt{L}^2_{t,\dot\tht(t)}(\cB^{\f34})}+\lam \| w^2_\Th\|_{\wt{L}^2_{t,\dot{\tht}^3(t)}(\cB^{\f54})}\bigr),
\end{align*}
and
\begin{align*}
\e\int_0^t\bigl|\bigl(&\D_k^{\rm h}({v}\p_y v)_\Th\ |\ \D_k^{\rm
h}w_\Th^2\bigr)_{L^2}\bigr|\,dt'\\
\lesssim & d_k^22^{-k}\bigl\|(u,\e v)_\Th\bigr\|_{\wt{L}^2_{t,\dot\tht^3(t)}(\cB^{\f54})}\bigl(\|(w^2_\Phi,\p_tw^2_\Th)\|_{\wt{L}^2_{t,\dot\tht(t)}(\cB^{\f34})}+\lam \| w^2_\Th\|_{\wt{L}^2_{t,\dot{\tht}^3(t)}(\cB^{\f54})}\bigr).
\end{align*}

It follows from \eqref{S5eq12} that
\begin{align*}
\int_0^t\bigl|\bigl(&\D_k^{\rm h}({v}\p_yw^2)_\Th\ |\ \D_k^{\rm
h}(w^2,\p_tw^2)_\Phi\bigr)_{L^2}\bigr|\,dt'\\
&\lesssim d_k^22^{-k}\|u_\Th\|^{\f12}_{L^\infty_t(\cB^{2})}
\|\p_yw^2_\Th\|_{\wt{L}^2_t(\cB^{\f12})}\bigl(\|(w^2_\Th,\p_tw^2_\Th)\|_{\wt{L}^2_{t,\dot\tht(t)}(\cB^{\f34})}
+\lam\|w^2_\Th\|_{\wt{L}^2_{t,\dot\tht^3(t)}(\cB^{\f54})}\bigr).
\end{align*}

And \eqref{S5eq11} ensures that
\begin{align*}
\int_0^t\bigl|\bigl(\D_k^{\rm h}&({w^2}\p_xu)_\Th\ |\ \D_k^{\rm
h}(w^2,\p_tw^2)_\Phi\bigr)_{L^2}\bigr|\,dt'\\
\lesssim &d_k^22^{-k}\bigl(\|u_\Th\|^{\f12}_{L^\infty_t(\cB^{2})}
\|\p_yw^2_\Th\|_{\wt{L}^2_t(\cB^{\f12})}+
\|w^2_\Th\|_{\wt{L}^2_{t,\dot\tht^3(t)}(\cB^{\f54})}\bigr)\\
&\qquad\qquad\qquad\times\bigl(\|(w^2_\Th,\p_tw^2_\Th)\|_{\wt{L}^2_{t,\dot\tht(t)}(\cB^{\f34})}
+\lam\|w^2_\Th\|_{\wt{L}^2_{t,\dot\tht^3(t)}(\cB^{\f54})}\bigr).
\end{align*}

As a result, it comes out
\beq \label{S5eq21}
\begin{split}
\e^2\int_0^t&\big|\big(\Delta_k^\h(v^\e\p_yv^\e)_\Th | \Delta_k^\h w^1_\Th\big)_{L^2}\big|\,dt'
\lesssim  \e d_k^22^{-k}\Bigl(\|(u,\e v, w^1,\e w^2)_\Th\bigr\|_{\wt{L}^2_{t,\dot\tht^3(t)}(\cB^{\f54})}\\
&+\e\|u_\Th\|^{\f12}_{L^\infty_t(\cB^{2})}
\|\p_yw^2_\Th\|_{\wt{L}^2_t(\cB^{\f12})}\Bigr)\bigl(\|(w^2_\Th,\p_tw^2_\Th)\|_{\wt{L}^2_{t,\dot\tht(t)}(\cB^{\f34})}
+\lam\|w^2_\Th\|_{\wt{L}^2_{t,\dot\tht^3(t)}(\cB^{\f54})}\bigr).
\end{split}
\eeq

Summing up (\ref{S5eq14}-\ref{S5eq21}) gives rise to \eqref{S5eq8}.
\end{proof}

\renewcommand{\theequation}{\thesection.\arabic{equation}}
\setcounter{equation}{0}
\appendix

\section{Littlewood-Paley theory and functional framework}\label{sect2}

In this section, we shall collect some basic facts on anisotropic  Littlewood-Paley theory,
which we have been used in this context.
Let us first recall from
\cite{BCD} that \beq
\begin{split}
&\Delta_k^{\rm h}a \eqdefa \cF^{-1}(\varphi(2^{-k}|\xi|)\widehat{a}),\qquad
S^{\rm h}_ka \eqdefa \cF^{-1}(\chi(2^{-k}|\xi|)\widehat{a}),
\end{split} \label{1.3a}\eeq where $\cF
a$ and $\widehat{a}$  denote the partial  Fourier transform of
the distribution $a$ with respect to $x$ variable,  that is, $
\widehat{a}(\xi,y)=\cF_{x\to\xi}(a)(\xi,y),$
  and $\chi(\tau),$ ~$\varphi(\tau)$ are
smooth functions such that
 \beno
&&\Supp \varphi \subset \Bigl\{\tau \in \R\,/\  \ \frac34 \leq
|\tau| \leq \frac83 \Bigr\}\andf \  \ \forall
 \tau>0\,,\ \sum_{j\in\Z}\varphi(2^{-j}\tau)=1,\\
&&\Supp \chi \subset \Bigl\{\tau \in \R\,/\  \ \ |\tau|  \leq
\frac43 \Bigr\}\quad \ \ \ \andf \  \ \, \chi(\tau)+ \sum_{j\geq
0}\varphi(2^{-j}\tau)=1.
 \eeno

\begin{Def}\label{def1.2}
{\sl  Let~$s\in \R$. For~$u\in {S}_h'(\cS),$ which
means that $u$ belongs to ~$S'(\cS)$ and
satisfies~$\lim_{k\to-\infty}\|S_k^{\rm h}u\|_{L^\infty}=0,$ we set
$$
\|u\|_{\cB^{s}}\eqdefa\big\|\big(2^{ks}\|\Delta_k^{\rm h}
u\|_{L^{2}}\big)_{k\in\Z}\bigr\|_{\ell ^{1}(\Z)}.
$$
\begin{itemize}

\item
For $s\leq \frac{1}{2}$, we define $ \cB^{s}(\cS)\eqdefa
\big\{u\in{S}_h'(\cS)\;\big|\; \|
u\|_{\cB^{s}}<\infty\big\}.$

\item
If $k$ is  a positive integer and if~$-\frac{1}{2}+k< s\leq
\frac{1}{2}+k$, then we define~$ \cB^{s}(\cS)$  as the subset
of distributions $u$ in~${S}_h'(\cS)$ such that
$\p_x^k u$ belongs to~$ \cB^{s-k}(\cS).$
\end{itemize}
}
\end{Def}

In  order to obtain a better description of the regularizing effect
of the diffusion equation, we need to use Chemin-Lerner
type spaces $\widetilde{L}^{\lambda}_T(\cB^{s}(\cS))$ from \cite{CL95}.
\begin{Def}\label{def2.2}
{\sl Let $p\in[1,\,+\infty]$ and $T\in]0,\,+\infty]$. We define
$\widetilde{L}^{p}_T(\cB^{s}(\cS))$ as the completion of
$C([0,T]; \,S(\cS))$ by the norm
$$
\|a\|_{\widetilde{L}^{p}_T(\cB^{s})} \eqdefa \sum_{k\in\Z}2^{ks}
\Big(\int_0^T\|\Delta_k^{\rm h}\,a(t) \|_{L^2}^{p}\,
dt\Big)^{\frac{1}{p}}
$$
with the usual change if $p=\infty.$ }
\end{Def}

In order to overcome the difficulty that one can not use Gronwall
type argument in the framework of Chemin-Lerner type spaces, we  need to use the time-weighted
Chemin-Lerner norm, which was introduced by us in
\cite{PZ1}.

\begin{Def}\label{def1.1} {\sl Let $f(t)\in L^1_{\mbox{loc}}(\R^+)$
be a nonnegative function. We define \beq \label{1.4}
\|a\|_{\wt{L}^p_{t,f}(\cB^{s})}\eqdefa
\sum_{k\in\Z}2^{ks}\Bigl(\int_0^t f(t')\|\D_k^{\rm
h}a(t')\|_{L^2}^p\,dt'\Bigr)^{\f1p}. \eeq}
\end{Def}

 \medbreak
For the convenience of the readers, we recall the following anisotropic
Bernstein type lemma from \cite{CZ1, Pa02}.

\begin{lem} \label{lem:Bern}
 {\sl Let $\bold{B}_{\rm h}$ be a ball
of~$\R_{\rm h}$, and~$\bold{C}_{\rm h}$  a ring of~$\R_{\rm
h}$; let~$1\leq p_2\leq p_1\leq \infty$ and ~$1\leq q\leq \infty.$
Then there holds:

\smallbreak\noindent If the support of~$\wh a$ is included
in~$2^k\bold{B}_{\rm h}$, then
\[
\|\partial_{x}^N a\|_{L^{p_1}_{\rm h}(L^{q}_{\rm v})} \lesssim
2^{k\left(N+\f 1 {p_2}-\f 1 {p_1}\right)}
\|a\|_{L^{p_2}_{\rm h}(L^{q}_{\rm v})}.
\]

\smallbreak\noindent If the support of~$\wh a$ is included
in~$2^k\bold{C}_{\rm h}$, then
\[
\|a\|_{L^{p_1}_{\rm h}(L^{q}_{\rm v})} \lesssim
2^{-kN} \|\partial_{x}^N a\|_{L^{p_1}_{\rm
h}(L^{q}_{\rm v})}.
\]
}
\end{lem}

In this context, we  constantly use  Bony's decomposition (see \cite{Bo}) in
the horizontal variable: \ben\label{Bony} fg=T^{\rm h}_fg+T^{\rm
h}_{g}f+R^{\rm h}(f,g), \een where \beno T^{\rm h}_fg\eqdefa \sum_kS^{\rm
h}_{k-1}f\Delta_k^{\rm h}g,\andf R^{\rm
h}(f,g)\eqdefa\sum_k{\Delta}_k^{\rm h}f\widetilde{\Delta}_{k}^{\rm h}g \eeno
with $\widetilde{\Delta}_k^{\rm h}g\eqdefa
\displaystyle\sum_{|k-k'|\le 1}\Delta_{k'}^{\rm h}g$.

\renewcommand{\theequation}{\thesection.\arabic{equation}}
\setcounter{equation}{0}
\section{The proof of \eqref{S7eq5}}\label{sectb}

The goal of this section is to present the proof of \eqref{S7eq5}.

\begin{proof}[Proof of \eqref{S7eq5}] For simplicity, we denote $\frak{W}(t,x)\eqdefa \om(t,x,1)-\om(t,x,0).$
In view of \eqref{S7eq2}, we write
\begin{align*}
\int_{\R} \ekt \D_k^{\h} \p_x p_\Phi& \cdot \ekt  \D_k^\h \bigl(\p_t\fw+\f12\fw\bigr)_\Phi\,dx
=\bigl(\f12-\fk\bigr)\int_{\R}  \bigl(\ekt\D_k^\h\fw_\Phi\bigr)^2\, dx\\
&+\f 12\f{d}{dt}\int_{\R}  \bigl(\ekt\D_k^\h\fw_\Phi\bigr)^2\, dx
+\lambda\dot\theta\int_{\R}  \bigl(\ekt\D_k^\h|D_x|^{\f 14}\fw_\Phi\bigr)^2\, dx\\
&-\f12\int_{\R} \Bigl(\ekt\int_0^1\D_k^{\h}\p_x(u^2)_\Phi\,dy\Bigr)\cdot\ekt\D_k^\h\bigl(\p_t\fw+\f12\fw\bigr)_\Phi
\,dx.
\end{align*}
Since we shall not handle the estimate of the term $\p_t \p_y \om_\Phi$ and $\p_yp=0,$ we get, by integrating the above equation over $[0,t]$ and using integration
by parts, that
\beq \label{S7eq6}
\begin{split}
{\rm I}_k(t)=&\bigl(\f12-\fk\bigr)\int_0^t\int_{\R}  \bigl(\ekt\D_k^\h\fw_\Phi\bigr)^2\, dx\,dt'+\f12\int_{\R}  \bigl(\ekt\D_k^\h\fw_\Phi\bigr)^2\, dx\\
&
-\f12\int_{\R}  \bigl( e^{\de|D_x|^{\f12}}\D_k^\h\fw_0
\bigr)^2\, dx+\lambda\int_0^t \dot\theta(t')\int_{\R}  \bigl(\ekt|D_x|^{\f 14}\D_k^\h\fw_\Phi\bigr)^2\, dx\,dt'\\
&	-\f12\int_{\R} \Bigl(\ekt\int_0^1\D_k^{\h}\p_x(u^2)_\Phi\,dy\Bigr)\cdot\ekt\D_k^\h\fw_\Phi\,dx\\
&+\f 12\int_{\R} \Bigl(\int_0^1\D_k^{\h} e^{\de|D_x|^{\f12}}\p_x(u^2_0)\,dy\Bigr)\cdot e^{\de|D_x|^{\f12}}\D_k^\h\fw_0\,dx\\
&-\f\lambda2\int_0^t \dot\theta(t')\int_{\R} \Bigl(\ektp\int_0^1\D_k^{\h}\p_x(u^2)_\Phi\,dy\Bigr)\ektp|D_x|^{\f 12}\D_k^\h\fw_\Phi\,dx\,dt'\\
&+\bigl(\fk-\f14\bigr)\int_0^t \int_{\R} \Bigl(\ektp\int_0^1\D_k^{\h}\p_x(u^2)_\Phi\,dy\Bigr)\ektp\D_k^\h\fw_\Phi\,dx\,dt'\\
&+\int_0^t \int_{\R} \Bigl(\ekt\int_0^1\D_k^{\h}\p_x(u \p_t u)_\Phi\,dy\Bigr)\cdot\ekt\D_k^\h\fw_\Phi\,dx\,dt'\eqdefa \sum_{i=1}^9{\rm I}_k^i.
\end{split}
\eeq

Observing that
$$ \int_0^1 \Delta_k^{\h} \om_\Phi(t,x,y) dy = 0, $$
 for any fixed $(t,x)\in \mathbb{R}^+ \times \mathbb{R}$, there exists $Y_k(t,x)$
so that $\Delta_k^{\h} \om_\Phi(t,x,Y_k(t,x))= 0.$ So that we have
\begin{align*}
\bigl(\Delta_k^{\h} \om_\Phi(t,x,y)\bigr)^2 =&\bigl(\Delta_k^{\h} \om_\Phi(t,x,y)\bigr)^2-\bigl(\Delta_k^{\h} \om_\Phi(t,x,Y_k(t,x))\bigr)^2\\
\leq &\int_0^1\bigl|\p_y\bigl(\Delta_k^{\h} \om_\Phi(t,x,y)\bigr)^2\bigr|\,dy\leq 2\| \Delta_k^{\h} \om_\Phi(t,x,\cdot) \|_{L^2_{\rm v}} \|\Delta_k^{\h} \pa_y\om_\Phi(t,x,\cdot)\|_{L^2_{\rm v}},
\end{align*}
from which, we infer
\beq \label{S7eq7}
 \|\Delta_k^{\h} \om_\Phi(t) \|^2_{L^2_\h(L^\infty_{\rm v})} \leq 2\| \Delta_k^{\h} \om_\Phi \|_{L^2} \| \Delta_k^{\h} \pa_y\om_\Phi \|_{L^2}.
 \eeq

Thanks to \eqref{S7eq7}, we deduce that
\begin{align*}
|{\rm I}_k^1|\leq \f12\|\ektp \D_k^\h \om_\Phi\|_{L^2_t(L^\infty_{\rm v}(L^2_\h))}^2
\lesssim & \| \ektp\Delta_k^{\h} \om_\Phi \|_{L^2_t(L^2)} \|\ektp\Delta_k^{\h} \pa_y\om_\Phi\|_{L^2_t(L^2)}\\
\lesssim &
d_k^22^{-2k\fs}\|\ektp \om_\Phi\|_{\wt{L}^2_t(\cB^{\fs})}\|\ektp \p_y\om_\Phi\|_{\wt{L}^2_t(\cB^{\fs})}.
\end{align*}
Along the same line, one has
\begin{align*}
|{\rm I}_k^2|\leq \f12\|\ektp \D_k^\h \om_\Phi\|_{L^\infty_t(L^\infty_{\rm v}(L^2_\h))}^2
\lesssim & \| \ektp\Delta_k^{\h} \om_\Phi \|_{L^\infty_t(L^2)} \|\ektp\Delta_k^{\h} \pa_y\om_\Phi\|_{L^\infty_t(L^2)}\\
\lesssim &
d_k^22^{-2k\fs}\|\ektp \om_\Phi\|_{\wt{L}^\infty_t(\cB^{\fs})}\|\ektp \p_y\om_\Phi\|_{\wt{L}^\infty_t(\cB^{\fs})},
\end{align*}
and
\begin{align*}
|{\rm I}_k^3|
\lesssim
d_k^22^{-2k\fs}\bigl\|e^{\de |D_x|^{\f12}} \om_0\bigr\|_{\cB^{\fs}}\bigl\|e^{\de |D_x|^{\f12}} \p_y\om_0\bigr\|_{\cB^{\fs}}.
\end{align*}

Furthermore, we get, by a similar derivation of \eqref{S7eq7}, that
\beno
|{\rm I}_k^4|
\lesssim\lam
d_k^22^{-2k\fs}\|\ektp \om_\Phi\|_{\wt{L}^2_{t,\dot{\tht}(t)}(\cB^{\fs+\f14})}\|\ektp \p_y\om_\Phi\|_{\wt{L}^2_{t,\dot{\tht}(t)}(\cB^{\fs+\f14})}.
\eeno

To handle the estimate of ${\rm I}_k^5$ to ${\rm I}_k^8,$ we use Bony's decomposition \eqref{Bony} for $u^2$ in the horizontal variable
to write
\beno
u^2=2T^\h_uu+R^\h(u,u). \eeno
We first get, by applying  \eqref{S3eq8} and \eqref{eq4.3}, that
\begin{align*}
\bigl|\int_{\R}& \Bigl(\ekt\int_0^1\D_k^\h\p_x(R^\h(u,u))_\Phi\,dy\Bigr)\cdot\ekt\D_k^\h \fw_\Phi\,dx\bigr|\\
\lesssim &2^{k}\sum_{k'\geq k-3}\|\wt{\D}_{k'}^\h u_\Phi(t)\|_{L^\infty} \| \ekt\Delta_{k'}^\h u_\Phi(t) \|_{L^2} \|\ekt\Delta_k^\h \om_\Phi(t)\|_{L^\infty_{\rm v}(L^2_\h)}\\
\lesssim &\de 2^{k}\sum_{k'\geq k-3} \|e^{\f{\fk}2t}\Delta_{k'}^\h u_\Phi(t) \|_{L^2} \|e^{\f{\fk}2t}\Delta_{k'}^\h \om_\Phi(t)\|_{L^2}^{\f12}\|e^{{\fk}t}\Delta_k^\h \p_y \om_\Phi(t)\|_{L^2}^{\f12}\\
\lesssim &\de d_k2^{k\left(\f{3}4-\fs\right)}\|e^{\f{\fk}2 t'}u_\Phi\|_{\wt{L}^\infty_t(\cB^{\fs+\f34})}\|e^{\f{\fk}2t'}\om_\Phi\|_{\wt{L}^\infty_t(\cB^{\fs+\f12})}^{\f12}
\|e^{{\fk}t}\p_y\om_\Phi\|_{\wt{L}^\infty_t(\cB^{\fs})}^{\f12}\Bigl(\sum_{k'\geq k-3}d_{k'}2^{-k'\left(\fs+\f34\right)}\Bigr)\\
\lesssim &\de d^2_k2^{-2k\fs}\|e^{\f{\fk}2t'}u_\Phi\|_{\wt{L}^\infty_t(\cB^{\fs+\f34})}\|e^{\f{\fk}2t'}\om_\Phi\|_{\wt{L}^\infty_t(\cB^{\fs+\f12})}^{\f12}
\|e^{{\fk}t}\p_y\om_\Phi\|_{\wt{L}^\infty_t(\cB^{\fs})}^{\f12}.
\end{align*}
The same estimate holds for
$$ \int_{\R}\Bigl(\ekt\int_0^1\D_k^\h\p_x(T^\h_uu)_\Phi\,dy\Bigr)\cdot\ekt\D_k^\h \fw_\Phi\,dx.$$
Hence we obtain
\beno
|{\rm I}_k^5|
\lesssim \de d^2_k2^{-2k\fs}\|e^{\f{\fk}2t'}u_\Phi\|_{\wt{L}^\infty_t(\cB^{\fs+\f34})}\|e^{\f{\fk}2t'}\om_\Phi\|_{\wt{L}^\infty_t(\cB^{\fs+\f12})}^{\f12}
\|e^{{\fk}t}\p_y\om_\Phi\|_{\wt{L}^\infty_t(\cB^{\fs})}^{\f12}.
\eeno

While by applying the law of product in Besov spaces and \eqref{S7eq7} that
\begin{align*}
|{\rm I}_k^6|
\lesssim &
\| e^{\de|D_x|^{\f12}}\p_x\D_k^\h(u^2_0)\|_{L^2}\| e^{\de|D_x|^{\f12}}\D_k^\h\om_0\|_{L^\infty_{\rm v}(L^2_\h)}\\
\lesssim & 2^k
\| e^{\de|D_x|^{\f12}}\D_k^\h(u_0^2)\|_{L^2}\| e^{\de|D_x|^{\f12}}\D_k^\h\om_0\|_{L^2}^{\f12}\| e^{\de|D_x|^{\f12}}\D_k^\h\p_y\om_0\|_{L^2}^{\f12}\\
\lesssim & d_k^22^{-2k\fs}\bigl\|e^{\de|D_x|^{\f12}}u_0\bigr\|_{\cB^{\f12}}^{\f12}\bigl\|e^{\de|D_x|^{\f12}}\om_0\bigr\|_{\cB^{\f12}}^{\f12}
\bigl\|e^{\de|D_x|^{\f12}}u_0\bigr\|_{\cB^{\fs+1}}
\bigl\|e^{\de|D_x|^{\f12}}\om_0\bigr\|_{\cB^{\fs}}^{\f12}\bigl\|e^{\de|D_x|^{\f12}}\p_y\om_0\bigr\|_{\cB^{\fs}}^{\f12}.
\end{align*}

Similarly by applying Lemma \ref{lem:Bern}, \eqref{S3eq8} and \eqref{eq4.3}, we find
\begin{align*}
\bigl|\int_0^t& \dot\theta(t')\int_{\R} \Bigl(\ekt\int_0^1\D_k^{\h}\p_x(R^\h(u,u))_\Phi\,dy\Bigr)\cdot\ekt|D_x|^{\f 12}\D_k^\h\fw_\Phi\,dx\,dt'\bigr|\\
\lesssim &2^{\f{3k}2}\sum_{k'\geq k-3}\int_0^t \dot\theta(t') \|\wt{\D}_k^\h u_\Phi(t')\|_{L^2_\h(L^\infty_{\rm v})} \| \ektp\Delta_k^{\h} u_\Phi(t') \|_{L^2} \|\ektp|D_x|^{\f12} \Delta_k^{\h}\om_\Phi\|_{L^\infty_{\rm v}(L^2_\h)}\,dt'\\
\lesssim &2^{2k}\sum_{k'\geq k-3}2^{-\f{k'}2}\int_0^t\dot\tht^3(t') \| \ektp\Delta_k^{\h} u_\Phi(t') \|_{L^2} \|\ektp\Delta_k^{\h} \om_\Phi\|_{L^\infty_{\rm v}(L^2_\h)}\,dt',
\end{align*}
from which and \eqref{S7eq7}, we deduce that
\begin{align*}
\bigl|\int_0^t& \dot\theta(t')\int_{\R} \Bigl(\ekt\int_0^1\D_k^{\h}\p_x(R^\h(u,u))_\Phi\,dy\Bigr)\cdot\ekt|D_x|^{\f 12}\D_k^\h\fw_\Phi\,dx\,dt'\bigr|\\
\lesssim & 2^{2k}\sum_{k'\geq k-3}2^{-\f{k'}2}\Bigl(\int_0^t\dot\tht^3(t') \| \ektp\Delta_{k'}^{\h} u_\Phi(t') \|_{L^2}^2\,dt'\Bigr)^{\f12}
\\
&\qquad\qquad\times\Bigl(\int_0^t\dot\tht^3(t') \| \ektp\Delta_k^{\h} \om_\Phi(t') \|_{L^2}^2\,dt'\Bigr)^{\f14}\Bigl(\int_0^t\dot\tht^3(t') \| \ektp\Delta_k^{\h} \p_y\om_\Phi(t') \|_{L^2}^2\,dt'\Bigr)^{\f14}\\
\lesssim & \de^{\f14} d_k2^{k\left(\f32-\fs\right)} \|\ektp u_\Phi\|_{\wt{L}^2_{t,\dot{\tht}^3(t)}(\cB^{\fs+1})}\|\ektp \om_\Phi\|_{\wt{L}^2_{t,\dot{\tht}^3(t)}(\cB^{\fs+\f34})}^{\f12}
\|\ektp \p_y\om_\Phi\|_{\wt{L}^2_{t,\dot{\tht}(t)}(\cB^{\fs+\f14})}^{\f12}\\
&\qquad\times\Bigl(\sum_{k'\geq k-3} d_{k'}2^{-k'\left(\fs+\f32\right)}\Bigr)\\
\lesssim & \de^{\f14}d_k^22^{-2k\fs}\|\ektp u_\Phi\|_{\wt{L}^2_{t,\dot{\tht}^3(t)}(\cB^{\fs+1})}\|\ektp \om_\Phi\|_{\wt{L}^2_{t,\dot{\tht}^3(t)}(\cB^{\fs+\f34})}^{\f12}
\|\ektp \p_y\om_\Phi\|_{\wt{L}^2_{t,\dot{\tht}(t)}(\cB^{\fs+\f14})}^{\f12}.
\end{align*}
The same estimate holds for $$\int_0^t \dot\theta(t')\int_{\R} \Bigl(\ekt\int_0^1\D_k^{\h}\p_x(T^\h_uu)_\Phi\,dy\Bigr)\cdot\ekt |D_x|^{\f 12} \D_k^\h\fw_\Phi\,dx\,dt'.$$
Therefore, we obtain
\beno
|{\rm I}_k^7|
\lesssim \lam \de^{\f14}
d_k^22^{-2k\fs}\|\ektp u_\Phi\|_{\wt{L}^2_{t,\dot{\tht}^3(t)}(\cB^{\fs+1})}\|\ektp \om_\Phi\|_{\wt{L}^2_{t,\dot{\tht}^3(t)}(\cB^{\fs+\f34})}^{\f12}
\|\ektp \p_y\om_\Phi\|_{\wt{L}^2_{t,\dot{\tht}(t)}(\cB^{\fs+\f14})}^{\f12}.
\eeno

Along the same line, we have
\beq \label{S7eq8}
\begin{split}
\bigl|\int_0^t& \int_{\R} \Bigl(\ektp\int_0^1\D_k^{\h}\p_x(R^\h(u,u))_\Phi\,dy\Bigr)\cdot\ekt\D_k^\h\fw_\Phi\,dx\,dt'\bigr|\\
\lesssim &2^{k}\sum_{k'\geq k-3}\int_0^t \|\wt{\D}_{k'}^\h u_\Phi(t')\|_{L^\infty} \| \ektp\Delta_{k'}^{\h} u_\Phi(t') \|_{L^2} \|\ektp\Delta_k^{\h} \om_\Phi\|_{L^\infty_{\rm v}(L^2_\h)}\,dt'\\
\lesssim & 2^{k}\sum_{k'\geq k-3}\Bigl(\int_0^t\dot\tht^3(t') \| \ektp\Delta_{k'}^{\h} u_\Phi(t') \|_{L^2}^2\,dt'\Bigr)^{\f12}
\\
&\qquad\qquad\times\Bigl(\int_0^t\dot\tht(t') \| \ektp\Delta_{k'}^{\h} \om_\Phi(t') \|_{L^2}^2\,dt'\Bigr)^{\f14}\Bigl(\int_0^t\dot\tht(t') \| \ektp\Delta_k^{\h} \p_y\om_\Phi(t') \|_{L^2}^2\,dt'\Bigr)^{\f14}\\
\lesssim & d_k^22^{-2k\fs}\|\ektp u_\Phi\|_{\wt{L}^2_{t,\dot{\tht}^3(t)}(\cB^{\fs+\f34})}\|\ektp \om_\Phi\|_{\wt{L}^2_{t,\dot{\tht}(t)}(\cB^{\fs+\f14})}^{\f12}
\|\ektp \p_y\om_\Phi\|_{\wt{L}^2_{t,\dot{\tht}(t)}(\cB^{\fs+\f14})}^{\f12}.
\end{split} \eeq
The same estimate holds for
$$
\int_0^t\int_{\R} \Bigl(\ektp\int_0^1\D_k^{\h}\p_x(T^\h_uu)_\Phi\,dy\Bigr)\cdot\ekt\D_k^\h\fw_\Phi\,dx\,dt'.$$
So that there holds
\beno
|{\rm I}_k^8|
\lesssim
d_k^22^{-2k\fs}\|\ektp u_\Phi\|_{\wt{L}^2_{t,\dot{\tht}^3(t)}(\cB^{\fs+\f34})}\|\ektp \om_\Phi\|_{\wt{L}^2_{t,\dot{\tht}(t)}(\cB^{\fs+\f14})}^{\f12}
\|\ektp \p_y\om_\Phi\|_{\wt{L}^2_{t,\dot{\tht}(t)}(\cB^{\fs+\f14})}^{\f12}.
\eeno

Finally to deal with  the estimate of ${\rm I}_k^9,$ we use Bony's decomposition \eqref{Bony} for $uu_t$ in the horizontal variable
to write
\beno
u u_t=T^\h_uu_t+T^\h_{u_t}u+R^\h(u,u_t). \eeno
We first observe from \eqref{S3eq8} that
\begin{align*}
\bigl|\int_0^t& \int_{\R} \Bigl(\ektp\int_0^1\D_k^{\h}\p_x(T^\h_uu_t)_\Phi\,dy\Bigr)\cdot\ekt\D_k^\h\fw_\Phi\,dx\,dt'\bigr|\\
\lesssim &2^{k}\sum_{|k'-k|\leq 4}\int_0^t \|S_{k'-1}^\h u_\Phi(t')\|_{L^\infty} \| \ektp\Delta_k^{\h} (\p_tu)_\Phi(t') \|_{L^2} \|\ektp\Delta_k^{\h} \om_\Phi(t')\|_{L^\infty_{\rm v}(L^2_\h)}\,dt'\\
\lesssim & 2^{k}\sum_{|k'-k|\leq 4}\Bigl(\int_0^t\dot\tht^2(t') \| \ektp\Delta_k^{\h} (\p_tu)_\Phi(t') \|_{L^2}^2\,dt'\Bigr)^{\f12}
\\
&\qquad\qquad\times\Bigl(\int_0^t\dot\tht^3(t') \| \ektp\Delta_k^{\h} \om_\Phi(t') \|_{L^2}^2\,dt'\Bigr)^{\f14}\Bigl(\int_0^t\dot\tht(t') \| \ektp\Delta_k^{\h} \p_y\om_\Phi(t') \|_{L^2}^2\,dt'\Bigr)^{\f14}\\
\lesssim & d_k^22^{-2k\fs}\|\ektp (\p_tu)_\Phi\|_{\wt{L}^2_{t,\dot{\tht}^2(t)}(\cB^{\fs+\f12})}\|\ektp \om_\Phi\|_{\wt{L}^2_{t,\dot{\tht}^3(t)}(\cB^{\fs+\f34})}^{\f12}
\|\ektp \p_y\om_\Phi\|_{\wt{L}^2_{t,\dot{\tht}(t)}(\cB^{\fs+\f14})}^{\f12}.
\end{align*}
Notice from \eqref{S4eq10} and \eqref{S4eq-1} that
\begin{align*}
\bigl|\int_0^t& \int_{\R} \Bigl(\ektp\int_0^1\D_k^{\h}\p_x(R^\h(u,u_t))_\Phi\,dy\Bigr)\cdot\ekt\D_k^\h\fw_\Phi\,dx\,dt'\bigr|\\
\lesssim &2^{k}\sum_{k'\geq k-3}\int_0^t \|\wt{\D}_{k'}^\h (\p_tu)_\Phi(t')\|_{L^\infty_\h(L^2_{\rm v})} \| \ektp\Delta_{k'}^{\h} u_\Phi(t') \|_{L^2} \|\ektp\Delta_k^{\h} \om_\Phi(t')\|_{L^\infty_{\rm v}(L^2_\h)}\,dt'\\
\lesssim &2^{k}\sum_{k'\geq k-3}\int_0^t \dot{\tht}^2(t') \| \ektp\Delta_{k'}^{\h} u_\Phi(t') \|_{L^2}
\|\ektp\Delta_k^{\h} \om_\Phi(t')\|_{L^\infty_{\rm v}(L^2_\h)}\,dt',
\end{align*}
from which, we deduce by a similar derivation of \eqref{S7eq8} that
\begin{align*}
\bigl|\int_0^t& \int_{\R} \Bigl(\ektp\int_0^1\D_k^{\h}\p_x(R^\h(u,u_t))_\Phi\,dy\Bigr)\cdot\ekt\D_k^\h\fw_\Phi\,dx\,dt'\bigr|\\
\lesssim & d_k^22^{-2k\fs}\|\ektp u_\Phi\|_{\wt{L}^2_{t,\dot{\tht}^3(t)}(\cB^{\fs+\f34})}\|\ektp \om_\Phi\|_{\wt{L}^2_{t,\dot{\tht}(t)}(\cB^{\fs+\f14})}^{\f12}
\|\ektp \p_y\om_\Phi\|_{\wt{L}^2_{t,\dot{\tht}(t)}(\cB^{\fs+\f14})}^{\f12}.
\end{align*}
The same estimate holds for
$$
\int_0^t \int_{\R} \Bigl(\ektp\int_0^1\D_k^{\h}\p_x(T^\h_{u_t}u)_\Phi\,dy\Bigr)\cdot\ekt\D_k^\h\fw_\Phi\,dx\,dt'.
$$
As a result, it comes out
\begin{align*}
|{\rm I}_k^9|
\lesssim
d_k^22^{-2k\fs}\Bigl(&\|\ektp (\p_tu)_\Phi\|_{\wt{L}^2_{t,\dot{\tht}(t)}(\cB^{\fs+\f12})}\|\ektp \om_\Phi\|_{\wt{L}^2_{t,\dot{\tht}^3(t)}(\cB^{\fs+\f34})}^{\f12}\\
&+\|\ektp u_\Phi\|_{\wt{L}^2_{t,\dot{\tht}^3(t)}(\cB^{\fs+\f34})}\|\ektp \om_\Phi\|_{\wt{L}^2_{t,\dot{\tht}(t)}(\cB^{\fs+\f14})}^{\f12}\Bigr)
\|\ektp \p_y\om_\Phi\|_{\wt{L}^2_{t,\dot{\tht}(t)}(\cB^{\fs+\f14})}^{\f12}.
\end{align*}

By inserting the above estimates into \eqref{S7eq6}, we arrive at \eqref{S7eq5}.
\end{proof}

\section*{Acknowledgments}
 Both authors are supported by K.C.Wong Education Foundation.
M. Paicu was partially supported by the Agence Nationale de
la Recherche, Project IFSMACS, grant ANR-15-CE40-0010.  Ping Zhang is partially supported by  NSF of China under Grants   11731007, 12031006 and 11688101.

\end{document}